\documentclass[12pt,draft]{article}

\usepackage{amsmath}
\usepackage{amsthm}
\usepackage{amssymb}
\usepackage{amscd}
\usepackage{a4wide}
\usepackage{epic}
\usepackage{eepic}

\theoremstyle{definition}
\newtheorem{theorem}{Theorem}[section]
\newtheorem{prop}[theorem]{Proposition}
\newtheorem{lemma}[theorem]{Lemma}
\newtheorem{corollary}[theorem]{Corollary}
\newtheorem{definition}[theorem]{Definition}
\newtheorem{example}[theorem]{Example}
\newtheorem{remark}[theorem]{Remark}
\newtheorem{conjecture}[theorem]{Conjecture}

\newtheorem{claim}[theorem]{Claim}

\numberwithin{equation}{section}
\newenvironment{demo}[1]{%
  \trivlist
  \item[\hskip\labelsep
        {\it #1.}]
}{%
\hfill\qedsymbol
  \endtrivlist
}

\newcommand\Nat{\mathbb{N}}
\newcommand\Int{\mathbb{Z}}
\newcommand\Rat{\mathbb{Q}}
\newcommand\Comp{\mathbb{C}}

\newcommand\Strict{\mathcal{S}}
\newcommand\FF{\mathcal{F}}
\newcommand\GG{\mathcal{G}}
\newcommand\PP{\mathbb{P}}
\newcommand\Pf{\operatorname{Pf}}

\newcommand\Symp{\mathbf{Sp}}

\newcommand\sgn{\operatorname{sgn}}
\newcommand\ord{\operatorname{ord}}
\newcommand\inv{\operatorname{inv}}

\newcommand\vectx{\boldsymbol{x}}
\newcommand\vecty{\boldsymbol{y}}
\newcommand\vecta{\boldsymbol{a}}
\newcommand\vectb{\boldsymbol{b}}

\newcommand\trans{{}^t\!}
\newcommand\vectone{\boldsymbol{1}}
\newcommand\vectzero{\boldsymbol{0}}

\newcommand\ep{\varepsilon}
\renewcommand\tilde{\widetilde}
\renewcommand\hat{\widehat}

\allowdisplaybreaks

\title{
A Generalization of Schur's $P$- and $Q$-Functions
}

\author{
Soichi Okada%
\footnote{
Graduate School of Mathematics, Nagoya University, 
Furo-cho, Chikusa-ku, Nagoya 464-8602, Japan, 
{\tt okada@math.nagoya-u.ac.jp}
}
\footnote{
This work was partially supported by 
JSPS Grants-in-Aid for Scientific Research No.~24340003, No.~15K13425 and No.~18K03208.
}
}

\date{
\it Dedicated to Christian Krattenthaler on the occasion of his 60th birthday
}

\begin{document}

\maketitle

\begin{abstract}
We introduce and study a generalization of Schur's $P$-/$Q$-functions 
associated to a polynomial sequence, 
which can be viewed as ``Macdonald's ninth variation'' for $P$-/$Q$-functions.
This variation includes as special cases Schur's $P$-/$Q$-functions, 
Ivanov's factorial $P$-/$Q$-functions and the $t=-1$ specialization of Hall--Littlewood functions 
associated to the classical root systems.
We establish several identities and properties such as generalizations 
of Schur's original definition of Schur's $Q$-functions, Cauchy-type identity, 
J\'ozefiak--Pragacz--Nimmo formula for skew $Q$-functions, 
and Pieri-type rule for multiplication.
\end{abstract}

\setcounter{tocdepth}{1}
\tableofcontents

\section{%
Introduction
}

Schur ($S$-)functions and Schur $P$-/$Q$-functions are two important families 
of symmetric functions, and they appear in several parallel situations.
For example, 
in the representation theory of the symmetric groups, 
Schur functions describe the characters of irreducible linear representations,  
while Schur $Q$-functions describe the characters of irreducible projective representations 
(see \cite{Schur}).
In the cohomology theory, 
Schur functions represent the Schubert classes of Grassmannians, 
while Schur $Q$-functions represent the Schubert classes of Lagrangian Grassmannians 
(see \cite{Pragacz}).
Also some identities for Schur functions have their counterparts for Schur $P$-/$Q$-functions.

There are several generalizations, variations or deformations of Schur functions, 
such as Hall--Littlewood functions, Macdonald functions and factorial Schur functions.
The generalization relevant to this paper is 
Macdonald's ninth variation (\cite{Macdonald92}, see also \cite{NNSY}) 
associated to a polynomial sequence, which is defined as follows.

Let $\FF = \{ f_d \}_{d=0}^\infty$ be a sequence of polynomials $f_d(u) \in K[u]$, 
where $K$ is a ground field of characteristic $0$, such that $\deg f_d = d$ for $d \ge 0$.
Given a partition $\lambda$ of length $l \le n$,
we define the \emph{generalized Schur function} $s^\FF_\lambda(x_1, \dots, x_n)$ 
as the ratio of two alternants:
\begin{equation}
\label{eq:genS}
s^\FF_\lambda(x_1, \dots, x_n)
 =
\frac{ \det \big( f_{\lambda_j+n-j}(x_i) \big)_{1 \le i, j \le n} }
     { \det \big( f_{n-j}(x_i) \big)_{1 \le i, j \le n} },
\end{equation}
where $\lambda_{l+1} = \dots = \lambda_n = 0$.
The original Schur functions $s_\lambda(\vectx)$ are recovered by setting 
$f_d(u) = u^d$ for $d \ge 0$.
And the factorial Schur functions $s_\lambda(\vectx|\vecta)$ with factorial parameters 
$\vecta = (a_0, a_1, \dots)$ are obtained by taking $f_d(u) = (u|\vecta)^d 
= \prod_{i=0}^{d-1} (u - a_i)$.
Moreover, classical group characters are special cases of generalized Schur functions.
For example, if the polynomial sequence $\FF = \{ f_d \}_{d=0}^\infty$ is defined by 
$$
f_d( x + x^{-1} ) = \frac{ x^{d+1} - x^{-d-1} }{ x - x^{-1} }
\quad(d \ge 0),
$$
then it is not difficult to see that the generalized Schur function 
$s^\FF_\lambda(x_1+x_1^{-1}, \dots, x_n+x_n^{-1})$ 
equals to the irreducible character of the symplectic group $\Symp_{2n}(\Comp)$ 
with highest weight $\lambda$.

Generalized Schur functions share many of the same properties as the original Schur functions.
For example, they satisfy the modified Jacobi--Trudi identity and the Giambelli identity:
\begin{align*}
s^\FF_\lambda(x_1, \dots, x_n)
 &=
\det \big( s^\FF_{(\lambda_i-i+j)} (x_j, \dots, x_n) \big)_{1 \le i, j \le l}
\\
 &=
\det \big( s^\FF_{(\alpha_i|\beta_j)} (x_1, \dots, x_n) \big)_{1 \le i, j \le r},
\end{align*}
where $\lambda$ is a partition of length $l \le n$ 
and $\lambda = (\alpha_1, \dots, \alpha_r|\beta_1, \dots, \beta_r)$ in the Frobenius notation.

The aim of this paper is to introduce and study the ``ninth variation'' of Schur $P$-/$Q$-functions, 
which we call generalized $P$-functions associated to polynomial sequences. 
We define generalized $P$-functions in terms of Nimmo-type formula 
and derive Pfaffian identities and basic properties 
by following a linear algebraic approach similar to \cite{Okada}.

We use the following terminologies on polynomial sequences.

\begin{definition}
\label{def:F}
Let $\FF = \{ f_d \}_{d=0}^\infty$ be a sequence of polynomials $f_d(u) \in K[u]$.
We say that $\FF$ is \emph{admissible} if it satisfies the conditions
\begin{equation}
\label{eq:cond}
f_0(t) = 1,
\quad
\deg f_d = d \quad(d \ge 1).
\end{equation}
And an admissible sequence $\FF$ is called \emph{constant-term free} 
if $f_d(0) = 0$ for any $d \ge 1$.
\end{definition}

In this article, 
a \emph{partition of length $l$} is a weakly decreasing sequence 
$\lambda = (\lambda_1, \dots, \lambda_l)$ of \emph{positive} integers.
We write $l = l(\lambda)$ and $|\lambda| = \sum_{i=1}^l \lambda_i$.
A partition $\lambda$ of length $l$ is called \emph{strict} 
if $\lambda_1 > \dots > \lambda_l$.
The empty sequence $\emptyset$ is the unique strict partition of length $0$.

For a sequence $\vectx = (x_1, \dots, x_n)$ of $n$ indeterminates, 
we put
\begin{equation}
\label{eq:A}
A(\vectx) = \left( \frac{ x_j - x_i }{ x_j + x_i } \right)_{1 \le i, j \le n},
\quad
\Delta(\vectx) = \prod_{1 \le i < j \le n} \frac{ x_j - x_i }{ x_j + x_i }.
\end{equation}
Now we give a definition of generalized Schur $P$-functions associated to polynomial sequences 
in terms of Nimmo-type formula (see \cite[(A13)]{Nimmo}).

\begin{definition}
\label{def:genP}
For an admissible sequence $\FF = \{ f_d \}_{d=0}^\infty$ of polynomials 
and a sequence $\alpha = (\alpha_1, \dots, \alpha_r)$ of nonnegative integers, 
let $V^{\FF}_\alpha(\vectx)$ be the $n \times r$ matrix given by
$$
V^{\FF}_\alpha(\vectx)
 = 
\big( f_{\alpha_j}(x_i) \big)_{1 \le i \le n, 1 \le j \le r}.
$$
Given a strict partition $\lambda$ of length $l$,
we define the corresponding \emph{generalized $P$-function} $P^{\FF}_\lambda(\vectx)$ 
associated to $\FF$ by putting
\begin{equation}
\label{eq:Nimmo}
P^\FF_\lambda(\vectx)
 =
\begin{cases}
 \dfrac{ 1 }{ \Delta(\vectx) }
 \Pf \begin{pmatrix}
  A(\vectx) & V^\FF_\lambda(\vectx) \\
  - \trans V^\FF_\lambda(\vectx) & O
 \end{pmatrix}
 &\text{if $n+l$ is even,}
\\
 \dfrac{ 1 }{ \Delta(\vectx) }
 \Pf \begin{pmatrix}
  A(\vectx) & V^\FF_{\lambda^0}(\vectx) \\
  - \trans V^\FF_{\lambda^0}(\vectx) & O
 \end{pmatrix}
 &\text{if $n+l$ is odd,} \\
\end{cases}
\end{equation}
where $\lambda^0 = (\lambda_1, \dots, \lambda_l, 0)$.
We simply write $V_\alpha(\vectx)$ and $P_\lambda(\vectx)$ for 
$V^\FF_\alpha(\vectx)$ and $P^\FF_\lambda(\vectx)$ 
if there is no confusion, e.g., in the proofs.
\end{definition}

Note that (see Proposition~\ref{prop:Schur-Pf})
$$
\Delta(\vectx)
 =
\begin{cases}
 \Pf A(\vectx)
  &\text{if $n$ is even,} \\
 \Pf \begin{pmatrix}
  A(\vectx) & \vectone \\
  -\trans \vectone & 0
 \end{pmatrix}
  &\text{if $n$ is odd,}
\end{cases}
$$
where $\vectone$ is the all-one column vector of appropriate size.
Hence our definition (\ref{eq:Nimmo}) can be regarded as a counterpart of 
the definition (\ref{eq:genS}) of generalized Schur functions.

\begin{example}
\begin{enumerate}
\item[(1)]
It follows from Nimmo's formula \cite[(A13)]{Nimmo} that 
we recover the original Schur $P$-function $P_\lambda(\vectx)$ and Schur $Q$-function $Q_\lambda(\vectx)$ 
by setting $f_d(u) = u^d$ and $f_d(u) = 2 u^d$ respectively.
\item[(2)]
It follows from Nimmo-type formula \cite[Theorem~3.2]{Ivanov2} that 
Ivanov's factorial $P$-function $P_\lambda(\vectx|\vecta)$ and $Q$-function $Q_\lambda(\vectx|\vecta)$ 
are obtained by taking $f_d(u) = (u|\vecta)^d$ and $f_d(u) = 2 (u|\vecta)^d$ respectively.
\item[(3)]
As we will see in Section~7, our generalized $P$-functions include the $t=-1$ specializations 
of Hall--Littlewood functions associated to the root system of type $B$, $C$ and $D$.
\end{enumerate}
\end{example}

Ikeda--Naruse \cite{IN13} and Nakagawa--Naruse \cite{NN} introduced other generalizations 
of factorial $P$- and $Q$-functions from the viewpoint of Schubert calculus.

The organization and main results of this paper are as follows.
In Section~2, we relate our definition of generalized $P$-functions (Definition~\ref{def:genP}) 
with generalizations of two other definitions of Schur $P$-/$Q$-functions.
Namely we prove that $P^\FF_\lambda(\vectx)$ is also obtained by setting $t=-1$ 
in the generalized Hall--Littlewood function associated to a polynomial sequence 
(see Theorem~\ref{thm:N=HL}), 
and that $P^\FF_\lambda(\vectx)$ is expressed as the Pfaffian of the skew-symmetric matrix 
with entries $P^\FF_{(\lambda_i,\lambda_j)}(\vectx)$ (see Theorem~\ref{thm:Schur}).
In Section~3, we introduce the notion of generalized dual $P$-functions $\hat{P}^\FF_\lambda(\vectx)$ 
and prove the Cauchy-type identity.
In Section~4, we define generalized skew $P$-functions $P^\FF_{\lambda/\mu,p}(\vectx)$ 
in terms of J\'ozefiak--Pragacz--Nimmo tpye Pfaffian and prove that $P^\FF_{\lambda/\mu,p}(\vectx)$ 
appears as the coefficient of $P^\FF_\mu(\vecty)$ in the expansion of $P^\FF_\lambda(\vectx,\vecty)$ 
(see Theorem~\ref{thm:JP}).
In Section~5, we consider the modified Pieri coefficients in the expansion of the product 
$P^\FF_\mu(\vectx) \cdot Q_{(r)}(\vectx)$ and obtain a determinant formula 
for the generating function of modified Pieri coefficients (see Theorem~\ref{thm:Pieri}).  
Section~6 focuses on Ivanov's factorial $P$-/$Q$-functions.
We derive a determinant formula of the factorial skew $P$-function 
in one variable (see Theorem~\ref{thm:fac-skewP}),
and an explicit product formula for the generating function of modified Pieri coefficients 
(see Theorem~\ref{thm:fac-Pieri}).
In Section~7, we show that the Hall--Littlewood functions at $t = -1$ 
associated to the classical root systems can be written as 
generalized $P$-functions associated to certain polynomial sequences (see Theorem~\ref{thm:N=HL-BCD}).
Appendix~A collects some Schur-type Pfaffian evaluations and useful formulas.
\section{%
Several expressions of generalized $\boldsymbol{P}$-functions
}

In this section, we give several expressions of generalized $P$-functions 
associated to an admissible polynomial sequence, and study their basic properties.

\subsection{%
Hall--Littlewood-type expression
}

In this subsection we prove that our generalized $P$-functions 
are obtained as the $t=-1$ specialization of Hall--Littlewood-type functions.

We begin with the following proposition.

\begin{prop}
\label{prop:genP}
Let $\FF$ be an admissible sequence of polynomials and $\vectx = (x_1, \dots, x_n)$.
Then we have
\begin{enumerate}
\item[(1)]
For the empty partition $\emptyset$, we have $P^\FF_\emptyset(\vectx) = 1$.
\item[(2)]
If $\lambda$ is a strict partition of length $>n$, then we have $P^\FF_\lambda(\vectx) = 0$.
\end{enumerate}
\end{prop}

\begin{demo}{Proof}
By using the definition of (\ref{eq:Nimmo}), 
we can derive (1) from the Pfaffian evaluations (\ref{eq:Schur-Pf1}) and (\ref{eq:Schur-Pf2}), 
and (2) from Proposition~\ref{prop:Pf-Laplace}.
\end{demo}

We define a generalization of Hall--Littlewood polynomials associated to 
an admissible polynomial sequence.

\begin{definition}
\label{def:genHL}
Let $n$ be a positive integer and $\vectx = (x_1, \dots, x_n)$.
Given a partition $\lambda$ of length $l \le n$, 
we regard $\lambda$ as a sequence $(\lambda_1, \dots, \lambda_l, 0, \dots, 0)$ of length $n$, 
and define a polynomial $v^{(n)}_\lambda(t)$ by putting
$$
v^{(n)}_\lambda(t)
 = 
\prod_{k \ge 0} [m_k]_t!,
$$
where $m_k = \{ i : 1 \le i \le n, \ \lambda_i = k \}$ and $[m]_t! = \prod_{j=1}^m (1-t^j)/(1-t)$.
For an admissible polynomial sequence $\FF$ and a partition $\lambda$ of length $\le n$, 
we define the \emph{generalized Hall--Littlewood function} $\PP^\FF_\lambda(\vectx;t)$ corresponding 
to $\lambda$ by putting
\begin{equation}
\label{eq:genHL}
\PP^{\FF}_\lambda(\vectx;t)
 =
\frac{1}{v^{(n)}_\lambda(t)}
\sum_{w \in S_n}
 w \left(
  \prod_{i=1}^n f_{\lambda_i}(x_i)
  \prod_{1 \le i < j \le n} \frac{x_i - t x_j}{x_i - x_j}
 \right),
\end{equation}
where $S_n$ is the symmetric group acting on $K(t)[x_1, \dots, x_n]$ by permuting variables.
\end{definition}

Setting $f_d(u) = u^d$ for $d \ge 0$, we recover the original Hall--Littlewood polynomials.
The following is the main theorem of this subsection.

\begin{theorem}
\label{thm:N=HL}
For an admissible sequence $\FF$ and a strict partition $\lambda$ of length $l \le n$, 
we have
\begin{equation}
\label{eq:N=HL}
P^\FF_\lambda(\vectx) = \PP^\FF_\lambda(\vectx;-1).
\end{equation}
\end{theorem}

Note that Equation~(\ref{eq:N=HL}) with $f_d(u) = u^d$ is the definition of Schur $P$-function 
adopted in \cite[III.8]{Macdonald95}.
For the sake of completeness and the later use, we give a proof of this theorem, 
which follows the argument in \cite[Appendix]{Nimmo}.
As a first step, we show the following lemma:

\begin{lemma}
\label{lem:HL}
For a strict partition $\lambda$ of length $l \le n$, we have
\begin{align}
\PP^\FF_\lambda(\vectx;-1)
 &=
\sum_{u \in S_n/S_{n-l}}
 u \left(
  \prod_{i=1}^l f_{\lambda_i}(x_i)
  \prod_{\substack{ 1 \le i < j \le n \\ i \le l }} \frac{ x_i + x_j }{ x_i - x_j }
 \right)
\label{eq:HL1}
\\
 &=
\frac{ 1 }{ (n-l)! }
\sum_{u \in S_n}
 u \left(
  \prod_{i=1}^l f_{\lambda_i}(x_i)
  \prod_{\substack{ 1 \le i < j \le n \\ i \le l }} \frac{ x_i + x_j }{ x_i - x_j }
 \right),
\label{eq:HL2}
\end{align}
where $S_{n-l}$ is the symmetric group on the last $n-l$ variables $x_{l+1}, \dots, x_n$.
\end{lemma}

\begin{demo}{Proof}
Since $f_0(u) = 1$ and the product $\prod_{1 \le i < j \le n, i \le l} (x_i - t x_j)/(x_i - x_j)$ 
is invariant under $S_{n-l}$, we have
\begin{multline*}
\PP_\lambda(\vectx;t)
\\
 =
\frac{ 1 }
     { v^{(n)}_\lambda(t) }
\sum_{ w' \in S_n/S_{n-l} }
w' \left(
 \prod_{i=1}^l f_{\lambda_i}(x_i)
 \prod_{\substack{ 1 \le i < j \le n \\ i \le l }} \frac{ x_i - t x_j }{ x_i - x_j }
 \sum_{w'' \in S_{n-l}}
  w'' \left(
   \prod_{l+1 \le i < j \le n} \frac{ x_i - t x_j }{ x_i - x_j }
  \right)
\right).
\end{multline*}
By using (see \cite[Theorem~2.8]{Macdonald72})
$$
\sum_{w'' \in S_{n-l}}
 w'' \left(
  \prod_{l+1 \le i < j \le n} \frac{ x_i - t x_j }{ x_i - x_j }
 \right)
 =
[n-l]_t!
 =
v^{(n)}_\lambda(t),
$$
we have
$$
\PP_\lambda(\vectx;t)
 =
\sum_{ w' \in S_n/S_{n-l} }
 w' \left(
  \prod_{i=1}^l f_{\lambda_i}(x_i)
  \prod_{\substack{ 1 \le i < j \le n \\ i \le l }} \frac{ x_i - t x_j }{ x_i - x_j }
 \right).
$$
By specializing $t=-1$, we obtain (\ref{eq:HL1}), from which (\ref{eq:HL2}) follows.
\end{demo}

\begin{demo}{Proof of Theorem~\ref{thm:N=HL}}
Since $\Delta(\vectx)$ is alternating in $x_1, \dots, x_n$, it follows from (\ref{eq:HL1}) that
$$
\PP_\lambda(\vectx;-1)
 =
\frac{ (-1)^{\binom{n}{2} + \binom{n-l}{2}} }
     { \Delta(\vectx) }
\sum_{v \in S_n/S_{n-l}} 
 \sgn(v) v \left(
  \prod_{i=1}^l f_{\lambda_i}(x_i)
  \prod_{l+1 \le i < j \le n} \frac{ x_j - x_i }{ x_j + x_i }
 \right).
$$
Since $\prod_{l+1 \le i < j \le n} (x_j - x_i)/(x_j + x_i)$ is invariant under the symmetric group $S_l$ 
acting on the first $l$ variables $x_1, \dots, x_l$, we have
\begin{align*}
&
\PP_\lambda(\vectx;-1)
\\
 &=
\frac{ (-1)^{\binom{n}{2} + \binom{n-l}{2}} }
     { \Delta(\vectx) }
\sum_{v' \in S_n/(S_l \times S_{n-l})} 
 \sgn(v') v' \left(
  \sum_{v'' \in S_l} \sgn(v'')
   v'' \left(
    \prod_{i=1}^l f_{\lambda_i}(x_i)
   \right)
  \prod_{l+1 \le i < j \le n} \frac{ x_j - x_i }{ x_j + x_i }
 \right)
\\
 &=
\frac{ (-1)^{\binom{n}{2} + \binom{n-l}{2}} }
     { \Delta(\vectx) }
\sum_{v' \in S_n/(S_l \times S_{n-l})} 
 \sgn(v') v' \left(
  \det \left( f_{\lambda_j}(x_i) \right)_{1 \le i, j \le l}
  \prod_{l+1 \le i < j \le n} \frac{ x_j - x_i }{ x_j + x_i }
 \right).
\end{align*}
We take $R = \{ u \in S_n : u(1) < \dots < u(l), u(l+1) < \dots < u(n) \}$ 
as a complete set of coset representatives of $S_n/(S_l \times S_{n-l})$.
We note that the correspondence $u \mapsto \{ u(l+1), \dots, u(n) \}$ gives a bijection 
between the coset representatives $R$ and the set $\binom{[n]}{n-l}$ of all $(n-l)$-element subsets 
of $[n] = \{ 1, \dots, n \}$.

First we consider the case where $n-l$ is even.
In this case, by using Schur's Pfaffian evaluation (\ref{eq:Schur-Pf1}), we have
$$
\PP_\lambda(\vectx;-1)
 =
\frac{ (-1)^{\binom{n}{2} + \binom{n-l}{2}} }
     { \Delta(\vectx) }
\sum_{u \in S_n/S_l \times S_{n-l}}
 \sgn(u) u \Big(
  \det V_\lambda(\vectx_{[l]})
  \Pf A(\vectx_{[n] \setminus [l]})
 \Big),
$$
where $[l] = \{ 1, \dots, l \}$, $[n] \setminus [l] = \{ l+1, \dots, n \}$ 
and $\vectx_J = (x_{j_1}, \dots, x_{j_m})$ for $J = \{ j_1, \dots, j_m \}$ 
with $j_1 < \dots < j_m$.
On the other hand, by applying Proposition~\ref{prop:Pf-Laplace} (a Pfaffian version of the Laplace expansion) 
to the matrices $Z = A(\vectx)$ and $W = V_\lambda(\vectx)$, we obtain
$$
\Pf \begin{pmatrix}
 A(\vectx) & V_\lambda(\vectx) \\
 -\trans V_\lambda(\vectx) & O
\end{pmatrix}
=
\sum_{I \in \binom{[n]}{n-l}} (-1)^{\Sigma(I) - \binom{n}{2}}
\Pf A(\vectx_I) \det V_\lambda(\vectx_{[n] \setminus I}),
$$
where $\Sigma(I) = \sum_{i \in I} i$. 
Since $n-l$ is even, we can see that, if $u \in R$ corresponds to $I \in \binom{[n]}{n-l}$, 
then the inversion number of $u$ is given by
$$
\inv(u)
 =
\binom{n+1}{2} - \binom{l+1}{2} - \Sigma(I)
 \equiv
\binom{n}{2} - \binom{l}{2} + \Sigma(I)
\quad \bmod 2.
$$
Hence we have
$$
\Pf \begin{pmatrix}
 A(\vectx) & V_\lambda(\vectx) \\
 -\trans V_\lambda(\vectx) & O
\end{pmatrix}
 =
\sum_{u \in S_n/(S_l \times S_{n-l})}
 (-1)^{\binom{l}{2}} \sgn(u) u \Big(
  \det V_\lambda(\vectx_{[l]}) \Pf A(\vectx_{[n] \setminus [l]}) 
 \Big),
$$
and
$$
\PP_\lambda(\vectx;-1)
 =
\frac{ (-1)^{ \binom{n}{2} + \binom{n-l}{2} + \binom{l}{2} } }
     { \Delta(\vectx) }
\Pf \begin{pmatrix}
 A(\vectx) & V_\lambda(\vectx) \\
 -\trans V_\lambda(\vectx) & O
\end{pmatrix}.
$$
Now we can use the relation $\binom{n}{2} - \binom{l}{2} - \binom{n-l}{2} = l (n-l) \equiv 0 \bmod 2$ 
to complete the proof of (\ref{eq:N=HL}) in the case where $n-l$ is even.

Next we consider the case where $n-l$ is odd.
In this case, by using (\ref{eq:Schur-Pf2}), we see that
$$
\PP_\lambda(\vectx;-1)
 =
\frac{ (-1)^{ \binom{n}{2} + \binom{n-l}{2}} }
     { \Delta(\vectx) }
\sum_{u \in S_n/(S_l \times S_{n-l})}
 \sgn(u) u \left(
  \det V_\lambda(\vectx_{[l]})
  \Pf \begin{pmatrix}
   A(\vectx_{[n] \setminus [l]}) & \vectone \\
   - \trans\vectone & 0
  \end{pmatrix}
 \right),
$$
where $\vectone$ is the all-one column vector.
On the other hand, by applying Proposition~\ref{prop:Pf-Laplace} to the matrices
$$
Z = \begin{pmatrix}
 A(\vectx) & \vectone \\
 -\trans\vectone & 0
\end{pmatrix},
\quad
W = \begin{pmatrix}
 V_\lambda(\vectx) \\ O
\end{pmatrix},
$$
we see that
$$
\Pf \begin{pmatrix}
 A(\vectx) & \vectone & V_\lambda(\vectx) \\
 -\trans\vectone & 0 & O \\
 -\trans V_\lambda(\vectx) & O & O
\end{pmatrix}
 =
\sum_I (-1)^{\Sigma(I) - \binom{n+1}{2}}
\Pf Z(I) \det W([n+1] \setminus I;[l]),
$$
where $I$ runs over all $(n+1-l)$-element subsets of $[n+1]$.
If $n+1 \not\in I$, then we have
$$
\det W([n+1] \setminus I;[l])
 =
\det \begin{pmatrix}
 V_\lambda(\vectx_{[n] \setminus I}) \\ O
\end{pmatrix}
 =
0.
$$
Hence we have
$$
\Pf \begin{pmatrix}
 A(\vectx) & \vectone & V_\lambda(\vectx) \\
 -\trans\vectone & 0 & O \\
 -\trans V_\lambda(\vectx) & O & O
\end{pmatrix}
 =
\sum_{I \in \binom{[n]}{n-l}}
 (-1)^{\Sigma(I \cup \{ n+1 \}) - \binom{n+1}{2}}
\Pf \begin{pmatrix}
 A(\vectx_I) & \vectone \\
 -\trans\vectone & 0
\end{pmatrix}
\det V_\lambda(\vectx_{[n] \setminus I}).
$$
Since $n-l$ is odd, we see that, 
if $u \in R$ corresponds to $I \in \binom{[n]}{n-l}$, then we have
$$
\inv (u)
 =
\binom{n+1}{2} - \binom{l+1}{2} - \Sigma(I)
 \equiv
\binom{n+1}{2} - \binom{l}{2} - \Sigma(I \cup \{ n+1 \})
\bmod 2.
$$
Also, by permuting rows/columns we have
$$
\Pf \begin{pmatrix}
 A(\vectx) & \vectone & V_\lambda(\vectx) \\
 -\trans\vectone & 0 & O \\
 -\trans V_\lambda(\vectx) & O & O
\end{pmatrix}
 =
(-1)^l \Pf \begin{pmatrix}
 A(\vectx) & V_{\lambda^0}(\vectx) \\
 -\trans V_{\lambda^0}(\vectx) & O
\end{pmatrix}.
$$
Hence we have
$$
\Pf \begin{pmatrix}
 A(\vectx) & V_{\lambda^0}(\vectx) \\
 -\trans V_{\lambda^0}(\vectx) & O
\end{pmatrix}
 =
\sum_{u \in S_n/(S_l \times S_{n-l})}
 (-1)^{\binom{l}{2} + l} \sgn(u) u \Big(
  \Pf A(\vectx_{[n] \setminus [l]}) \det V_\lambda(\vectx_{[l]})
 \Big),
$$
and
$$
\PP_\lambda(\vectx;-1)
 =
\frac{ (-1)^{ \binom{n}{2} + \binom{n-l}{2} + \binom{l}{2} + l } }
     { \Delta(\vectx) }
\Pf \begin{pmatrix}
 A(\vectx) & V_{\lambda^0}(\vectx) \\
 -\trans V_{\lambda^0}(\vectx) & O
\end{pmatrix}.
$$
Now we can complete the proof in the case where $n-l$ is odd by using the congruence relation
$\binom{n}{2} - \binom{l}{2} - \binom{n-l}{2} = l (n-l) \equiv l \bmod 2$.
\end{demo}

By combining Theorem~\ref{thm:N=HL} and Lemma~\ref{lem:HL}, we obtain

\begin{corollary}
\label{cor:N=HL}
For a strict partition $\lambda$ of length $l \le n$, we have
\begin{align}
P^\FF_\lambda(\vectx)
 &=
\sum_{u \in S_n/S_{n-l}}
 u \left(
  \prod_{i=1}^l f_{\lambda_i}(x_i)
  \prod_{\substack{ 1 \le i < j \le n \\ i \le l }} \frac{ x_i + x_j }{ x_i - x_j }
 \right)
\label{eq:HL3}
\\
 &=
\frac{ 1 }{ (n-l)! }
\sum_{u \in S_n}
 u \left(
  \prod_{i=1}^l f_{\lambda_i}(x_i)
  \prod_{\substack{ 1 \le i < j \le n \\ i \le l }} \frac{ x_i + x_j }{ x_i - x_j }
 \right).
\label{eq:HL4}
\end{align}
\end{corollary}

\subsection{%
Schur-type Pfaffian formula
}

In this subsection we use the definition (\ref{eq:Nimmo}) 
and a Pfaffian version of Sylvester formula (Proposition~\ref{prop:Pf-Sylvester}) 
to derive the Schur-type Pfaffian formula for $P^\FF_\lambda(\vectx)$, 
which generalizes (a part of) Schur's original definition of Schur $Q$-functions \cite[\S~35]{Schur} 
and a similar formula for factorial $Q$-functions \cite[Theorem~9.1]{Ivanov2}.
We use the following conventions:
\begin{gather}
\label{eq:convention1}
P^\FF_{(0)} (\vectx) = 1,
\\
\label{eq:convention2}
P^\FF_{(s,r)}(\vectx) = - P^\FF_{(r,s)}(\vectx),
\quad
P^\FF_{(r,0)}(\vectx) = - P^\FF_{(0,r)}(\vectx) = P^\FF_{(r)}(\vectx),
\quad
P^\FF_{(0,0)}(\vectx) = 0,
\end{gather}
where $r$ and $s$ are positive integers.

\begin{theorem}
\label{thm:Schur}
Let $\FF$ be an admissible sequence.
For a sequence $\alpha = (\alpha_1, \dots, \alpha_r)$ of non-negative integers, 
let $S^\FF_\alpha(\vectx)$ be the $r \times r$ skew-symmetric matrix defined by
\begin{equation}
\label{eq:S}
S^\FF_\alpha(\vectx)
 =
\Big( P^\FF_{(\alpha_i, \alpha_j)}(\vectx) \Big)_{1 \le i, j \le r}.
\end{equation}
Then, for a strict partition $\lambda$ of length $l$, we have
\begin{equation}
\label{eq:Schur}
P^\FF_\lambda(\vectx)
 = 
\begin{cases}
 \Pf S^\FF_\lambda(\vectx) &\text{if $l$ is even,} \\
 \Pf S^\FF_{\lambda^0}(\vectx) & \text{if $l$ is odd,}
\end{cases}
\end{equation}
where $\lambda = (\lambda_1, \dots, \lambda_l)$ and $\lambda^0 = (\lambda_1, \dots, \lambda_l,0)$.
\end{theorem}

In order to prove this theorem, we can use the same argument as in \cite[Theorem~4.1 (3) and Remark~4.3]{Okada}.
As we will see in Proposition~\ref{prop:stability}, the generalized $P$-functions 
do not have the stability property, so we cannot reduce the proof to the case where $n$ is even.

\begin{demo}{Proof}
By applying Proposition~\ref{prop:Pf-Sylvester} to the matrix $X$ given by
$$
X = \begin{cases}
 \begin{pmatrix}
  A(\vectx) & V_\lambda(\vectx) \\
  - \trans V_\lambda(\vectx) & O
 \end{pmatrix}
 &\text{if $n$ is even and $l$ is even,} \\
 \begin{pmatrix}
  A(\vectx) & V_{\lambda^0}(\vectx) \\
  - \trans V_{\lambda^0}(\vectx) & O
 \end{pmatrix}
 &\text{if $n$ is even and $l$ is odd,} \\
\begin{pmatrix}
 A(\vectx) & \vectone & V_\lambda(\vectx) \\
 -\trans\vectone & 0 & O \\
 -\trans V_\lambda(\vectx) & O & O
\end{pmatrix}
 &\text{if $n$ is odd and $l$ is even,}
\\
\begin{pmatrix}
 A(\vectx) & \vectone & V_\lambda(\vectx) & 0 \\
 -\trans\vectone & 0 & O & -1 \\
 -\trans V_\lambda(\vectx) & O & O & 0 \\
 0 & 1 & 0 & 0
\end{pmatrix}
 &\text{if $n$ is odd and $l$ is odd.}
\end{cases}
$$
If $n$ is even, then we have
$$
\frac{ \Pf X }
     { \Pf X([n]) }
 =
P_\lambda(\vectx),
\quad
\frac{ \Pf X([n] \cup \{ n+i, n+j \}) }
     { \Pf X([n]) }
= P_{(\lambda_i,\lambda_j)}(\vectx).
$$
If $n$ is odd and $l$ is even, then by permuting rows/columns, we see that
$$
\frac{ \Pf X }
     { \Pf X([n+1]) }
= P_\lambda(\vectx),
\quad
\frac{ \Pf X([n+1] \cup \{ n+1+i, n+1+j \}) }
     { \Pf X([n+1]) }
= P_{(\lambda_i,\lambda_j)}(\vectx).
$$
If $n$ is odd and $l$ is odd, then by expanding the Pfaffians along the last row/column, we have
$$
\frac{ \Pf X }
     { \Pf X([n+1]) }
= P_\lambda(\vectx),
\quad
\frac{ \Pf X([n+1] \cup \{ n+1+i, n+1+(l+1) \}) }
     { \Pf X([n+1]) }
= P_{(\lambda_i)}(\vectx),
$$
and by permuting rows/columns we see that
$$
\frac{ \Pf X([n+1] \cup \{ n+1+i, n+1+j \}) }
     { \Pf X([n+1]) }
= P_{(\lambda_i,\lambda_j)}(\vectx).
$$
Now Theorem~\ref{thm:Schur} follows immediately from Proposition~\ref{prop:Pf-Sylvester}.
\end{demo}

\subsection{%
Stability
}

The Schur $P$-functions have the stability property (see \cite[III, (2.5)]{Macdonald95})
$$
P_\lambda(x_1, \dots, x_n, 0) = P_\lambda(x_1, \dots, x_n).
$$
Our generalizations $P^\FF_\lambda(\vectx)$ do not have the stability property in general.
For example we can show that
$$
P^\FF_{(r)}(x_1, \dots, x_n, 0)
 =
P^\FF_{(r)}(x_1, \dots, x_n) + (-1)^n f_r(0)
$$
for $r \ge 1$.
The following ``mod $2$ stability property'' was given by 
\cite[Proposition~8.1]{IN09} for factorial $P$-functions.

\begin{prop}
\label{prop:stability}
Let $\FF$ be an admissible sequence, and $\lambda$ a strict partition.
\begin{enumerate}
\item[(1)]
In general, we have
$$
P^\FF_\lambda(x_1, \dots, x_n, 0, 0) = P^\FF_\lambda(x_1, \dots, x_n).
$$
\item[(2)]
If $\FF$ is constant-term free, then we have
$$
P^\FF_\lambda(x_1, \dots, x_n, 0) = P^\FF_\lambda(x_1, \dots, x_n).
$$
\end{enumerate}
\end{prop}

\begin{demo}{Proof}
(1)
Let $\vectx = (x_1, \dots, x_n)$ and $\tilde{\vectx} = (x_1, \dots, x_n, x_{n+1})$.
It follows from the definition (\ref{eq:Nimmo}) that
$$
P_\lambda(x_1, \dots, x_n, x_{n+1}, 0)
 =
\frac{ 1 }
     { (-1)^{n+1} \Delta(\tilde{\vectx}) }
\Pf \begin{pmatrix}
 A(\tilde{\vectx}) & - \vectone & V_{\lambda^*}(\tilde{\vectx}) \\
 \trans\vectone & 0 & V_{\lambda^*}(0) \\
 -\trans V_{\lambda^*}(\tilde{\vectx}) & -\trans V_{\lambda^*}(0) & O
\end{pmatrix},
$$
where $\lambda^* = \lambda$ or $\lambda^0$ according whether $n+l$ is even or odd, 
and $V_\alpha(0)$ is the row vector $\left( f_{\alpha_1}(0), \dots, f_{\alpha_r}(0) \right)$.
Hence, if we put $x_{n+1} = 0$ in the above formula, we have
$$
P_\lambda(x_1, \dots, x_n, 0, 0)
 =
\frac{ 1 }
     { (-1)^{n+1} \cdot (-1)^n \Delta(\vectx) }
\Pf \begin{pmatrix}
 A(\vectx) & -\vectone & -\vectone & V_{\lambda^*}(\vectx) \\
 \trans\vectone & 0 & -1 & V_{\lambda^*}(0) \\
 \trans\vectone & 1 & 0 & V_{\lambda^*}(0) \\
 -\trans V_{\lambda^*}(\vectx) & -\trans V_{\lambda^*}(0) & -\trans V_{\lambda^*}(0) & O
\end{pmatrix}.
$$
By subtracting the $(n+1)$st column/row from the $(n+2)$nd column/row 
and then by expanding the resulting Pfaffian along the $(n+2)$nd columns/row, 
we see that
$$
P_\lambda(x_1, \dots, x_n, 0, 0)
 =
\frac{ 1 }
     { (-1)^{n+1} \cdot (-1)^n \Delta(\vectx) }
\cdot (-1)
\Pf \begin{pmatrix}
 A(\vectx) & V_{\lambda^*}(\vectx) \\
 -\trans V_{\lambda^*}(\vectx) & O
\end{pmatrix}
 =
P_\lambda(\vectx).
$$

(2)
By the definition (\ref{eq:Nimmo}), we have
$$
P_\lambda(x_1, \dots, x_n, 0)
 =
\frac{ 1 }
     { (-1)^n \Delta(\vectx) }
\Pf \begin{pmatrix}
 A(\vectx) & -\vectone & V_{\lambda^*}(\vectx) \\
 \trans \vectone & 0 & V_{\lambda^*}(0) \\
 -\trans V_{\lambda^*}(\vectx) & -\trans V_{\lambda^*}(0) & O
\end{pmatrix},
$$
where $\lambda^* = \lambda^0$ or $\lambda$ according whether $n+l$ is even or odd.
If $n+l$ is even, by adding the $(n+1)$st row/column to the last row/column 
and then by expanding the resulting Pfaffian along the last row/column, 
we see that $P_\lambda(\vectx,0) = P_\lambda(\vectx)$.
If $n+l$ is odd, then by permuting rows/columns, we obtain 
$P_\lambda(\vectx,0) = P_\lambda(\vectx)$.
\end{demo}

\subsection{%
Relation with generalized Schur functions
}

We conclude this section by proving a relation between generalized $P$-functions and generalized Schur functions.

\begin{prop}
Let $\FF$ be an admissible sequence and $\vectx = (x_1, \dots, x_n)$.
For a partition $\mu$ of length $m \le n$, 
let $\mu+\delta_n$ be the strict partition obtained from 
$(\mu_1+n-1, \mu_2+n-2, \dots, \mu_{n-1}+1, \mu_n)$ by removing $0$s, where $\mu_{m+1} = \dots = \mu_n = 0$.
Then we have
$$
P^\FF_{\mu+\delta_n}(\vectx)
 =
\prod_{1 \le i < j \le n} (x_i + x_j)
\cdot
s^\FF_\mu(\vectx),
$$
where $s^\FF_\mu(\vectx)$ is the generalized Schur function given by (\ref{eq:genS}).
\end{prop}

\begin{demo}{Proof}
Since the strict partition $\mu+\delta_n$ has length $n-1$ or $n$, 
it follows from (\ref{eq:HL4}) that
\begin{align*}
P^\FF_{\mu+\delta_n}(\vectx)
 &=
\sum_{w \in S_n}
 w \left(
  \prod_{i=1}^n f_{\mu_i+n-i}(x_i)
  \prod_{1 \le i < j \le n} \frac{x_i+x_j}{x_i-x_j}
 \right)
\\
 &=
\prod_{1 \le i < j \le n} \frac{x_i+x_j}{x_i-x_j}
\sum_{w \in S_n}
 \sgn(w) w \left(
  \prod_{i=1}^n f_{\mu_i+n-i}(x_i)
 \right)
\\
 &=
\prod_{1 \le i < j \le n} (x_i + x_j)
\cdot
s^\FF_\mu(\vectx).
\end{align*}
\end{demo}

\section{%
Dual $\boldsymbol{P}$-functions and Cauchy-type identity
}

In this section, we introduce the dual of generalized $P$-functions 
and prove the Cauchy-type identity for generalized $P$-functions.

\subsection{%
Dual sequences
}

For a nonzero formal power series $g(v) = \sum_{i=0}^\infty b_i v^i \in K[[v]]$, 
the \emph{order} $\ord g$ of $g$ is defined to be the minimum integer $k$ 
such that $b_k \neq 0$.
Let $\langle \quad, \quad \rangle : K[u] \times K[[v]] \to K$ be 
the non-degenerate bilinear pairing defined by
$$
\langle u^i, v^j \rangle
 =
\begin{cases}
 1 &\text{if $i=j=0$,} \\
 1/2 &\text{if $i=j>0$,} \\
 0 &\text{if $i \neq j$.}
\end{cases}
$$

\begin{lemma}
\label{lem:dual}
Let $\FF = \{ f_d \}_{d=0}^\infty$ be an admissible sequence of polynomials.
\begin{enumerate}
\item[(1)]
Let $\hat{\FF} = \{ \hat{f}_d \}_{d=0}^\infty$ be a sequence of formal power series 
$\hat{f}_d(v) \in K[[v]]$ satisfying $\ord \hat{f}_d = d$ for $d \ge 0$.
Then $\langle f_k, \hat{f}_l \rangle = \delta_{k,l}$ for any $k$, $l \ge 0$ 
if and only if
\begin{equation}
\label{eq:dual}
\sum_{k=0}^\infty f_k(u) \hat{f}_k(v) = \frac{ 1 + uv }{ 1 - uv }.
\end{equation}
\item[(2)]
There exists a unique sequence $\hat{\FF} = \{ \hat{f}_d \}_{d=0}^\infty$ satisfying 
$\ord \hat{f}_d = d$ for $d \ge 0$ and the equivalent conditions in (1).
We call such a sequence $\hat{\FF}$ the \emph{dual} of $\FF$.
\item[(3)]
If $\hat{\FF} = \{ \hat{f}_d \}_{d=0}^\infty$ is the dual of $\FF$, then 
$\FF$ is constant-term free if and only if $\hat{f}_0 = 1$.
\end{enumerate}
\end{lemma}

\begin{demo}{Proof}
(1), (2)
We write $f_d(u) = \sum_{i \ge 0} a_{d,i} u^i$ 
and $\hat{f}_d(v) = \sum_{i \ge 0} b_{d,i} v^i$. 
Since $\deg f_d = d$ and $\ord \hat{f}_d = d$, we have 
$a_{d,i} = 0$ for $i \ge {d+1}$ and $b_{d,i} = 0$ for $i \le {d-1}$.
We define $b'_{d,i}$ by putting
$$
b'_{d,i}
 =
\begin{cases}
 b_{d,0} &\text{if $i=0$,} \\
 b_{d,i}/2 &\text{if $i > 0$.}
\end{cases}
$$

Fix a nonnegative integer $N$ and consider two $(N+1) \times (N+1)$ matrices 
$A = (a_{i,j})_{0 \le i, j \le N}$ and $B' = (b'_{i,j})_{0 \le i, j \le N}$.
Then it is easy to see that 
\begin{enumerate}
\item[(a)]
$\langle f_k, \hat{f}_l \rangle = \delta_{k,l}$ for all $0 \le k, l \le N$ 
if and only if $A \cdot \trans B' = I_{N+1}$.
\item[(b)]
$\sum_{k \ge 0} f_k(u) \hat{f}_k(v) = (1+uv)/(1-uv)$ in the quotient ring 
$K[[u,v]]/(u^{N+1}, v^{N+1})$ if and only if $\trans A \cdot B' = I_{N+1}$.
\end{enumerate}
The claims (1) and (2) follow from this observation.

(3)
If $f_d(0) = 0$ for any $d \ge 1$, then by substituting $u = 0$ in (\ref{eq:dual}) 
we obtain $\hat{f}_0(v) = 1$.
Conversely, if $\hat{f}_0 = 1$, then $\langle f_d, \hat{f_0} \rangle = \delta_{d,0}$ 
is equal to the constant term of $f_d$.
\end{demo}

For example, if $f_d(u) = u^d$ for $d \ge 0$, then the dual of $\FF$ is given by
$$
\hat{f}_d(v)
 = 
\begin{cases}
 1 &\text{if $d = 0$,} \\
 2 v^d &\text{if $d \ge 1$.}
\end{cases}
$$
See Lemma~\ref{lem:fac-dual} for the dual of the sequence 
$\{ (u|\vecta)^d \}_{d=0}^\infty$ of factorial monomials.

\subsection{%
Generating functions of generalized $\boldsymbol{P}$-functions
}

For a sequence of variables $\vectx = (x_1, \dots, x_n)$ and another variable $z$, we put
\begin{equation}
\label{eq:Pi}
\Pi_z(\vectx) = \prod_{i=1}^n \frac{ 1 + x_i z }{ 1 - x_i z}.
\end{equation}
Then the generating functions of Schur $Q$-functions $Q_{(r)}(\vectx)$ and $Q_{(r,s)}(\vectx)$ 
are expressed as
\begin{gather}
\label{eq:GF-SchurQ1}
\sum_{r \ge 0} Q_{(r)}(\vectx) z^r
 = 
\Pi_z(\vectx),
\\
\label{eq:GF-SchurQ2}
\sum_{r,s \ge 0} Q_{(r,s)}(\vectx) z^r w^s
 =
\frac{z-w}{z+w}
\big( \Pi_z(\vectx) \Pi_w(\vectx) - 1 \big),
\end{gather}
respectively (see \cite[III, (8.1)]{Macdonald95} and \cite[p.~117]{Stembridge}).
We can generalize these generating functions in terms of the dual sequence $\hat{\FF}$.

\begin{prop}
\label{prop:GF}
Let $\FF = \{ f_d \}_{d=0}^\infty$ be an admissible sequence of polynomials, 
and $\hat{\FF} = \{ \hat{f}_d \}_{d=0}^\infty$ the dual sequence of $\FF$.
Then, under the convention (\ref{eq:convention1}) and (\ref{eq:convention2}), we have
\begin{enumerate}
\item[(1)]
The generating function of generalized $P$-functions $P^\FF_{(r)}(\vectx)$ is given by
\begin{equation}
\label{eq:GF1}
\sum_{r \ge 0} P^\FF_{(r)}(\vectx) \hat{f}_r(z)
 =
\begin{cases}
 \Pi_z(\vectx)
 &\text{if $n$ is odd,}
\\
 \Pi_z(\vectx) + \hat{f}_0(z) - 1
 &\text{if $n$ is even,}
\end{cases}
\end{equation}
under the convention (\ref{eq:convention1}).
\item[(2)]
The generating function of generalized $P$-functions $P^\FF_{(r,s)}(\vectx)$ is given by 
\begin{multline}
\label{eq:GF2}
\sum_{r, s \ge 0} P^\FF_{(r,s)}(\vectx) \hat{f}_r(z) \hat{f}_s(w)
 \\
=
\left\{
\begin{array}{l}
 \dfrac{z-w}{z+w} \big(
  \Pi_z(\vectx) \Pi_w(\vectx) -1
 \big)
 + (\hat{f}_0(w)-1) \Pi_z(\vectx)
 - (\hat{f}_0(z)-1) \Pi_w(\vectx)
\\
\hfill\text{if $n$ is odd,}
\\
 \dfrac{z-w}{z+w} \big(
  \Pi_z(\vectx) \Pi_w(\vectx) -1
 \big)
 \hfill\text{if $n$ is even,}
\end{array}
\right.
\end{multline}
under the convention (\ref{eq:convention2}).
\end{enumerate}
\end{prop}

If we set $f_d(u) = (u|\vecta)^d$ (factorial monomial) with $a_0 = 0$, 
then it follows from Lemma~\ref{lem:fac-dual} that the formulas 
(\ref{eq:GF1}) and (\ref{eq:GF2}) reduce to the formulas given in 
\cite[Theorem~8.2]{Ivanov2} and \cite[Theorem~8.4]{Ivanov2} for factorial $P$-functions 
respectively.

\begin{demo}{Proof}
The idea of the proof is similar to that of \cite[Theorem~4.1 (1) and (2)]{Okada}.
Let $B_z(\vectx)$ be the column vector with $i$th entry $(1 + x_i z)/(1 - x_i z)$.
By (\ref{eq:dual}) we have
$$
\sum_{r \ge 0} \hat{f}_r(z) V_{(r)}(\vectx) = B_z(\vectx),
\quad
\sum_{r \ge 1} \hat{f}_r(z) V_{(r)}(\vectx) = B_z(\vectx) - \hat{f}_0(z) \vectone.
$$

(1)
If $n$ is odd, then by the definition (\ref{eq:Nimmo}) and the multilinearity of Pfaffians, 
we have
$$
\sum_{r \ge 0} P_{(r)}(\vectx) \hat{f}_r(z)
 =
\frac{ 1 }
     { \Delta(\vectx) }
\Pf \begin{pmatrix}
 A(\vectx) & B_z(\vectx) \\
 -\trans B_z(\vectx) & O
\end{pmatrix},
$$
which equals to $\Pi_z(\vectx)$ by (\ref{eq:Schur-Pf3}).
Similarly, if $n$ is even, then we have
$$
\sum_{r \ge 0} P_{(r)}(\vectx) \hat{f}_r(z)
 =
\hat{f}_0(z)
 +
\frac{ 1 }
     { \Delta(\vectx) }
\Pf \begin{pmatrix}
 A(\vectx) & B_z(\vectx) - \hat{f}_0(z) \vectone & \vectone \\
 -\trans B_z(\vectx) + \hat{f}_0(z) \trans\vectone & 0 & 0 \\
 -\trans\vectone & 0 & 0
\end{pmatrix}.
$$
By adding the last row/column multiplied by $\hat{f}_0(z)$ to the second to the last row/column, 
and then by using (\ref{eq:Schur-Pf4}) with $w=0$, 
we obtain
$$
\sum_{r \ge 0} P_{(r)}(\vectx) \hat{f}_r(z)
 =
\hat{f}_0(z)
 +
\frac{ 1 }
     { \Delta(\vectx) }
\Pf \begin{pmatrix}
 A(\vectx) & B_z(\vectx) & \vectone \\
 -\trans B_z(\vectx) & 0 & 0 \\
 -\trans\vectone & 0 & 0
\end{pmatrix}
 =
\hat{f}_0(z)
 +
\left( \Pi_z(\vectx) - 1 \right).
$$

(2)
If $n$ is even, then we have
$$
\sum_{r, s \ge 0} P_{(r,s)}(\vectx) \hat{f}_r(z) \hat{f}_s(w)
 =
\frac{ 1 }
     { \Delta(\vectx) }
\Pf \begin{pmatrix}
 A(\vectx) & B_z(\vectx) & B_w(\vectx) \\
 -\trans B_z(\vectx) & 0 & 0 \\
 -\trans B_w(\vectx) & 0 & 0
\end{pmatrix},
$$
and this equals to $\Pi_z(\vectx) \Pi_w(\vectx) - 1$ by (\ref{eq:Schur-Pf4}).
If $n$ is odd, then we have
\begin{align*}
&
\sum_{r,s \ge 0} P_{(r,s)}(\vectx) \hat{f}_r(z) \hat{f}_s(w)
\\
 &=
\sum_{r,s > 0} P_{(r,s)}(\vectx) \hat{f}_r(z) \hat{f}_s(w)
 +
\sum_{r \ge 0} P_{(r,0)}(\vectx) \hat{f}_r(z) \hat{f}_0(w)
 +
\sum_{s \ge 0} P_{(0,s)}(\vectx) \hat{f}_0(z) \hat{f}_s(w)
\\
 &=
\frac{1}{\Delta(\vectx)}
\Pf \begin{pmatrix}
 A(\vectx) & B_z(\vectx) - \hat{f}_0(z) \vectone & B_w(\vectx) - \hat{f}_0(w) \vectone & \vectone \\
 -\trans B_z(\vectx) + \hat{f}_0(z) \trans\vectone & 0 & 0 & 0 \\
 -\trans B_w(\vectx) + \hat{f}_0(w) \trans\vectone & 0 & 0 & 0 \\
 -\trans\vectone & 0 & 0 & 0
\end{pmatrix}
\\
&\quad
+
\hat{f}_0(w)
\cdot
\frac{1}{\Delta(\vectx)}
\Pf \begin{pmatrix}
 A(\vectx) & B_z(\vectx) - \hat{f}_0(z) \vectone \\
 -\trans B_z(\vectx) + \hat{f}_0(z) \trans\vectone & 0
\end{pmatrix}
\\
&\quad
-
\hat{f}_0(z)
\cdot
\frac{1}{\Delta(\vectx)}
\Pf \begin{pmatrix}
 A(\vectx) & B_w(\vectx) - \hat{f}_0(w) \vectone \\
 -\trans B_w(\vectx) + \hat{f}_0(w) \trans\vectone & 0
\end{pmatrix}.
\end{align*}
By adding the last row/column multiplied by $\hat{f}_0(z)$ (resp. $\hat{f}_0(w)$) 
to the third (resp. second) to the last row/column in the first Pfaffian, 
we see that
\begin{multline*}
\Pf \begin{pmatrix}
 A(\vectx) & B_z(\vectx) - \hat{f}_0(z) \vectone & B_w(\vectx) - \hat{f}_0(w) \vectone & \vectone \\
 -\trans B_z(\vectx) + \hat{f}_0(z) \trans\vectone & 0 & 0 & 0 \\
 -\trans B_w(\vectx) + \hat{f}_0(w) \trans\vectone & 0 & 0 & 0 \\
 -\trans\vectone & 0 & 0 & 0
\end{pmatrix}
\\
 =
\Pf \begin{pmatrix}
 A(\vectx) & B_z(\vectx) & B_w(\vectx) & \vectone \\
 -\trans B_z(\vectx) & 0 & 0 & 0 \\
 -\trans B_w(\vectx) & 0 & 0 & 0 \\
 -\trans\vectone & 0 & 0 & 0
\end{pmatrix}.
\end{multline*}
By using the multilinearity of Pfaffians, we have
\begin{multline*}
\Pf \begin{pmatrix}
 A(\vectx) & B_z(\vectx) - \hat{f}_0(z) \vectone \\
 -\trans B_z(\vectx) + \hat{f}_0(z) \trans\vectone & 0
\end{pmatrix}
\\
=
\Pf \begin{pmatrix}
 A(\vectx) & B_z(\vectx) \\
 -\trans B_z(\vectx) & 0
\end{pmatrix}
-
\hat{f}_0(z)
\Pf \begin{pmatrix}
 A(\vectx) & \vectone \\
 -\trans\vectone & 0
\end{pmatrix},
\end{multline*}
\begin{multline*}
\Pf \begin{pmatrix}
 A(\vectx) & B_w(\vectx) - \hat{f}_0(w) \vectone \\
 -\trans B_w(\vectx) + \hat{f}_0(w) \trans\vectone & 0
\end{pmatrix}
\\
 =
\Pf \begin{pmatrix}
 A(\vectx) & B_w(\vectx) \\
 -\trans B_w(\vectx) & 0
\end{pmatrix}
-
\hat{f}_0(w)
\Pf \begin{pmatrix}
 A(\vectx) & \vectone \\
 -\trans\vectone & 0
\end{pmatrix}.
\end{multline*}
Hence we can use (\ref{eq:Schur-Pf5}), (\ref{eq:Schur-Pf3}) and (\ref{eq:Schur-Pf2}) 
to evaluate these Pfaffians and complete the proof.
\end{demo}

\subsection{%
Dual $\boldsymbol{P}$-functions and Cauchy-type identity
}

In this section we introduce the dual $P$-functions and prove the Cauchy-type identity.

\begin{definition}
\label{def:dualP}
Let $\FF = \{ f_d \}_{d=0}^\infty$ be an admissible sequence 
and $\hat{\FF} = \{ \hat{f}_d \}_{d=0}^\infty$ the dual of $\FF$.
For a sequence $\alpha = (\alpha_1, \dots, \alpha_r)$ of nonnegative integers, 
let $\hat{V}^\FF_\alpha(\vectx)$ be the $n \times r$ matrix given by
$$
\hat{V}^\FF_\alpha(\vectx) = \left( \hat{f}_{\alpha_j}(x_i) \right)_{1 \le i \le n, 1 \le j \le r}.
$$
Given a strict partition $\lambda$ of length $l$,
we define the \emph{generalized dual $P$-functions} $\hat{P}^{\FF}_\lambda(\vectx)$ 
by putting
\begin{equation}
\label{eq:dualNimmo}
\hat{P}^\FF_\lambda(\vectx)
 =
\begin{cases}
 \dfrac{ 1 }{ \Delta(\vectx) }
 \Pf \begin{pmatrix}
  A(\vectx) & \hat{V}^\FF_\lambda(\vectx) \\
  - \trans \hat{V}^\FF_\lambda(\vectx) & O
 \end{pmatrix}
 &\text{if $n+l$ is even,}
\\
 \dfrac{ 1 }{ \Delta(\vectx) }
 \Pf \begin{pmatrix}
  A(\vectx) & \hat{V}^\FF_{\lambda^0}(\vectx) \\
  - \trans \hat{V}^\FF_{\lambda^0}(\vectx) & O
 \end{pmatrix}
 &\text{if $n+l$ is odd,} \\
\end{cases}
\end{equation}
where $\lambda^0 = (\lambda_1, \dots, \lambda_l, 0)$.
\end{definition}

If $f_d(u) = u^d$ for $d \ge 0$, then $\hat{f}_d(v) = 2 v^d$ for $d \ge 0$.
Hence $P^\FF_\lambda (\vectx) = P_\lambda(\vectx)$ is the original Schur $P$-function, 
and the dual $P$-function $\hat{P}^\FF_\lambda(\vectx) = Q_\lambda(\vectx)$ 
is the original Schur $Q$-function.

By the same arguments as in Section~2, we obtain the following Proposition.

\begin{prop}
Suppose that $\FF = \{ f_d \}_{d=0}^\infty$ be a constant-term free admissible sequence 
with dual $\hat{\FF} = \{ \hat{f}_d \}_{d=0}^\infty$.
Then we have
\begin{enumerate}
\item[(1)]
For a strict partition of length $l \le n$, we have
\begin{equation}
\label{eq:dualHL}
\hat{P}^\FF_\lambda(\vectx)
 =
\frac{ 1 }{ (n-l)! }
\sum_{u \in S_n}
 u \left(
  \prod_{i=1}^l \hat{f}_{\lambda_i}(x_i)
  \prod_{\substack{ 1 \le i < j \le n \\ i \le l }} \frac{ x_i + x_j }{ x_i - x_j }
 \right).
\end{equation}
\item[(2)]
For a sequence $\alpha = (\alpha_1, \dots, \alpha_r)$ of non-negative integers, 
let $\hat{S}^\FF_\alpha(\vectx)$ be the $r \times r$ skew-symmetric matrix defined by
$$
\hat{S}^\FF_\alpha(\vectx)
 =
\left( \hat{P}^\FF_{(\alpha_i, \alpha_j)}(\vectx) \right)_{1 \le i, j \le r},
$$
where we use the same convention as (\ref{eq:convention2}).
Then, for a strict partition $\lambda$ of length $l$, we have
\begin{equation}
\label{eq:dualSchur}
\hat{P}^\FF_\lambda(\vectx)
 = 
\begin{cases}
 \Pf \hat{S}^\FF_\lambda(\vectx) &\text{if $l$ is even,} \\
 \Pf \hat{S}^\FF_{\lambda^0}(\vectx) & \text{if $l$ is odd.}
\end{cases}
\end{equation}
\item[(3)]
For a strict partition $\lambda$, we have
\begin{equation}
\label{eq:dualstability}
\hat{P}^\FF_\lambda(x_1, \dots, x_n, 0) = \hat{P}^\FF_\lambda(x_1, \dots, x_n).
\end{equation}
\end{enumerate}
\end{prop}

\begin{demo}{Proof}
By Lemma~\ref{lem:dual} (3), we have $\hat{f}_0 = 1$ for a constant-term free admissible sequence $\FF$.
Hence the proofs of Theorems~\ref{thm:N=HL}, \ref{thm:Schur} and Proposition~\ref{prop:stability} (2) 
work literally in this dual setting.
\end{demo}

For Schur's $P$- and $Q$-functions we have the following Cauchy-type identity 
(see \cite[III, (8.13)]{Macdonald95}, and \cite[Theorem~5.1]{Okada} for a linear algebraic proof):
$$
\sum_\lambda
 P_\lambda(\vectx) Q_\lambda(\vecty)
=
\prod_{i,j=1}^n
 \frac{ 1 + x_i y_j }
      { 1 - x_i y_j },
$$
where $\vectx = (x_1, \dots, x_n)$, $\vecty = (y_1, \dots, y_n)$ 
and the summation is taken over all strict partitions of length $\le n$.
We can use the notion of dual $P$-functions to formulate the Cauchy-type identity 
for generalized $P$-functions.

\begin{theorem}
\label{thm:Cauchy}
Let $\FF$ be an admissible sequence of polynomials, 
and let $\vectx = (x_1, \dots, x_n)$ and $\vecty = (y_1, \dots, y_n)$ be 
two sets of indeterminates.
Then we have
$$
\sum_\lambda P^\FF_\lambda(\vectx) \hat{P}^\FF_\lambda(\vecty)
 =
\prod_{i,j=1}^n
 \frac{ 1 + x_i y_j }
      { 1 - x_i y_j },
$$
where $\lambda$ runs over all strict partitions.
\end{theorem}

\begin{demo}{Proof}
We use the same argument as in the proof of \cite[Theorem~5.1]{Okada}.
Apply a Pfaffian version of Cauchy--Binet formula (\ref{eq:Pf-CB2}) 
to the matrices
$$
A = A(\vectx),
\quad
B = A(\vecty),
\quad
S = \Bigl( f_k(x_i) \Bigr)_{1 \le i \le n, k \ge 0},
\quad
T = \Bigl( \hat{f}_k(y_i) \Bigr)_{1 \le i \le n, k \ge 0}.
$$
Strict partitions $\lambda$ are in bijection with subsets 
of $\Nat$ satisfying $\# I = n \bmod 2$ 
via the correspondence $\lambda \mapsto 
I = \{ \lambda_1, \dots, \lambda_{l(\lambda)} \}$ or $\{ \lambda_1, \dots, \lambda_{l(\lambda)}, 0 \}$. 
And we have
\begin{align*}
P^\FF_\lambda(\vectx)
 &=
\frac{(-1)^{\binom{\# I}{2}}}{\Delta(\vectx)} 
\Pf \begin{pmatrix}
 A(\vectx) & S([n];I) \\
 -\trans S([n];I) & O
\end{pmatrix},
\\
\hat{P}^\FF_\lambda(\vecty)
 &=
\frac{(-1)^{\binom{\# I}{2}}}{\Delta(\vecty)} 
\Pf \begin{pmatrix}
 A(\vecty) & T([n];I) \\
 -\trans T([n];I) & O
\end{pmatrix},
\end{align*}
where $S([n];I)$ and $T([n];I)$ are the submatrices of $S$ and $T$ 
obtained by picking up columns indexed by $I$ respectively.
Since the $(i,j)$ entry of $S \trans T$ is equal to $(1+x_iy_j)/(1-x_iy_j)$ 
by (\ref{eq:dual}), we can complete the proof by using the Pfaffian evaluation (\ref{eq:Schur-Pf3}).
\end{demo}
\section{%
Generalized skew $\boldsymbol{P}$-functions
}

In this section, we introduce generalized skew $P$-function 
in terms of J\'ozefiak--Pragacz--Nimmo-type Pfaffians, and study their properties.

\subsection{%
J\'ozefiak--Pragacz--Nimmo-type formula
}

First we define generalized skew $P$-functions associated to an admissible sequence.

\begin{definition}
\label{def:skewP}
Let $\FF = \{ f_d \}_{d=0}^\infty$ be an admissible sequence.
For a pair of nonnegative integers $r$ and $k$, we define a symmetric polynomial $R^\FF_{r/k}(\vectx)$ 
by the relation
\begin{equation}
\label{eq:skewP1}
P^\FF_{(r)}(x_1, \dots, x_n, y)
 =
\sum_{k=0}^\infty R^\FF_{r/k}(x_1, \dots, x_n) f_k(y).
\end{equation}
For two sequences $\alpha = (\alpha_1, \dots, \alpha_r)$
 and $\beta = (\beta_1, \dots, \beta_s)$ of nonnegative integers, 
let $M_{\alpha/\beta}(\vectx)$ be the $r \times s$ matrix given by
\begin{equation}
\label{eq:M}
M^\FF_{\alpha/\beta}(\vectx)
 =
\big( R^\FF_{\alpha_i/\beta_{s+1-j}}(\vectx) \big)_{1 \le i \le r, 1 \le j \le s}.
\end{equation}
For a pair of strict partitions $\lambda$ of length $l$ and $\mu$ of length $m$ 
and a positive integer $p$, 
we define the \emph{generalized skew $P$-function} $P^\FF_{\lambda/\mu,p}(\vectx)$ by putting
\begin{equation}
\label{eq:skewP}
P^\FF_{\lambda/\mu,p}(\vectx)
 =
\begin{cases}
 \Pf \begin{pmatrix}
  S^\FF_\lambda(\vectx) & M^\FF_{\lambda/\mu}(\vectx) \\
  -\trans M^\FF_{\lambda/\mu}(\vectx) & O
  \end{pmatrix}
 &\text{if $l \equiv p$ and $m \equiv p \bmod 2$,}
\\
 \Pf \begin{pmatrix}
  S^\FF_\lambda(\vectx) & M^\FF_{\lambda/\mu^0}(\vectx) \\
  -\trans M^\FF_{\lambda/\mu^0}(\vectx) & O
 \end{pmatrix}
 &\text{if $l \equiv p$ and $m \not\equiv p \bmod 2$, }
\\
 \Pf \begin{pmatrix}
  S^\FF_{\lambda^0}(\vectx) & M^\FF_{\lambda^0/\mu}(\vectx) \\
  -\trans M^\FF_{\lambda^0/\mu}(\vectx) & O
 \end{pmatrix}
 &\text{if $l \not\equiv p$ and $m \equiv p \bmod 2$,}
\\
 \Pf \begin{pmatrix}
  S^\FF_{\lambda^0}(\vectx) & M^\FF_{\lambda^0/\mu^0}(\vectx) \\
  -\trans M^\FF_{\lambda^0/\mu^0}(\vectx) & O
 \end{pmatrix}
 &\text{if $l \not\equiv p$ and $m \not\equiv p \bmod 2$,}
\end{cases}
\end{equation}
where $S^\FF_\alpha(\vectx)$ is the skew-symmetric matrix defined in (\ref{eq:S}), 
and $\lambda^0 = (\lambda_1, \dots, \lambda_l, 0)$, $\mu^0 = (\mu_1, \dots, \mu_m, 0)$.
\end{definition}

The main result of this section is the following theorem, which is a generalization of 
J\'ozefiak--Pragacz--Nimmo formula for Schur $P$/$Q$ functions 
(see \cite[Theorem~1]{JP} and \cite[(2.22)]{Nimmo}).

\begin{theorem}
\label{thm:JP}
For two sets of variables $\vectx = (x_1, \dots, x_n)$ and $\vecty = (y_1, \dots, y_p)$, 
we have
\begin{equation}
\label{eq:JP}
P^\FF_\lambda(\vectx, \vecty)
 =
\sum_\mu
 P^\FF_{\lambda/\mu,p}(\vectx) P^\FF_\mu(\vecty),
\end{equation}
where $\mu$ runs over all strict partitions.
\end{theorem}

We postpone the proof of this theorem to the next subsection.
Before the proof, we derive several properties of generalized skew $P$-functions 
from the definition (\ref{eq:skewP}).
We begin with the following property of $R^\FF_{r/k}(\vectx)$.

\begin{lemma}
\label{lem:R}
For a positive integer $r$, the generalized $P$-function $P^\FF_{(r)}(\vectx,y)$ 
has degree at most $r$ in $y$.
Hence we have $R^\FF_{r/k}(\vectx) = 0$ unless $r \ge k$.
\end{lemma}

\begin{demo}{Proof}
The coefficient of $z^r$ in $\Pi_z(\vectx,y) = \prod_{i=1}^n (1+x_iz)/(1-x_iz) \cdot (1+yz)/(1-yz)$ 
has degree at most $r$ in $y$.
On the other hand, since $\ord \hat{f}_r = r$, 
the coefficient of $z^r$ in $\sum_{r \ge 0} P^\FF_{(r)}(\vectx,y) \hat{f}_r(z)$ is 
a linear combination of $P^\FF_{(0)}(\vectx,y), \dots, P^\FF_{(r)}(\vectx,y)$ with nonzero 
coefficient for $P^\FF_{(r)}(\vectx,y)$.
Hence, by using (\ref{eq:GF1}) and the induction on $r$, we can conclude 
that $P^\FF_{(r)}(\vectx,y)$ has degree at most $r$ in $y$.
\end{demo}

This lemma can be used to prove the following vanishing property of generalized skew $P$-functions.
For a strict partition $\lambda$, we define its shifted diagram $S(\lambda)$ by putting
$$
S(\lambda) = \{ (i,j) \in \Int^2 : i \le j \le \lambda_i +i-1 \}.
$$
For two strict partitions $\lambda$ and $\mu$, 
we write $\lambda \supset \mu$ if $S(\lambda) \supset S(\mu)$.

\begin{prop}
\label{prop:skewP-1}
Let $\FF$ be an admissible sequence.
For two strict partitions $\lambda$ and $\mu$, we have $P^\FF_{\lambda/\mu,p}(\vectx) = 0$ 
unless $\lambda \supset \mu$.
\end{prop}

\begin{demo}{Proof}
Suppose that there exists an index $k$ such that $\lambda_k < \mu_k$.
Then we have, if $i \ge k$ and $j \le k$, 
then $\lambda_i \le \lambda_k < \mu_k \le \mu_j$ 
and $R_{\lambda_i/\mu_j}(\vectx) = 0$ Lemma~\ref{lem:R}.
Hence the skew-symmetric matrices $X$ appearing in the definition (\ref{eq:skewP}) 
of $P^\FF_{\lambda/\mu,p}(\vectx)$ are of the form
$$
X = \begin{pmatrix}
 Z & W \\
 -\trans W & O
\end{pmatrix},
$$
where $W$ has $k$ columns and at most $(k-1)$ nonzero rows.
Since all the $k \times k$ minors of $W$ vanish, 
it follows from the Laplace expansion (Proposition~\ref{prop:Pf-Laplace}) 
that $P^\FF_{\lambda/\mu,p}(\vectx) = \Pf X = 0$.
\end{demo}

If the polynomial sequence $\FF$ is constant-term free, 
then the skew $P$-function $P^\FF_{\lambda/\mu,p}$ is independent of $p$ 
and some formulas have simple forms.

\begin{prop}
\label{prop:skewP-2}
Suppose that $\FF$ is a constant-term free admissible sequence.
\begin{enumerate}
\item[(1)]
$R^\FF_{0/0}(\vectx) = 1$ and $R^\FF_{r/0}(\vectx) = P^\FF_{(r)}(\vectx)$ for a positive integer $r$.
\item[(2)]
For two strict partitions $\lambda$ of length $l$ and $\mu$ of length $m$, we have
$$
P^\FF_{\lambda/\mu,p}(\vectx)
 =
\left\{
\begin{array}{l}
 \Pf \begin{pmatrix}
  S^\FF_\lambda(\vectx) & M^\FF_{\lambda/\mu}(\vectx) \\
  -\trans M^\FF_{\lambda/\mu}(\vectx) & O
 \end{pmatrix}
 \hfill\text{if $l \equiv m \bmod 2$,}
\\
 \Pf \begin{pmatrix}
  S^\FF_\lambda(\vectx) & M^\FF_{\lambda/\mu^0}(\vectx) \\
  -\trans M^\FF_{\lambda/\mu^0}(\vectx) & O
 \end{pmatrix}
 =
 \Pf \begin{pmatrix}
  S^\FF_{\lambda^0}(\vectx) & M^\FF_{\lambda^0/\mu}(\vectx) \\
  -\trans M^\FF_{\lambda^0/\mu}(\vectx) & O
 \end{pmatrix}
\\
 \hfill\text{if $l \not\equiv m \bmod 2$.}
\end{array}
\right.
$$
In particular, $P^\FF_{\lambda/\mu,p}(\vectx)$ is independent of $p$, 
so we write $P^\FF_{\lambda/\mu}(\vectx)$ in this case.
\end{enumerate}
\end{prop}

\begin{demo}{Proof}
(1)
Substituting $y = 0$ in the definition (\ref{eq:skewP1}), we obtain
$$
P_{(r)}(\vectx,0) = \sum_{k=0}^r R_{r/k}(\vectx) f_k(0).
$$
Since $P_{(r)}(\vectx,0) = P_{(r)}(\vectx)$ by the stability property (Proposition~\ref{prop:stability} (2)), 
and $f_r(0) = \delta_{r,0}$ by assumption, 
we have $R_{r/0}(\vectx) = P_{(r)}(\vectx)$.

(2)
By using (1), we see that
$$
\begin{pmatrix}
 S_\lambda(\vectx) & M_{\lambda/\mu^0}(\vectx) \\
 -\trans M_{\lambda/\mu^0}(\vectx) & O 
\end{pmatrix}
 =
\begin{pmatrix}
 S_{\lambda^0}(\vectx) & M_{\lambda^0/\mu}(\vectx) \\
 -\trans M_{\lambda^0/\mu}(\vectx) & O 
\end{pmatrix}.
$$
It remains to show that, if $l(\lambda) \not\equiv p$ and $l(\mu) \not\equiv p \bmod 2$, then
$$
\Pf \begin{pmatrix}
 S_{\lambda^0}(\vectx) & M_{\lambda^0/\mu^0}(\vectx) \\
 -\trans M_{\lambda^0/\mu^0}(\vectx) & O
\end{pmatrix}
 =
\Pf \begin{pmatrix}
 S_\lambda(\vectx) & M_{\lambda/\mu}(\vectx) \\
 -\trans M_{\lambda/\mu}(\vectx) & O
\end{pmatrix}.
$$
Let $T_\lambda(\vectx)$ (resp. $T_{\lambda/0}(\vectx)$) be the column vector 
with $i$th entry $P_{(\lambda_i)}(\vectx)$ (resp. $R_{\lambda_i/0}(\vectx)$).
Since $R_{0/0}(\vectx) = 1$ and $R_{0/k}(\vectx) = 0$ for $k \ge 1$, we have
$$
\begin{pmatrix}
 S_{\lambda^0}(\vectx) & M_{\lambda^0/\mu^0}(\vectx) \\
 -\trans M_{\lambda^0/\mu^0}(\vectx) & O 
\end{pmatrix}
=
\begin{pmatrix}
 S_\lambda(\vectx) & T_\lambda(\vectx) & T_{\lambda/0}(\vectx) & M_{\lambda/\mu}(\vectx) \\
 -\trans T_\lambda(\vectx) & 0 & 1 & O \\
 -\trans T_{\lambda/0}(\vectx) & -1 & 0 & O \\
 -\trans M_{\lambda/\mu}(\vectx) & O & O & O
\end{pmatrix}.
$$
Since $T_\lambda(\vectx) = T_{\lambda/0}(\vectx)$ by (1), we add the $(l+1)$st row/column 
multiplied by $-1$ to the $(l+2)$nd row/column 
and expanding the resulting Pfaffian along the $(l+2)$nd row/column to obtain
$$
\Pf \begin{pmatrix}
 S_{\lambda^0}(\vectx) & M_{\lambda^0/\mu^0}(\vectx) \\
 -\trans M_{\lambda^0/\mu^0}(\vectx) & O
\end{pmatrix}
 =
\Pf \begin{pmatrix}
 S_\lambda(\vectx) & M_{\lambda/\mu}(\vectx) \\
 -\trans M_{\lambda/\mu}(\vectx) & O
\end{pmatrix}.
$$
\end{demo}

If $\vectx$ consists of a single variable, then we have the following proposition.

\begin{prop}
\label{prop:skewP-3}
Let $\FF$ be an admissible sequence and $\vectx = (x)$ a single variable.
\begin{enumerate}
\item[(1)]
For two strict partitions $\lambda$ and $\mu$, 
we have $P^\FF_{\lambda/\mu,p}(x) = 0$ unless $l(\lambda) - l(\mu) \le 2$.
\item[(2)]
If $\FF$ is constant-term free, 
then we have $P^\FF_{\lambda/\mu}(x) = 0$ unless $l(\lambda) - l(\mu) \le 1$, and
\begin{equation}
\label{eq:skewP-x}
P^\FF_{\lambda/\mu}(x)
 =
\det \big( 
 R^\FF_{\lambda_i/\mu_j}(x)
\big)_{1 \le i, j \le l(\lambda)}.
\end{equation}
\end{enumerate}
\end{prop}

\begin{demo}{Proof}
Put $l = l(\lambda)$ and $m = l(\mu)$.

(1)
Since $P_{(\lambda_i,\lambda_j)}(x) = 0$ for $\lambda_i > \lambda_j > 0$ by Proposition~\ref{prop:genP}, we have
$$
P_{\lambda/\mu,p}(x)
 =
\begin{cases}
\Pf \begin{pmatrix}
 O_l & M_{\lambda/\mu}(x) \\
 -\trans M_{\lambda/\mu}(x) & O_m 
\end{pmatrix}
 &\text{if $l \equiv p$ and $m \equiv p \bmod 2$,}
\\
\Pf \begin{pmatrix} 
 O_l & M_{\lambda/\mu^0}(x) \\
 -\trans M_{\lambda/\mu^0}(x) & O_{m+1}
\end{pmatrix}
 &\text{if $l \equiv p$ and $m \not\equiv p \bmod 2$,}
\\
\Pf \begin{pmatrix}
 O_l & T_\lambda(x) & M_{\lambda/\mu}(x) \\
 -\trans T_\lambda(x) & 0 & O \\
 -\trans M_{\lambda/\mu}(x) & O & O_m
\end{pmatrix}
 &\text{if $l \not\equiv p$ and $m \equiv p \bmod 2$,}
\\
\Pf \begin{pmatrix}
 O_l & T_\lambda(x) & T_{\lambda/0}(x) & M_{\lambda/\mu}(x) \\
 -\trans T_\lambda(x) & 0 & 1 & O \\
 -\trans T_{\lambda/0}(x) & -1 & 0 & O \\
 -\trans M_{\lambda/\mu}(x) & O & O & O_m
\end{pmatrix}
 &\text{if $l \not\equiv p$ and $m \not\equiv p \bmod 2$,}
\end{cases}
$$
where $T_\lambda(x)$ (resp. $T_{\lambda/0}(x)$) is the column vector 
with $i$th entry $P_{(\lambda_i)}(x)$ (resp. $R_{\lambda_i/0}(x)$).
By using Proposition~\ref{prop:Pf-Laplace}, we see that 
$P_{\lambda/\mu}(x) = 0$ unless
$$
\begin{cases}
 l = m &\text{if $l \equiv p$ and $m \equiv p \bmod 2$,} \\
 l = m+1 &\text{if $l \equiv p$ and $m \not\equiv p \bmod 2$,} \\
 l = m+1 &\text{if $l \not\equiv p$ and $m \equiv p \bmod 2$,} \\
 l = m \text{ or } m+2 &\text{if $l \not\equiv p$ and $m \not\equiv p \bmod 2$.}
\end{cases}
$$
Here we note that $m \le l \le m+2$ if and only if $l \le (l+m+2)/2$ and $m+1 \le (l+m+2)/2$.

(2)
By using Proposition~\ref{prop:skewP-2} (2) and Proposition~\ref{prop:Pf-Laplace}, 
we have
$$
P^\FF_{\lambda/\mu,p}(x)
 =
\begin{cases}
 \Pf \begin{pmatrix}
  O_l & M^\FF_{\lambda/\mu}(x) \\
  -\trans M^\FF_{\lambda/\mu}(x) & O
 \end{pmatrix}
 =
 (-1)^{\binom{l}{2}} \det M_{\lambda/\mu}(x)
 &\text{if $l=m$,}
\\
 \Pf \begin{pmatrix}
  O_l & M^\FF_{\lambda/\mu^0}(\vectx) \\
  -\trans M^\FF_{\lambda/\mu^0}(\vectx) & O
 \end{pmatrix}
 =
 (-1)^{\binom{l}{2}} \det M_{\lambda/\mu^0}(x)
 &\text{if $l=m+1$.}
\end{cases}
$$
By permuting columns of $M_{\lambda/\mu}(x)$ and $M_{\lambda/\mu^0}(x)$, 
we obtain (\ref{eq:skewP-x}).
\end{demo}

\subsection{%
Proof of Theorem~\ref{thm:JP}
}

We give a proof of Theorem~\ref{thm:JP} by using the same idea 
as in the proof of \cite[Theorem~6.1]{Okada} for Schur $Q$-functions.
A key is played by the following proposition, which interpolates 
the Nimmo-type formula (\ref{eq:Nimmo}) ($n=0$ case) and 
the Schur-type formula (\ref{eq:Schur}) ($p=0$ case).

\begin{prop}
\label{prop:NS}
Let $\FF$ be an admissible polynomial sequence and 
$\vectx = (x_1, \dots, x_n)$ and $\vecty = (y_1, \dots, y_p)$ two sequence of indeterminates.
For a non-negative integer sequence $\alpha = (\alpha_1, \dots, \alpha_r)$, 
let $N^\FF_\alpha(\vectx|\vecty)$ be the $r \times p$ matrix defined by
$$
N^\FF_\alpha(\vectx|\vecty)
 =
\big(
 P^\FF_{\alpha_i}(\vectx,y_j)
\big)_{1 \le i \le r, 1 \le j \le p}
$$
Then, for a strict partition $\lambda$ of length $l$, we have
\begin{equation}
\label{eq:NS}
P^\FF_\lambda(\vectx,\vecty)
 =
\begin{cases}
\dfrac{ 1 }{ \Delta(\vecty) }
 \Pf \begin{pmatrix}
  S^\FF_\lambda(\vectx) & N^\FF_\lambda(\vectx|\vecty) \\
  -\trans N^\FF_\lambda(\vectx|\vecty) & A(\vecty)
 \end{pmatrix}
 &\text{if $l+p$ is even,} \\
\dfrac{ 1 }{ \Delta(\vecty) }
 \Pf \begin{pmatrix}
  S^\FF_{\lambda^0}(\vectx) & N^\FF_{\lambda^0}(\vectx|\vecty) \\
  -\trans N^\FF_{\lambda^0}(\vectx|\vecty) & A(\vecty)
 \end{pmatrix}
 &\text{if $l+p$ is odd,}
\end{cases}
\end{equation}
where $\lambda^0 = (\lambda_1, \dots, \lambda_l, 0)$, 
and $S^\FF_\alpha(\vectx)$ and $A(\vecty)$ are given by (\ref{eq:S}) and (\ref{eq:A}) respectively.
\end{prop}

We denote by $P'_\lambda(\vectx|\vecty)$ the right hand side of (\ref{eq:NS}).
First we prove that $P'_\lambda(\vectx|\vecty)$ satisfy the Schur-type Pfaffian formula:

\begin{lemma}
\label{lem:Schur'}
For a strict partition $\lambda$ of length $l$, we have
$$
P'_\lambda (\vectx|\vecty)
 =
\begin{cases}
\Pf \left( P'_{(\lambda_i,\lambda_j)}(\vectx|\vecty) \right)_{1 \le i,j \le l}
 &\text{if $l$ is even,}
\\
\Pf \left( P'_{(\lambda_i,\lambda_j)}(\vectx|\vecty) \right)_{1 \le i,j \le l+1}
 &\text{if $l$ is odd,}
\end{cases}
$$
where $\lambda_{l+1} = 0$ if $l$ is odd, 
and we use the similar convention to (\ref{eq:convention2}).
\end{lemma}

\begin{demo}{Proof}
The proof is similar to that of Theorem~\ref{thm:Schur}, so we omit the detail.
We apply a Pfaffian analogue of Sylevester formula (Proposition~\ref{prop:Pf-Sylvester}) to 
the following matrix $X$ (after permuting rows/columns):
$$
X =
\begin{cases}
\begin{pmatrix}
 S_\lambda(\vectx) & N_\lambda(\vectx|\vecty) \\
 -\trans N_\lambda(\vectx|\vecty) & A(\vecty)
\end{pmatrix}
 &\text{if $l$ is even and $p$ is even,}
\\
\begin{pmatrix}
 S_{\lambda^0}(\vectx) & N_{\lambda^0}(\vectx|\vecty) \\
 -\trans N_{\lambda^0}(\vectx|\vecty) & A(\vecty)
\end{pmatrix}
 &\text{if $l$ is odd and $p$ is even,}
\\
\begin{pmatrix}
 S_\lambda(\vectx) & T_\lambda(\vectx) & N_\lambda(\vectx|\vecty) \\
 -\trans T_\lambda(\vectx) & 0 & \trans\vectone \\
 -\trans N_\lambda(\vectx|\vecty) & -\vectone & A(\vecty)
\end{pmatrix}
 &\text{if $l$ is even and $p$ is even,}
\\
\begin{pmatrix}
 S_\lambda(\vectx) & T_\lambda(\vectx) & \vectzero & N_\lambda(\vectx|\vecty) \\
 -\trans T_\lambda(\vectx) & 0 & 1 & O \\
 0 & -1 & 0 & \trans\vectone \\
 -\trans N_\lambda(\vectx|\vecty) & O & -\trans\vectone & A(\vecty)
\end{pmatrix}
 &\text{if $l$ is odd and $p$ is odd,}
\end{cases}
$$
where $T_\lambda(\vectx)$ is the column vector with $i$th entry $P_{(\lambda_i)}(\vectx)$.
\end{demo}

\begin{demo}{Proof of Theorem~\ref{prop:NS}}
By comparing Theorem~\ref{thm:Schur} with Lemma~\ref{lem:Schur'}, 
the proof of Theorem~\ref{prop:NS} is reduced to showing
$$
P_{(r)}(\vectx,\vecty) = P'_{(r)}(\vectx|\vecty),
\quad
P_{(r,s)}(\vectx,\vecty) = P'_{(r,s)}(\vectx|\vecty).
$$
By considering the generating functions and using Proposition~\ref{prop:GF}, 
it is enough to prove
\begin{equation}
\label{eq:GF'1}
\sum_{r \ge 0} P'_{(r)}(\vectx|\vecty) \hat{f}_r(z)
 =
\begin{cases}
 \Pi_z(\vectx,\vecty)
 &\text{if $n+p$ is odd},
\\
 \Pi_z(\vectx,\vecty) + \hat{f}_0(z) - 1
 &\text{if $n+p$ is even},
\end{cases}
\end{equation}
and
\begin{multline}
\label{eq:GF'2}
\sum_{r, s \ge 0} P'_{(r,s)}(\vectx|\vecty) \hat{f}_r(z) \hat{f}_s(w)
\\
=
\left\{
\begin{array}{l}
 \dfrac{z-w}{z+w} \big(
  \Pi_z(\vectx,\vecty) \Pi_w(\vectx,\vecty) -1
 \big)
 + (\hat{f}_0(w)-1) \Pi_z(\vectx,\vecty)
 - (\hat{f}_0(z)-1) \Pi_w(\vectx,\vecty)
\\
\hfill\text{if $n+p$ is odd},
\\
 \dfrac{z-w}{z+w} \big(
  \Pi_z(\vectx,\vecty) \Pi_w(\vectx,\vecty) -1
 \big)
\hfill\text{if $n+p$ is even}.
\end{array}
\right.
\end{multline}

First we prove (\ref{eq:GF'1}).
We put
$$
F_z(\vectx) = \sum_{r \ge 0} P_{(r)}(\vectx) \hat{f}_r(\vectx),
\quad
F^+_z(\vectx) = \sum_{r \ge 1} P_{(r)}(\vectx) \hat{f}_r(\vectx) = F_z(\vectx) - \hat{f}_0(z).
$$
Let $U_z(\vectx|\vecty)$ (resp. $U^+_z(\vectx|\vecty)$) be the row vector with $j$th entry 
$\sum_{r \ge 0} P_{(r)}(\vectx,y_j)\hat{f}_r(z)$ 
(resp.$\sum_{r \ge 1} P_{(r)}(\vectx,y_j)\hat{f}_r(z)$).
Then by the definition of $P'_{(r)}(\vectx)$ and the multilinearity of Pfaffians, we have
$$
\sum_{r \ge 0} P'_{(r)}(\vectx|\vecty) \hat{f}_r(z)
 =
\begin{cases}
 \dfrac{ 1 }
       { \Delta(\vecty) }
 \Pf \begin{pmatrix}
  0 & U_z(\vectx|\vecty) \\
  -\trans U_z(\vectx|\vecty) & A(\vecty)
 \end{pmatrix}
&\text{if $p$ is odd,}
\\
 \hat{f}_0(z)
  +
 \dfrac{ 1 }
       { \Delta(\vecty) }
 \Pf \begin{pmatrix}
  0 & F_z^+(\vectx) & U_z^+(\vectx|\vecty) \\
  - F^+_z(\vectx) & 0             & \trans\vectone \\
  - \trans U^+z(\vectx|\vecty) & -\vectone & A(\vecty)
 \end{pmatrix}
&\text{if $p$ is even.}
\end{cases}
$$
By adding the 2nd row/column multiplied by $\hat{f}_0(z)$ to the 1st row/column,
we see that
$$
\Pf \begin{pmatrix}
  0 & F_z^+(\vectx) & U_z^+(\vectx|\vecty) \\
  - F^+_z(\vectx) & 0             & \trans\vectone \\
  - \trans U^+z(\vectx|\vecty) & -\vectone & A(\vecty)
\end{pmatrix}
 =
\Pf \begin{pmatrix}
  0 & F_z^+(\vectx) & U_z(\vectx|\vecty) \\
  - F^+_z(\vectx) & 0             & \trans\vectone \\
  - \trans U_z(\vectx|\vecty) & -\vectone & A(\vecty)
\end{pmatrix}.
$$
By Proposition~\ref{prop:GF}, we have
\begin{gather}
\label{eq:GF-U}
U_z(\vectx|\vecty)
 =
\begin{cases}
 \Pi_z(\vectx) \trans B_z(\vecty) &\text{if $n$ is even}, \\
 \Pi_z(\vectx) \trans B_z(\vecty) + (\hat{f}_0(z)-1) \trans\vectone &\text{if $n$ is odd}.
\end{cases}
\\
\label{eq:GF-F}
F^+_z(\vectx)
 =
\begin{cases}
 \Pi_z(\vectx) - 1 &\text{if $n$ is even,} \\
 \Pi_z(\vectx) - \hat{f}_0(z) &\text{if $n$ is odd,}
\end{cases}
\end{gather}
where $B_z(\vectx)$ is the column vector with $i$th entry $(1+x_iz)/(1-x_iz)$.
Now we divide into 4 cases according to the parity of $p$ and $n$.

If $p$ is odd and $n$ is even, then by using (\ref{eq:Schur-Pf3}) we have
$$
\sum_{r \ge 0} P'_{(r)}(\vectx|\vecty) \hat{f}_r(z)
 =
\frac{ 1 }
     { \Delta(\vecty) }
\cdot
\Pi_z(\vectx)
\Pf \begin{pmatrix}
 0 & \trans B_z(\vecty) \\
 -B_z(\vecty) & A(\vecty)
\end{pmatrix}
 =
\Pi_z(\vectx) \Pi_z(\vecty).
$$
If $p$ is odd and $n$ is odd, then by using the multilinearity of Pfaffians and 
(\ref{eq:Schur-Pf3}), (\ref{eq:Schur-Pf2}), we have
\begin{align*}
\sum_{r \ge 0} P'_{(r)}(\vectx|\vecty) \hat{f}_r(z)
 &=
\frac{ 1 }
     { \Delta(\vecty) }
\left[
\Pi_z(\vectx)
\Pf \begin{pmatrix}
 0 & \trans B_z(\vecty) \\
 -B_z(\vecty) & A(\vecty)
\end{pmatrix}
 +
(\hat{f}_0(z)-1)
\Pf \begin{pmatrix}
 0 & \trans\vectone \\
 -\vectone & A(\vecty)
\end{pmatrix}
\right]
\\
&
 =
\Pi_z(\vectx,\vecty) + \hat{f}_0(z) - 1.
\end{align*}
If $p$ is even and $n$ is even, 
then by using the multilinearity and the expansion of Pfaffians, 
(\ref{eq:Schur-Pf3}) (with $w=0$) and (\ref{eq:Schur-Pf2}), we have
\begin{align*}
&
\sum_{r \ge 0} P'_{(r)}(\vectx|\vecty) \hat{f}_r(z)
\\
 &=
\hat{f}_0(z)
\\
&\quad
 +
\frac{ 1 }
     { \Delta(\vecty) }
\left[
 (\Pi_z(\vectx)-1)
\Pf \begin{pmatrix}
 0 & 1 & \trans\vectzero \\
 -1 & 0                & \trans\vectone \\
 0 & -\vectone                & A(\vecty)
\end{pmatrix}
 +
\Pi_z(\vectx)
\Pf \begin{pmatrix}
 0 & 0 & \trans B_z(\vecty) \\
 0 & 0                & \trans\vectone \\
 -B_z(\vecty) & -\vectone & A(\vecty)
\end{pmatrix}
\right]
\\
 &=
\Pi_z(\vectx,\vecty) + \hat{f}_0(z) - 1.
\end{align*}
If $p$ is even and $n$ is odd, 
then similarly we have
\begin{align*}
&
\sum_{r \ge 0} P'_{(r)}(\vectx|\vecty) \hat{f}_r(z)
\\
 &=
\hat{f}_0(z)
\\
 &\quad
 +
\frac{ 1 }
     { \Delta(\vecty) }
\left[
 (\Pi_z(\vectx) - \hat{f}_0(z))
\Pf \begin{pmatrix}
 0 & 1 & 0 \\
 -1 & 0                      & \trans\vectone \\
 0 & -\vectone                      & A(\vecty)
\end{pmatrix}
 +
\Pi_z(\vectx)
\Pf \begin{pmatrix}
 0 & 0 & \trans B_z(\vecty) \\
 0 & 0 & \trans\vectone \\
 - B_z(\vecty) & -\vectone & A(\vecty)
\end{pmatrix}
\right]
\\
 &=
\Pi_z(\vectx,\vecty).
\end{align*}

Next we prove (\ref{eq:GF'2}).
We put
$$
G_{z,w}(\vectx)
 =
\sum_{r,s \ge 0} P_{(r,s)}(\vectx) \hat{f}_r(z) \hat{f}_s(w),
\quad
G^{++}_{z,w}(\vectx)
 =
\sum_{r,s \ge 1} P_{(r,s)}(\vectx) \hat{f}_r(z) \hat{f}_s(w).
$$
By the definition of $P'_{(r,s)}(\vectx)$ and the multilinearity of Pfaffians, 
we see that, if $p$ is even, then
$$
\sum_{r,s \ge 0} P'_{(r,s)}(\vectx|\vecty) \hat{f}_r(z) \hat{f}_s(w)
 =
\frac{ 1 }
     { \Delta(\vecty) }
\Pf \begin{pmatrix}
 0 & G_{z,w}(\vectx) & U_z(\vectx|\vecty) \\
 -G_{z,w}(\vectx) & 0               & U_w(\vectx|\vecty) \\
 -\trans U_z(\vectx|\vecty) & -\trans U_w(\vectx|\vecty) & A(\vecty)
\end{pmatrix},
$$
and, if $p$ is odd, then
\begin{align*}
&
\sum_{r,s \ge 0} P'_{(r,s)}(\vectx|\vecty) \hat{f}_r(z) \hat{f}_s(w)
\\
&=
\sum_{r,s \ge 1} P'_{(r,s)}(\vectx|\vecty) \hat{f}_r(z) \hat{f}_s(w)
 +
\sum_{r \ge 1} P'_{(r,0)}(\vectx|\vecty) \hat{f}_r(z) \hat{f}_0(w)
 +
\sum_{s \ge 1} P'_{(0,s)}(\vectx|\vecty) \hat{f}_0(z) \hat{f}_s(w)
\\
&=
\frac{ 1 }
     { \Delta(\vecty) }
\Pf \begin{pmatrix}
 0 & G_{z,w}^{++}(\vectx) & F_z^+(\vectx) & U_z^+(\vectx|\vecty) \\
 -G^{++}_{z,w}(\vectx) & 0                     & F_w^+(\vectx) & U_w^+(\vectx|\vecty) \\
 -F^+_z(\vectx) & -F^+_w(\vectx) & 0           & \trans\vectone \\
 -\trans U^+_z(\vectx|\vecty) & -\trans U^+_w(\vectx|\vecty) & -\vectone & A(\vecty)
\end{pmatrix}
\\
&\quad
+
\hat{f}_0(w)
\cdot
\frac{ 1 }
     { \Delta(\vecty) }
\Pf \begin{pmatrix}
 0 & U_z^+(\vectx|\vecty) \\
 -\trans U^+_z(\vectx|\vecty) & A(\vecty)
\end{pmatrix}
\\
&\quad
-
\hat{f}_0(z)
\cdot
\frac{ 1 }
     { \Delta(\vecty) }
\Pf \begin{pmatrix}
 0 & U_w^+(\vectx|\vecty) \\
 -\trans U^+_w(\vectx|\vecty) & A(\vecty)
\end{pmatrix}.
\end{align*}
Here we note
$$
G_{z,w}(\vectx)
 =
G^{++}_{z,w}(\vectx) + \hat{f}_0(w) F^+_z(\vectx) - \hat{f}_0(z) F^+_w(\vectx),
$$
and
$$
U_z(\vectx|\vecty) = U^+_z(\vectx|\vecty) + \hat{f}_0(z) \trans\vectone,
\quad
U_w(\vectx|\vecty) = U^+_w(\vectx|\vecty) + \hat{f}_0(w) \trans\vectone.
$$
By adding the 3rd row/column multiplied by $\hat{f}_0(z)$ (resp. $\hat{f}_0(w)$) 
to the 1st (resp. 2nd) row/column in the first Pfaffian, we see that
\begin{multline*}
\Pf \begin{pmatrix}
 0 & G_{z,w}^{++}(\vectx) & F_z^+(\vectx) & U_z^+(\vectx|\vecty) \\
 -G^{++}_{z,w}(\vectx) & 0                     & F_w^+(\vectx) & U_w^+(\vectx|\vecty) \\
 -F^+_z(\vectx) & -F^+_w(\vectx) & 0           & \trans\vectone \\
 -\trans U^+_z(\vectx|\vecty) & -\trans U^+_w(\vectx|\vecty) & -\vectone & A(\vecty)
\end{pmatrix}
\\
 =
\Pf \begin{pmatrix}
 0 & G_{z,w}(\vectx) & F_z^+(\vectx) & U_z(\vectx|\vecty) \\
 -G_{z,w}(\vectx) & 0                     & F_w^+(\vectx) & U_w(\vectx|\vecty) \\
 -F^+_z(\vectx) & -F^+_w(\vectx) & 0           & \trans\vectone \\
 -\trans U_z(\vectx|\vecty) & -\trans U_w(\vectx|\vecty) & -\vectone & A(\vecty)
\end{pmatrix}.
\end{multline*}
By Proposition~\ref{prop:GF}, we have
\begin{multline}
\label{eq:GF-G}
G_{z,w}(\vectx)
\\
 =
\begin{cases}
 \dfrac{z-w}{z+w} \big(
  \Pi_z(\vectx) \Pi_w(\vectx) -1
 \big)
 &\text{if $n$ is even,}
\\
\dfrac{z-w}{z+w} \big(
  \Pi_z(\vectx) \Pi_w(\vectx) -1
 \big)
 + (\hat{f}_0(w)-1) \Pi_z(\vectx)
 - (\hat{f}_0(z)-1) \Pi_w(\vectx)
 &\text{if $n$ is odd.}
\end{cases}
\end{multline}
Now we divide into 4 cases according to the parity of $p$ and $n$, 
and in each case we evaluate Pfaffians by using (\ref{eq:GF-U}), (\ref{eq:GF-F}) and (\ref{eq:GF-G}) 
together with Propositions~\ref{prop:Schur-Pf} and \ref{prop:Schur-Pf-var}.
The rest of the proof is done by straightforward computation, 
so we omit it.
\end{demo}

Now we are in position to give a proof of Theorem~\ref{thm:JP}.

\begin{demo}{Proof of Theorem~\ref{thm:JP}}
We put $l=l(\lambda)$.

First we consider the case where $l \equiv p \bmod 2$.
Then we apply a Pfaffian analogue of Cauchy-Binet formula (\ref{eq:Pf-CB1}) to 
the matrices given by
\begin{gather*}
A = S_\lambda(\vectx),
\quad
S = \big( R_{\lambda_i/k}(\vectx) \big)_{1 \le i \le l, k \ge 0},
\\
B = A(\vecty),
\quad
T = \big( f_k(x_i) \big)_{1 \le i \le p, k \ge 0}.
\end{gather*}
Strict partitions $\mu$ are in bijection with subsets 
of $\Nat$ satisfying $\# I = p \bmod 2$ 
via the correspondence $\mu \mapsto I = \{ \mu_1, \dots, \mu_{l(\mu)} \}$ or $\{ \mu_1, \dots, \mu_{l(\mu)}, 0 \}$,
and we see that
\begin{align*}
\Pf \begin{pmatrix}
 S_\lambda(\vectx) & S([l];I) \\
 -\trans S([l];I) & O
\end{pmatrix}
 &=
P_{\lambda/\mu,p}(\vectx),
\\
\Pf \begin{pmatrix}
 A(\vecty) & T([p];I) \\
 -\trans T([p];I) & O
\end{pmatrix}
 &=
(-1)^{\binom{\# I}{2}} \Delta(\vecty) P_\mu(\vecty).
\end{align*}
By (\ref{eq:skewP1}), the $(i,j)$ entry of $S \trans T$ is equal to
$$
\sum_{k \ge 0} R_{\lambda_i/k}(\vectx) f_k(y_j)
 =
P_{(\lambda_i)}(\vectx,y_j).
$$
Hence, by applying (\ref{eq:Pf-CB1}), we have
$$
\sum_\mu
P_{\lambda/\mu,p}(\vectx) P_\mu(\vecty)
 =
\frac{ 1 }
     { \Delta(\vecty) }
\Pf \begin{pmatrix}
 S_\lambda(\vectx) & N_\lambda(\vectx,\vecty) \\
 -\trans N_\lambda(\vectx,\vecty) & A(\vecty)
\end{pmatrix}.
$$

If $l \not\equiv p \bmod 2$, then we apply (\ref{eq:Pf-CB1}) to the matrices given by
\begin{gather*}
A = S_{\lambda^0}(\vectx),
\quad
S = \big( R_{\lambda_i/k}(\vectx) \big)_{1 \le i \le l+1, k \ge 0},
\\
B = A(\vecty),
\quad
T = \big( f_k(y_i) \big)_{1 \le i \le p, k \ge 0}.
\end{gather*}
Then by an argument similar to above, we obtain
$$
\sum_\mu
P_{\lambda/\mu,p}(\vectx) P_\mu(\vecty)
 =
\frac{ 1 }
     { \Delta(\vecty) }
\Pf \begin{pmatrix}
 S_{\lambda^0}(\vectx) & N_{\lambda^0}(\vectx,\vecty) \\
 -\trans N_{\lambda^0}(\vectx,\vecty) & A(\vecty)
\end{pmatrix}.
$$
Now the proof of Theorem~\ref{thm:JP} can be completed by using Proposition~\ref{prop:NS}.
\end{demo}

\section{%
Pieri-type rule
}

In this section we give a Pieri-type rule for the product of any generalized $P$-function $P^\FF_\lambda(\vectx)$ 
with Schur $Q$-function $Q_{(r)}(\vectx)$ corresponding to a one-row partition.

\subsection{%
The ring of Schur $\boldsymbol{P}$- and $\boldsymbol{Q}$-functions
}

Let $\Gamma^{(n)}$ be the subring of the ring of symmetric polynomials $\Lambda^{(n)} = K[x_1, \dots, x_n]^{S_n}$ 
defined by
$$
\Gamma^{(n)}
 = 
\{ 
 f \in K[x_1, \dots, x_n]^{S_n}
 :
 \text{$f(t,-t,x_3, \dots, x_n)$ is independent of $t$}
\}
$$
Then it is known that Schur $P$-functions $\{ P_\lambda(\vectx) : \lambda \in \Strict^{(n)} \}$ 
form a basis of $\Gamma^{(n)}$, 
where $\Strict^{(n)}$ is the set of all strict partitions of length $\le n$ 
(see \cite[Theorem~2.11]{Pragacz}).

We give a relation between two families of generalized $P$-functions 
associated to different admissible sequences.

\begin{prop}
\label{prop:PF_by_PG}
Let $\FF = \{ f_d \}_{d=0}^\infty$ and $\GG = \{ g_d \}_{d=0}^\infty$ be two admissible sequences.
For a strict partition $\lambda$ of length $l \le n$, the generalized $P$-function 
$P^{\FF}_\lambda(\vectx)$ associated to $\FF$ can be written as a $K$-linear combination 
of the generalized $P$-functions $P^{\GG}_\mu(\vectx)$ associated to $\GG$ 
in the following form:
$$
P^{\FF}_\lambda(\vectx)
 =
a_{\lambda,\lambda} P^{\GG}_\lambda(\vectx)
 +
\sum_{\mu \subsetneq \lambda} a_{\lambda,\mu} P^{\GG}_\mu(\vectx),
$$
where $a_{\lambda,\lambda} \neq 0$, and $\mu$ runs over all strict partitions 
satisfying $\mu \subsetneq \lambda$.
\end{prop}

\begin{demo}{Proof}
We write $f_k = \sum_{i=0}^\infty a_{k,l} g_l$ for $k \ge 0$.
Then by the assumption (\ref{eq:cond}) 
we have $a_{k,l} = 0$ for $k < l$, $a_{k,k} \neq 0$, and $a_{0,0} = 1$.

If $n+l$ is even, then by using the multilinear and alternating property of Pfaffians we have
\begin{align*}
P^{\FF}_\lambda(\vectx)
 &=
\frac{ 1 }
     { \Delta(\vectx) }
\sum_{\alpha \in \Nat^l}
 a_{\lambda_1,\alpha_1} \dots a_{\lambda_l,\alpha_l}
\Pf \begin{pmatrix}
 A(\vectx) & V^\GG_\alpha(\vectx) \\
 -\trans V^\GG_\alpha(\vectx) & O
\end{pmatrix}
\\
 &=
\sum_{\mu_1 > \dots > \mu_l \ge 0}
 \det \big( a_{\lambda_i,\mu_j} \big)_{1 \le i, j \le l}
 P^\GG_\mu(\vectx),
\end{align*}
where $\mu$ runs over all strict partitions of length $l-1$ or $l$.
Similarly, if $n+l$ is odd, then we have
$$
P^{\FF}_\lambda(\vectx)
 =
\frac{ 1 }
     { \Delta(\vectx) }
\sum_{\alpha \in \Nat^{l+1}}
 a_{\lambda_1,\alpha_1} \dots a_{\lambda_l,\alpha_l} a_{0,\alpha_{l+1}}
\Pf \begin{pmatrix}
 A(\vectx) & V^\GG_\alpha(\vectx) \\
 -\trans V^\GG_\alpha(\vectx) & O
\end{pmatrix}.
$$
Since $a_{0,l} = 0$ for $l > 0$, we see that
$$
P^{\FF}_\lambda(\vectx)
 =
\sum_{\mu_1 > \dots > \mu_l > 0}
 \det \big( a_{\lambda_i,\mu_j} \big)_{1 \le i, j \le l}
 P^\GG_\mu(\vectx),
$$
where $\mu$ runs over all strict partitions of length $l$.
For a strict partition $\mu$ of length $l-1$ or $l$, we put
$$
a_{\lambda,\mu} = \det \big( a_{\lambda_i,\mu_j} \big)_{1 \le i, j \le l},
$$
where $\mu_l = 0$ if $l(\mu) = l-1$.
We prove that $a_{\lambda,\lambda} \neq 0$ 
and $a_{\lambda,\mu} = 0$ unless $\lambda \supset \mu$.
Since the matrix $\big( a_{\lambda_i, \lambda_j} \big)_{1 \le i, j \le l}$ is upper-triangular, 
we have $a_{\lambda,\lambda} = \prod_{i=1}^l a_{\lambda_i,\lambda_i} \neq 0$.
If there is an index $k$ such that $\lambda_k < \mu_k$, 
then we have $\lambda_i \le \lambda_k < \mu_k \le \mu_j$ and $a_{\lambda_i,\mu_j} = 0$ 
for $i \ge k$ and $j \le k$ and thus $a_{\lambda,\mu} = 0$.
Hence we obtain the desired result.
\end{demo}

\begin{corollary}
\label{cor:basis}
The generalized $P$-functions $\{ P^{\FF}_\lambda(\vectx) : \lambda \in \Strict^{(n)} \}$ 
associated to a fixed sequence $\FF$ form a basis of $\Gamma^{(n)}$.
\end{corollary}

\subsection{%
Pieri-type rule
}

Let $q_r(\vectx) = Q_{(r)}(\vectx)$ be the Schur $Q$-function corresponding to a one-row partition $(r)$,  
and set $q_0(\vectx) = 1$.
It is known (see \cite[III.(8.5)]{Macdonald95}) that $\Gamma^{(n)}$ is generated by $q_r(\vectx)$ ($r \ge 1$).
Thus the algebra structure of $\Gamma^{(n)}$ is governed by the multiplication rule for $q_r(\vectx)$s.

\begin{theorem}
\label{thm:Pieri}
Let $\FF = \{ f_d \}_{d=0}^\infty$ be an admissible sequence.
We define formal power series $b^s_r(z)$ by the relation
\begin{equation}
\label{eq:Pieri-coeff1}
f_r(u) \cdot \frac{ 1 + u z }{ 1 - u z }
 =
\sum_{s=0}^\infty b^s_r(z) f_s(u).
\end{equation}
And we define the modified Pieri coefficients $c^\lambda_{\mu,r}$ by the relation
\begin{equation}
\label{eq:Pieri-coeff2}
P^\FF_\mu(\vectx) \cdot q_r(\vectx)
 =
\sum_\lambda c^\lambda_{\mu,r} P^\FF_\lambda(\vectx),
\end{equation}
where the summation is taken over all strict partitions.
Then the generating function of modified Pieri coefficients is given as follows:
\begin{equation}
\label{eq:Pieri}
\sum_{r=0}^\infty c^\lambda_{\mu,r} z^r
 =
\begin{cases}
 \det B^\lambda_\mu &\text{if $n+l(\mu)$ is even and $l(\lambda) = l(\mu)$,} \\
 \det B^{\lambda^0}_\mu &\text{if $n+l(\mu)$ is even and $l(\lambda) = l(\mu)-1$,} \\
 \det B^\lambda_{\mu^0} &\text{if $n+l(\mu)$ is even and $l(\lambda) = l(\mu)+1$,} \\
 \det B^{\lambda^0}_{\mu^0} &\text{if $n+l(\mu)$ is odd and $l(\lambda) = l(\mu)$,} \\
 \det B^\lambda_{\mu^0} &\text{if $n+l(\mu)$ is odd and $l(\lambda) = l(\mu)+1$,} \\
 0 &\text{otherwise,}
\end{cases}
\end{equation}
where $B^\alpha_\beta = \big( b^{\alpha_i}_{\beta_j}(z) \big)_{1 \le i, j \le r}$ 
for $\alpha = (\alpha_1, \dots, \alpha_r)$ and $\beta = (\beta_1, \dots, \beta_r)$.
\end{theorem}

If we put
$$
c^\lambda_\mu(z) = \sum_{r=0}^\infty c^\lambda_{\mu,r} z^r,
$$
then it follows from (\ref{eq:GF-SchurQ1}) and (\ref{eq:Pieri-coeff2}) that
\begin{equation}
\label{eq:Pieri-coeff3}
P^\FF_\mu(z) \cdot \Pi_z(\vectx)
 =
\sum_\lambda c^\lambda_\mu(z) P^\FF_\lambda(z).
\end{equation}
In order to prove the above theorem, 
we derive a Nimmo-type expression of the product $P^\FF_\mu(\vectx) \cdot \Pi_z(\vectx)$ 
in terms of Pfaffian.

\begin{lemma}
\label{lem:Pieri-Nimmo}
For a sequence $\alpha = (\alpha_1, \dots, \alpha_r)$ of nonnegative integers, we put
$$
W_\alpha(\vectx)
 =
\left(
 f_{\alpha_j}(x_i) \cdot \dfrac{1+x_iz}{1-x_iz}
\right)_{1 \le i \le n, 1 \le j \le r}.
$$
Then, for a strict partition $\mu$ of length $m$, we have
$$
P^\FF_\mu(\vectx) \cdot \Pi_z(\vectx)
 =
\begin{cases}
\dfrac{1}{\Delta(\vectx)}
 \Pf \begin{pmatrix}
  A(\vectx) & W_\mu(\vectx) & W_{(0)}(\vectx) & \vectone \\
  -\trans W_\mu(\vectx) & 0 & 0 & 0 \\
  -\trans W_{(0)}(\vectx) & 0 & 0 & 1 \\
  -\trans\vectone & 0 & -1 & 0
 \end{pmatrix}
 &\text{if $n+m$ is even,} \\
\dfrac{1}{\Delta(\vectx)}
 \Pf \begin{pmatrix}
  A(\vectx) & W_{\mu^0}(\vectx) \\
  -\trans W_{\mu^0}(\vectx) & 0
 \end{pmatrix}
 &\text{if $n+m$ is odd.} \\
\end{cases}
$$
\end{lemma}

\begin{demo}{Proof}
The method of the proof is the same as in the proof of Theorem~\ref{thm:N=HL}, 
so we will give a sketch of the proof.

Since $\Pi_z(\vectx) = \prod_{i=1}^n (1+x_iz)/(1-x_iz)$ is invariant under the symmetric group $S_n$, 
it follows from Corollary~\ref{cor:N=HL} that
\begin{align*}
&
P_\mu(\vectx) \cdot \Pi_z(\vectx)
\\
 &=
\sum_{w \in S_n/S_{n-m}}
 w \left(
  \prod_{i=1}^m f_{\mu_i}(x_i)
  \prod_{i=1}^n \frac{1+x_iz}{1-x_iz}
  \prod_{\substack{1 \le i < j \le n \\ i \le m}} \frac{x_i+x_j}{x_i-x_j}
 \right)
\\
&=
\frac{ (-1)^{\binom{n}{2} + \binom{n-m}{2}} }
     { \Delta_n(\vectx) }
\\
&\quad
\times
\sum_{w \in S_n/S_{n-m}}
 \sgn(w) w \left(
  \prod_{i=1}^m \left( f_{\mu_i}(x_i) \frac{1+x_iz}{1-x_iz} \right)
  \prod_{i=m+1}^n \frac{1+x_iz}{1-x_iz}
  \prod_{m+1 \le i < j \le n} \frac{x_j-x_i}{x_j+x_i}
\right).
\end{align*}
By using (\ref{eq:Schur-Pf3}) with $p=2$, $y_1=z$, $y_2 =0$, or $p=1$, $y_1=z$, 
we obtain
$$
P_\mu(\vectx) \cdot \Pi_z(\vectx)
 =
\dfrac{ (-1)^{\binom{n}{2} + \binom{n-m}{2}} }
      { \Delta(\vectx) }
\sum_{w \in S_n/S_m \times S_{n-m}}
 \sgn(w) w \Big(
  \det W_\mu(\vectx_{[m]})
  \Pf X
 \Big),
$$
where the skew-symmetric matrix $X$ is given by
$$
X
 = 
\begin{cases}
\begin{pmatrix}
 A(\vectx_{[n] \setminus [m]}) & W_{(0)}(\vectx_{[n] \setminus [m]}) & \vectone \\
 -\trans W_{(0)}(\vectx_{[n] \setminus [m]}) & 0                      & 1 \\
 -\trans\vectone         & -1                      & 0
\end{pmatrix}
&\text{if $n+m$ is even,}
\\
\begin{pmatrix}
 A(\vectx_{[n] \setminus [m]}) & W_{(0)}(\vectx_{[n] \setminus [m]}) \\
 -\trans W_{(0)}(\vectx_{[n] \setminus [m]}) & 0                      
\end{pmatrix}
&\text{if $n+m$ is odd.}
\end{cases}
$$
Now we can use a Pfaffian analogue of Laplace expansion (Proposition~\ref{prop:Pf-Laplace}) 
to complete the proof.
\end{demo}

\begin{demo}{Proof of Theorem~\ref{thm:Pieri}}
The argument is similar to that in the proof of Proposition~\ref{prop:PF_by_PG}.

First we consider the case where $n+m$ is even.
In this case, by using the multilinearity and the expansion along the last columns/row, 
we have
\begin{multline*}
\Pf \begin{pmatrix}
 A(\vectx) & W_\mu(\vectx) & W_{(0)}(\vectx) & \vectone \\
 -\trans W_\mu(\vectx)    & O & O & O \\
 -\trans W_{(0)}(\vectx)  & O & 0 & 1 \\
 -\trans\vectone   & O & -1 & 0
\end{pmatrix}
\\
 =
\Pf \begin{pmatrix}
 A(\vectx) & W_{\mu^0}(\vectx) & \vectone \\
 -\trans W_{\mu^0}(\vectx)   & O & O \\
 -\trans\vectone             & O & 0
\end{pmatrix}
+
(-1)^{n+m}
\Pf \begin{pmatrix}
 A(\vectx) & W_\mu(\vectx) \\
 -\trans W_\mu(\vectx) & 0
\end{pmatrix}.
\end{multline*}
Since $W_{(r)}(\vectx) = \sum_{s=0}^\infty b^s_r(z) V_{(s)}(\vectx)$ by (\ref{eq:Pieri-coeff1}), 
we can use the multilinear and alternating property to obtain
\begin{align*}
&
\Pf \begin{pmatrix}
 A(\vectx) & W_{\mu^0}(\vectx) & \vectone \\
 -\trans W_{\mu^0}(\vectx)  & O & O \\
 -\trans\vectone   & O & 0
\end{pmatrix}
\\
 &\quad
=
\sum_{\alpha \in \Nat^{m+1}}
 \prod_{i=1}^{m+1} b^{\alpha_i}_{\mu_i}(z)
\Pf \begin{pmatrix}
 A(\vectx) & V_\alpha(\vectx) & \vectone \\
 -\trans V_\alpha(\vectx) & O & O \\
 -\trans\vectone & O & 0
\end{pmatrix}
\\
 &\quad
=
\sum_{\lambda_1 > \dots > \lambda_{m+1} \ge 0}
 \det \big( b^{\lambda_i}_{\mu_j}(z) \big)_{1 \le i, j \le m+1}
\Pf \begin{pmatrix}
 A(\vectx) & V_\lambda(\vectx) & \vectone \\
 -\trans V_\lambda(\vectx)   & O & O \\
 -\trans\vectone         & O & 0
\end{pmatrix}.
\end{align*}
If $\lambda_{m+1} = 0$, then the last column of $V_\lambda(\vectx)$ coincides with $\vectone$, 
so the corresponding Pfaffian vanishes.
Hence we have
$$
\frac{ 1 }
     { \Delta(\vectx) }
\Pf \begin{pmatrix}
 A(\vectx) & W_{\mu^0}(\vectx) & \vectone \\
 -\trans W_{\mu^0}(\vectx)  & O & O \\
 -\trans\vectone   & O & 0
\end{pmatrix}
 =
\sum_\lambda
 \det B^\lambda_{\mu^0}
 P_\lambda(\vectx),
$$
where $\lambda$ runs over all strict partitions of length $m+1$.
Similarly we have
$$
\frac{ 1 }
     { \Delta(\vectx) }
\Pf \begin{pmatrix}
 A(\vectx) & W_\mu(\vectx) \\
 -\trans W_\mu(\vectx) & 0             
\end{pmatrix}
 =
\sum_\lambda
 \det B^{\lambda^*}_\mu
 P_\lambda(\vectx),
$$
where $\lambda$ runs over all strict partitions of length $m-1$ or $m$, 
and $\lambda^* = \lambda^0$ or $\lambda$.

The case where $n+m$ is odd can be treated in a similar manner.
\end{demo}

If $\FF$ is constant-term free, then we have a simpler formula.

\begin{corollary}
\label{cor:Pieri}
If $\FF$ is constant-term free, then we have
$$
c^\lambda_\mu(z)
 =
\begin{cases}
 \det B^\lambda_\mu &\text{if $l(\lambda) = l(\mu)$,} \\
 \det B^\lambda_{\mu^0} &\text{if $l(\lambda) = l(\mu)+1$,} \\
 0 &\text{otherwise,}
\end{cases}
$$
under the same notation as in Theorem~\ref{thm:Pieri}.
\end{corollary}

\begin{demo}{Proof}
By substituting $t=0$ in (\ref{eq:Pieri-coeff1}) and using the assumption $f_d(0) = \delta_{d,0}$, 
we see that $b^0_r(z) = \delta_{r,0}$.
Hence we have
$$
\det B^{\lambda^0}_{\mu^0} = \det B^\lambda_\mu,
\quad
\det B^{\lambda^0}_\mu = 0,
$$
and we obtain Corollary~\ref{cor:Pieri}.
\end{demo}

\section{%
Applications to Factorial $\boldsymbol{P}$- and $\boldsymbol{Q}$-functions
}

In this section we focus on Ivanov's factorial $P$- and $Q$-functions.

\subsection{%
Factorial $\boldsymbol{P}$- and $\boldsymbol{Q}$-functions
}

Recall the definition of Ivanov's factorial $P$- and $Q$-functions.
Let $\vecta = (a_0, a_1, \dots)$ be parameters, called \emph{factorial parameters}.
We define the factorial monomial $(u|\vecta)^d$ by putting
$$
(u|\vecta)^d = \prod_{i=0}^{d-1} (u - a_i).
$$
Then the \emph{factorial $P$-function} $P_\lambda(\vectx|\vecta)$ is 
defined to be the generalized $P$-function $P^\FF_\lambda(\vectx)$ 
associated to $\FF = \{ (u|\vecta)^d \}_{d=0}^\infty$ (see Definition~\ref{def:genP}).
And the \emph{factorial $Q$-function} $Q_\lambda(\vectx|\vecta)$ is defined by
$$
Q_\lambda(\vectx|\vecta) = 2^{l(\lambda)} P_\lambda(\vectx|\vecta).
$$
The factorial $Q$-function $Q_\lambda(\vectx|\vecta)$ is also 
the generalized $P$-function $P^{\FF'}_\lambda(\vectx)$ 
associated to the sequence $\FF' = \{ f'_d \}_{d=0}^\infty$ given by
$$
f'_d(u)
 = 
\begin{cases}
 1 &\text{if $d=0$,} \\
 2 (u|\vecta)^d &\text{if $d \ge 1$.}
\end{cases}
$$
Note that the sequence $\FF = \{ (u|\vecta)^d \}_{d=0}^\infty$ is constant-term free 
if and only if $a_0 = 0$.

Since the factorial $P$- and $Q$-functions are the special cases of our generalized 
$P$-functions, we can recover some formulas in \cite{Ivanov1} and \cite{Ivanov2} 
from the results of Section~2.
For example, we recover \cite[Theorem~9.1]{Ivanov2} without assumption $a_0 = 0$:
$$
Q_\lambda(\vectx|\vecta)
 =
\Pf \left(
 Q_{(\lambda_i,\lambda_j)}(\vectx|\vecta)
\right)_{1 \le i, j \le r},
$$
where $r = l(\lambda)$ if $l(\lambda)$ is even and $l(\lambda)+1$ if $l(\lambda)$ is odd, 
and $\lambda_{l(\lambda)+1} = 0$. 
(We use the convention (\ref{eq:convention2}).)

Next we compute explicitly the dual of $\FF = \{ (t|\vecta)^d \}_{d=0}^\infty$ introduced in Section~3.

\begin{lemma}
\label{lem:fac-dual}
Let $\hat{\FF} = \{ \hat{f}_d \}_{d=0}^\infty$ be the dual of $\FF = \{ (u|\vecta)^d \}_{d=0}^\infty$.
Then we have 
$$
\hat{f}_d(v)
 =
\begin{cases}
 \dfrac{ 1 + a_0 v }{ 1 - a_0 v } &\text{if $d = 0$,} \\
 \dfrac{ 2 v^d }{\prod_{i=0}^d (1 - a_i v) } &\text{if $d \ge 1$.}
\end{cases}
$$
\end{lemma}

\begin{demo}{Proof}
Put $f_d(u) = (u|\vecta)^d$.
The sequence $\{ \hat{f}_d \}_{d=0}^\infty$ is uniquely determined by the relation (\ref{eq:dual}).
Since $f_0 = 1$, we see that $\hat{f}_0(v)$ is determined by substituting $u = a_0$ 
in (\ref{eq:dual}), and obtain $\hat{f}_0 (v) = (1 + a_0 v)/(1 - a_0 v)$.
Let $r > 0$.
Since $f_k(a_r) = 0$ for $k > r$ and $f_r(a_r) \neq 0$, 
we see that $\hat{f}_r(t)$ is determined inductively by the relation
$$
\sum_{k=0}^r \hat{f}_k(v) f_k(a_r)
 =
\frac{ 1 + a_r v }
     { 1 - a_r v }.
$$
Hence it is enough to show 
\begin{equation}
\label{eq:fac-dual1}
\frac{ 1 + a_r v }
     { 1 - a_r v }
 =
\frac{ 1 + a_0 v }
     { 1 - a_0 v }
 +
\sum_{k=1}^r
 \frac{ 2 v^r }
      { \prod_{i=0}^k (1 - a_i v) }
 \prod_{j=0}^{k-1} (a_r - a_j).
\end{equation}
By using
$$
\frac{ 1 + a_r v }{ 1 - a_r v }
-
\frac{ 1 - a_0 v }{ 1 - a_0 v }
 =
\frac{ 2 v (a_r - a_0) }
     { (1 - a_0 v)(1 - a_r v) },
$$
and cancelling the common factor $2 z ( a_r - a_0)/(1 - a_0 z)$, 
we see that (\ref{eq:fac-dual1}) is equivalent to
\begin{equation}
\label{eq:fac-dual2}
\frac{ 1 }{ 1 - a_r v }
 =
\sum_{k=1}^r
 \frac{ v^{k-1} }
      { \prod_{i=1}^k (1 - a_i v) }
 \prod_{j=1}^{k-1} (a_r - a_j).
\end{equation}

We proceed by induction on $r$ to prove (\ref{eq:fac-dual2}).
The case $r = 1$ is trivial.
If $r>1$, then by the induction hypothesis with factorial parameters $(a_2, \dots, a_r)$, 
we have
$$
\frac{ 1 }{1 - a_r v}
 =
\sum_{k=2}^r
 \frac{ v^{k-2} }
      { \prod_{i=2}^k (1 - a_i v) }
\prod_{j=2}^{k-1} (a_r -a_j).
$$
Hence we have
\begin{align*}
\sum_{k=1}^r
 \frac{ v^{k-1} }
      { \prod_{i=1}^k (1 - a_i v) }
 \prod_{j=1}^{k-1} (a_r -a_i)
&=
\frac{ 1 }{ 1 - a_1 v }
 +
\frac{ v ( a_r - a_1 ) }
     { 1 - a_1 v }
\sum_{k=2}^r
 \frac{ v^{k-2} }
      { \prod_{i=2}^k (1 - a_i v) }
 \prod_{j=2}^{k-1} (a_r -a_j)
\\
&=
\frac{ 1 }{ 1 - a_1 v }
 +
\frac{ v ( a_r - a_1 ) }
     { 1 - a_1 v }
\cdot
\frac{ 1 }{ 1 - a_r v }
\\
&=
\frac{ 1 }
     { 1 - a_r v}.
\end{align*}
This completes the proof.
\end{demo}

We denote by $\hat{P}_\lambda(\vectx|\vecta)$ the dual $P$-function $\hat{P}^\FF_\lambda(\vectx)$ associated 
to $\FF = \{ (u|\vecta)^d \}_{d=0}^\infty$.
If the first parameter $a_0$ is equal to $0$, 
then we can recover Korotkikh's dual $P$-function 
(see \cite[Definition~2]{Korotkikh} and (\ref{eq:dualHL})) given by
$$
\hat{P}_\lambda(\vectx|\vecta)
 =
\frac{ 2^l }
     { (n-l)! }
\sum_{w \in S_n}
 w \left(
  \prod_{i=1}^l
   \frac{ x_i^{\lambda_i} }
        { \prod_{k=0}^{\lambda_i} (1 - a_k x_i) }
  \prod_{\substack{1 \le i < j \le n \\ i \le l}} \frac{x_i+x_j}{x_i-x_j}
 \right)
$$
and the Cauchy-type identity (see \cite[Theorem~8]{Korotkikh}):
$$
\sum_\lambda P_\lambda(\vectx|\vecta) \hat{P}_\lambda(\vectx|\vecta)
 =
\prod_{i,j=1}^n
 \frac{ 1 + x_i y_j }
      { 1 - x_i y_j }.
$$

\subsection{%
Factorial skew $\boldsymbol{P}$-functions
}

For two strict partitions $\lambda$ and $\mu$ and a positive integer $p$, 
we denote by $P_{\lambda/\mu,p}(\vectx|\vecta)$ the generalized skew $P$-functions 
$P^\FF_{\lambda/\mu,p}(\vectx)$ associated to $\FF = \{ (u|\vecta)^d \}_{d=0}^\infty$, 
and call it the \emph{factorial skew $P$-function} (see Definition~\ref{def:skewP}).
Since $(u|\vecta)^r$ depends only on $a_0, a_1, \dots, a_{r-1}$, 
it follows from the definition (\ref{eq:Nimmo}) that 
$P_{(r)}(\vectx|\vecta)$ also depends only on $a_0, a_1, \dots, a_{r-1}$.
So we write $P_{(r)}(\vectx|a_0, a_1, \dots, a_{r-1})$ for $P_{(r)}(\vectx|\vecta)$.

\begin{prop}
\label{prop:fac-skewP1}
We define $R_{r/k}(\vectx|\vecta)$ by the relation
$$
P_{(r)}(\vectx,y|\vecta)
 =
\sum_{k=0}^\infty
 R_{r/k}(\vectx|\vecta) (y|\vecta)^k.
$$
Then we have we have
\begin{equation}
\label{eq:fac-skewP1}
R_{r/k}(\vectx|\vecta)
 =
\begin{cases}
 P_{(r)}(\vectx|-a_0, a_1,\dots, a_{r-1}) 
 &\text{if $k=0$,}
\\
 P_{(r-k)}(\vectx|0,a_{k+1}, \dots, a_{r-1})
 &\text{if $1 \le k \le r-1$,}
\\
 1
 &\text{if $k=r$,}
\\
 0
 &\text{if $k>r$.}
\end{cases}
\end{equation}
In particular, if $a_0 = 0$, then we have
\begin{equation}
\label{eq:fac-skewP2}
R_{r/k}(\vectx|\vecta)
 =
P_{(r-k)}(\vectx|0, a_{k+1}, \dots, a_{r-1}).
\end{equation}
\end{prop}

If $\vectx = (x)$ consists of a single variable, 
then $P_{(r)}(x|\vecta) = (x|\vecta)^r$.
Hence we obtain

\begin{corollary}
\label{cor:fac-skewP1}
If $\vectx = (x)$ consists of a single variable, then we have
$$
R_{r/k}(x|\vecta)
 = 
\begin{cases}
 (x+a_0) \prod_{i=1}^{r-1} (x - a_i) &\text{if $k=0$,} \\
 x \prod_{i=k+1}^{r-1} (x - a_i) &\text{if $1 \le k \le r-1$,} \\
 1 &\text{if $k=r$,} \\
 0 &\text{if $k>r$.}
\end{cases}
$$
\end{corollary}

In the proof of this proposition, 
we need the following relations for elementary symmetric polynomials $e_r(\vectx)$.

\begin{lemma}
\label{lem:rel-e}
\begin{enumerate}
\item[(1)]
If $k > 0$ and $l > 0$ then we have
$$
\sum_{m=1}^{r-1} e_{m-k}(x_1, \dots, x_m) e_{r-m-l}(x_{m+2}, \dots, x_r)
 =
e_{r-k-l}(x_1, \dots, x_r).
$$
\item[(2)]
If $l > 0$, then we have
\begin{multline*}
2 \sum_{m=1}^{r-1} e_m(x_1, \dots, x_m) e_{r-m-l}(x_{m+2}, \dots, x_r)
 + e_{r-l}(-x_1, x_2, \dots, x_r)
\\
 =
e_{r-l}(x_1, x_2, \dots, x_r).
\end{multline*}
\end{enumerate}
\end{lemma}

\begin{demo}{Proof}
For $1 \le a < b \le r$, we put $[a,b] = \{ a, a+1, \dots, b \}$ and 
denote by $\binom{[a,b]}{p}$ the set of $p$-element subsets of $[a,b]$.
If we put $x_I = \prod_{i \in I} x_i$ for $I \subset [a,b]$, then we have
$$
e_p(x_a, \dots, x_b) = \sum_{I \in \binom{[a,b]}{p}} x_I.
$$ 

(1)
We define a map 
$$
\varphi : \bigsqcup_{m=1}^{r-1} \binom{[1,m]}{m-k} \times \binom{[m+2,r]}{r-m-l}
 \to \binom{[1,r]}{r-k-l}
$$
by $\varphi(I,J) = I \sqcup J$.
Given $K \in \binom{[1,r]}{r-k-l}$, 
let $m+1$ be $(k+1)$st smallest element in the $(k+l)$-element subset of $[r] \setminus K$ 
and put $I = K \cap [1,m]$ and $J = K \cap [m+2,r]$.
Then the correspondence $K \mapsto (I,J)$ gives the inverse map of $\varphi$, 
and we obtain the desired identity.

(2)
Since we have
$$
e_{r-l}(-x_1, x_2, \dots, x_r)
 = 
- \sum_{\substack{ I \in \binom{[1,r]}{r-l} \\ 1 \in I}} x_I
+ \sum_{\substack{ I \in \binom{[1,r]}{r-l} \\ 1 \not\in I}} x_I,
$$
it is enough to show that
$$
\sum_{m=1}^{r-1} e_m(x_1, \dots, x_m) e_{r-m-l}(x_{m+2}, \dots, x_r)
 =
\sum_{\substack{K \in \binom{[1,r]}{r-l} \\ 1 \in K}} x_K.
$$
We define a map
$$
\psi : 
\bigsqcup_{m=1}^{r-1} \binom{[m+2,r]}{r-m-l}
 \to 
\left\{
 K \in \binom{[1,r]}{r-l} : 1 \in K
\right\}
$$
by $\psi(J) = [1,m] \sqcup J$.
Given $K \in \binom{[1,r]}{r-l}$ with $1 \in K$, 
let $m$ be the maximum integer $m$ satisfying $[1,m] \subset K$, 
and put $J = K \setminus [1,m]$.
Then the correspondence $K \mapsto J$ give the inverse map of $\psi$, 
and we obtain the desired identity.
\end{demo}

Now we prove Proposition~\ref{prop:fac-skewP1}.

\begin{demo}{Proof of Proposition~\ref{prop:fac-skewP1}}
We need to show that
\begin{align}
P_{(r)}(\vectx,y|a_0, \dots, a_{r-1})
 &=
P_{(r)}(\vectx|-a_0, a_1, \dots,a_{r-1})
\notag
\\
&\quad
 +
\sum_{k=1}^{r-1}
 P_{(r-k)}(\vectx|0, a_{k+1}, \dots, a_{r-1}) (y|a_0, \dots, a_{k-1})^k
\notag
\\
&\quad
 +
(y|a_0, \dots, a_{r-1})^r.
\label{eq:fac-skewP3}
\end{align}
We compare the coefficients of $P_{(k)}(\vectx) y^l$ in the expansions of the both hand sides, 
where $P_{(k)}(\vectx)$ is the Schur $P$-function.
We denote by $a_{k,l}$ and $b_{k,l}$ the coefficients of $P_{(k)}(\vectx) y^l$ 
on the left and right hand sides respectively.

Plugging $(u|\vecta)^r = \sum_{m=0}^r (-1)^{r-m} e_{r-m}(a_0, \dots, a_{r-1}) u^m$ 
into the definition (\ref{eq:Nimmo}) and using the multilinearity of Pfaffians, 
we have
\begin{equation}
\label{eq:expansion1}
P_{(r)}(\vectx|\vecta)
 =
\begin{cases}
\sum_{m=1}^r (-1)^{r-m} e_{r-m}(a_0, \dots, a_{r-1}) P_{(m)}(\vectx) &\text{if $n$ is even,} \\
\sum_{m=0}^r (-1)^{r-m} e_{r-m}(a_0, \dots, a_{r-1}) P_{(m)}(\vectx) &\text{if $n$ is odd.}
\end{cases}
\end{equation}
Since $Q_{(r)}(\vectx) = 2 P_{(r)}(\vectx)$, 
it follows from (\ref{eq:GF-SchurQ1}) that
$$
1 + 2 \sum_{r=1}^\infty P_{(r)}(\vectx,y) z^r 
 =
\prod_{i=1}^n \frac{1 + x_i z}{1 - x_i z}
\cdot \frac{1 + y z}{1 - y z}
 =
\left( 
 1 + 2 \sum_{r=1}^\infty P_{(r)}(\vectx) z^r 
\right)
\left(
 1 + 2 \sum_{r=1}^\infty y^r z^r
\right).
$$
Equating the coefficient of $z^r$, we have
\begin{equation}
\label{eq:expansion2}
P_{(r)}(\vectx,y)
 =
P_{(r)}(\vectx) + 2 \sum_{h=1}^{r-1} P_{(r-h)}(\vectx) y^h + y^r.
\end{equation}

Using (\ref{eq:expansion1}) and (\ref{eq:expansion2}), we have
\begin{multline*}
P_{(r)}(\vectx,y|\vecta)
\\
=
\left\{
\begin{array}{l}
e_r(a_0, \dots, a_{r-1})
\\
 +
\sum_{m=1}^r (-1)^{r-m} e_{r-m}(a_0, \dots, a_{r-1})
\left(
 P_{(m)}(\vectx) + 2 \sum_{l=1}^{m-1} P_{(m-l)}(\vectx) y^l + y^m
\right)
\\
 \hfill\text{if $n$ is even,}
\\
\sum_{m=1}^r (-1)^{r-m} e_{r-m}(a_0, \dots, a_{r-1})
\left(
 P_{(m)}(\vectx) + 2 \sum_{l=1}^{m-1} P_{(m-l)}(\vectx) y^l + y^m
\right)
\\
\hfill\text{if $n$ is odd.}
\end{array}
\right.
\end{multline*}
Hence the coefficient $a_{k,l}$ of $P_{(k)}(\vectx) y^l$ in the left hand side of (\ref{eq:fac-skewP3}) 
is given by
$$
a_{k,l}
 =
\begin{cases}
 e_r(a_0, \dots, a_{r-1}) &\text{if $k=0$, $l=0$ and $n$ is even,} \\
 0                        &\text{if $k=0$, $l=0$ and $n$ is odd,} \\
 (-1)^{r-l} e_{r-l}(a_0, \dots, a_{r-1})
                          &\text{if $k=0$ and $l>0$,} \\
 2 (-1)^{r-k-l} e_{r-k-l}(a_0, \dots, a_{r-1})
                          &\text{if $k>0$.}
\end{cases}
$$
In a similar manner, we can compute the coefficient $b_{k,l}$ in the right hand side of (\ref{eq:fac-skewP3}) 
and see that
\begin{enumerate}
\item[(a)]
if $k=l=0$, then
$$
b_{0,0}
 = 
\begin{cases}
 e_r(a_0, \dots, a_{r-1}) &\text{if $n$ is even,} \\
 0                        &\text{if $n$ is odd,} \\
\end{cases}
$$
\item[(b)]
if $k=0$ and $l>0$, then
$$
b_{0,l} = (-1)^{r-l} e_{r-l}(a_0, \dots, a_{r-1}),
$$
\item[(c)]
if $k>0$ and $l=0$, then
$$
b_{k,0}
 =
(-1)^{r-k}
\left(
\begin{array}{l}
 e_{r-k}(-a_0, a_1, \dots, a_{r-1})
\\
 +
 \displaystyle\sum_{m=1}^{r-1} e_{r-m-k}(0, a_{m+1}, \dots, a_{r-1}) e_m(a_0, \dots, a_{m-1})
\end{array}
\right),
$$
\item[(d)]
if $k>0$ and $l>0$, then
$$
b_{k,l}
 =
(-1)^{r-k-l}
\sum_{m=1}^{r-1}
 e_{r-m-k}(0, a_{m+1}, \dots, a_{r-1}) e_{m-l}(a_0, \dots, a_{m-1}).
$$
\end{enumerate}
Now by using Lemma~\ref{lem:rel-e}, we see $a_{k,l} = b_{k,l}$ and obtain (\ref{eq:fac-skewP3}).
\end{demo}

By a similar argument to the proof of Proposition~\ref{prop:skewP-3}, 
we can derive the determinant formula for $P_{\lambda/\mu,p}(x|\vecta)$ for a single variable $x$.

\begin{theorem}
\label{thm:fac-skewP}
Let $\vecta = (a_0, a_1, \dots)$ be factorial parameters.
For two strict partitions $\lambda$ of lenth $l$ and $\mu$ of length $m$, 
the factorial skew $P$-function $P_{\lambda/\mu,p}(x|\vecta)$ in a single variable $x$ 
is given as follows:
\begin{enumerate}
\item[(1)]
We have $P_{\lambda/\mu,p}(x|\vecta) = 0$ unless $\lambda \supset \mu$ and $m = l$ or $l-1$.
\item[(2)]
If $\lambda \supset \mu$ and $m = l$ or $l-1$, then we have
$$
P_{\lambda/\mu,p}(x|\vecta)
 = 
\det \left( R_{\lambda_i/\mu_j}(x|\vecta) \right)_{1 \le i, j \le l}.
$$
\end{enumerate}
\end{theorem}

\begin{demo}{Proof}
It follows from Proposition~\ref{prop:skewP-1} that 
$P_{\lambda/\mu,p}(x|\vecta) = 0$ unless $\lambda \supset \mu$.
By a similar argument to the proof of Proposition~\ref{prop:skewP-2}, we can show the following:
\begin{enumerate}
\item[(a)]
if $l \equiv p$ and $m \equiv p \bmod 2$, then
$$
P_{\lambda/\mu,p}(x|\vecta)
 =
\begin{cases}
 \det \left( R_{\lambda_i/\mu_j}(x|\vecta) \right)_{1 \le i, j \le l}
 &\text{if $l=m$,}
\\
 0 &\text{otherwise.}
\end{cases}
$$
\item[(b)]
if $l \equiv p$ and $m \not\equiv p \bmod 2$, 
or if $l \not\equiv p$ and $m \equiv p \bmod 2$, then
$$
P_{\lambda/\mu,p}(x|\vecta)
 =
\begin{cases}
 \det \left( R_{\lambda_i/\mu_j}(x|\vecta) \right)_{1 \le i, j \le l}
 &\text{if $l=m+1$,}
\\
 0 &\text{otherwise,}
\end{cases}
$$
where $\mu_l = 0$.
\end{enumerate}
It remains to consider the case where $l \not\equiv p$ and $m \not\equiv p \bmod 2$.
By the definition (\ref{eq:skewP}) and the multilinearity of Pfaffians we have
\begin{align*}
P_{\lambda/\mu,p}(x|\vecta)
 &=
\Pf \begin{pmatrix}
 O_l & T_\lambda(\vectx) & T_{\lambda/0}(\vectx) & M_{\lambda/\mu}(\vectx) \\
 -\trans T_\lambda(\vectx) & 0 & 1 & O \\
 -\trans T_{\lambda/0}(\vectx) & -1 & 0 & O \\
 -\trans M_{\lambda/\mu}(\vectx) & O & O & O_m
\end{pmatrix}
\\
 &=
\Pf \begin{pmatrix}
 O_l & T_\lambda(x) & T_{\lambda/0}(x) & M_{\lambda/\mu}(x) \\
 -\trans T_\lambda(x) & 0 & 0 & O \\
 -\trans T_{\lambda/0}(x) & 0 & 0 & O \\
 -\trans M_{\lambda/\mu}(x) & O & O & O_m
\end{pmatrix}
\\
&\quad
+
\Pf \begin{pmatrix}
 O_l & T_\lambda(x) & 0 & M_{\lambda/\mu}(x) \\
 -\trans T_\lambda(x) & 0 & 1 & O \\
 0 & -1 & 0 & O \\
 -\trans M_{\lambda/\mu} & O & O & O_m
\end{pmatrix},
\end{align*}
where $T_\lambda(x)$ and $T_{\lambda/0}(x)$ is the column vector with $i$th entry 
$P_{(\lambda_i)}(x|\vecta)$ and $R_{\lambda_i/0}(x|\vecta)$ respectively.
Since $P_{(r)}(x|\vecta) = (x-a_0) \prod_{i=1}^{r-1} (x - a_i)$ and 
$R_{r/0}(x|\vecta) = (x+a_0) \prod_{i=1}^{r-1} (x - a_i)$ by Corollary~\ref{cor:fac-skewP1}, 
we see that $T_\lambda(x)$ and $T_{\lambda/0}(x)$ are linearly independent.
Hence the first Pfaffian vanishes.
By expanding the second Pfaffian along the $(l+2)$nd row/column, we have
$$
P_{\lambda/\mu,p}(x)
 =
\Pf \begin{pmatrix}
 O_l & M_{\lambda/\mu}(x) \\
 -\trans M_{\lambda/\mu}(x) & O_m
\end{pmatrix}
 =
\begin{cases}
 \det \left( R_{\lambda_i/\mu_j}(x) \right)_{1 \le i, j \le l} &\text{if $l=m$,} \\
 0 &\text{otherwise.}
\end{cases}
$$
This completes the proof.
\end{demo}

We can use this theorem to provide a lattice path proof to the tableau description 
of factorial $P$- and $Q$-functions 
given in \cite[Theorem~2.1]{CI} and \cite[Theorem~4.3]{Ivanov1} (the case where $a_0=0$).

\subsection{%
Modified Pieri coefficients
}

In this subsection, we give a combinatorial description to the modified Pieri coefficients 
for factorial $P$-functions.
Recall that the skew shifted diagram $S(\lambda/\mu) = S(\lambda) \setminus S(\mu)$ 
is called a \emph{border strip} 
if it is connected and contains no $2 \times 2$ block of cells.

\begin{theorem}
\label{thm:fac-Pieri}
Let $\vecta = (a_0, a_1, \dots)$ be factorial parameters.
We define the modified Pieri coefficient $c^\lambda_{\mu,r}(\vecta)$ by the relation
\begin{equation}
P_\mu(\vectx|\vecta) \cdot q_r(\vectx)
 =
\sum_\lambda c^\lambda_{\mu,r}(\vecta) P_\lambda(\vectx|\vecta),
\end{equation}
where $q_r(\vectx) = Q_{(r)}(\vectx)$ is Schur's $Q$-function, 
and $\lambda$ runs over all strict partitions.
For two strict partitions $\lambda$ and $\mu$, we consider the generating function 
of modified Pieri coefficients
$$
c^\lambda_\mu(z|\vecta)
 =
\sum_{r=0}^\infty c^\lambda_{\mu,r}(\vecta) z^r.
$$
Then we have
\begin{enumerate}
\item[(1)]
If the skew shifted diagram $S(\lambda/\mu)$ contains a $2 \times 2$ block of cells, 
then we have $c^\lambda_\mu(z|\vecta) = 0$.
\item[(2)]
Suppose that $S(\lambda/\mu)$ contains no $2 \times 2$ block of cells.
Let
$$
S(\lambda/\mu)
 =
\bigsqcup_{i=1}^r
 S((\lambda_{m(i)}, \dots, \lambda_{M(i)})/(\mu_{m(i)}, \dots, \mu_{M(i)}))
$$
be the decomposition of $S(\lambda)/S(\mu)$ into a disjoint union of border strips, 
where $m(1) \le M(1) < m(2) \le M(2) < \dots < m(r) \le M(r)$.
Then we have
\begin{equation}
\label{eq:fac-Pieri}
c^\lambda_\mu(z|\vecta)
 =
\prod_{k \in K}
 \frac{ 1 + a_{\mu_k} z }
      { 1 - a_{\mu_k} z }
\prod_{i=1}^r
 \frac{ 2 z^{\lambda_{m(i)} - \mu_{M(i)}} }
       { \prod_{j=\mu_{M(i)}}^{\lambda_{m(i)}} ( 1 - a_j z ) }
\end{equation}
where
$$
K = \begin{cases}
 \{ k : 1 \le k \le l(\mu), \, \lambda_k = \mu_k \} &\text{if $n+l(\mu)$ is even,} \\
 \{ k : 1 \le k \le l(\mu)+1, \, \lambda_k = \mu_k \} &\text{if $n+l(\mu)$ is odd.}
\end{cases}
$$
\item[(3)]
In particular, the modified Pieri coefficient $c^\lambda_{\mu,r}(\vecta)$ is a polynomial 
in the factorial parameters $a_0, a_1, \dots$ with nonnegative integer coefficients.
\end{enumerate}
\end{theorem}

\begin{example}
Let $\lambda = (8,6,4,3,2)$ and $\mu = (6,5,4,2,1)$.
Then the skew shifted Young diagram $S(\lambda/\mu)$ is decomposed into a disjoint union 
of border strips as follows:
\begin{align*}
S(\lambda/\mu)
 &=
\raisebox{-27.5pt}{
\setlength\unitlength{1.3pt}
\begin{picture}(80,50)
\dashline{1.5}(0,50)(60,50)
\put(60,50){\line(1,0){20}}
\dashline{1.5}(0,40)(60,40)
\put(60,40){\line(1,0){20}}
\dashline{1.5}(10,30)(60,30)
\put(60,30){\line(1,0){10}}
\dashline{1.5}(20,20)(50,20)
\put(50,20){\line(1,0){10}}
\dashline{1.5}(30,10)(50,10)
\put(50,10){\line(1,0){10}}
\dashline{1.5}(40,0)(50,0)
\put(50,0){\line(1,0){10}}
\dashline{1.5}(0,40)(0,50)
\dashline{1.5}(10,30)(10,50)
\dashline{1.5}(20,20)(20,50)
\dashline{1.5}(30,10)(30,50)
\dashline{1.5}(40,0)(40,50)
\put(50,0){\line(0,1){20}}
\dashline{1.5}(50,20)(50,50)
\put(60,0){\line(0,1){20}}
\dashline{1.5}(60,20)(60,30)
\put(60,30){\line(0,1){20}}
\put(70,30){\line(0,1){20}}
\put(80,40){\line(0,1){10}}
\end{picture}
}
\\
 &=
S((8,6)/(6,5)) \sqcup S((3,2)/(2,1)).
\end{align*}
Since we have
$$
K = \begin{cases}
 \{ 3 \} &\text{if $n$ is even,} \\
 \{ 3, 6 \} &\text{if $n$ is odd,}
\end{cases}
$$
we obtain
$$
c^{(8,6,4,3,2)}_{(6,5,4,2,1)}(z|\vecta)
 =
\begin{cases}
 \dfrac{ 1 + a_4 z }{ 1 - a_4 z }
 \cdot
 \dfrac{ 2 z^{8-5} }{ \prod_{i=5}^8 (1 - a_i z) }
 \cdot
 \dfrac{ 2 z^{3-1} }{ \prod_{i=1}^3 (1 - a_i z) }
 &\text{if $n$ is even,}
\\[10pt]
 \dfrac{ 1 + a_4 z }{ 1 - a_4 z }
 \cdot
 \dfrac{ 1 + a_0 z }{ 1 - a_0 z }
 \cdot
 \dfrac{ 2 z^{8-5} }{ \prod_{i=5}^8 (1 - a_i z) }
 \cdot
 \dfrac{ 2 z^{3-1} }{ \prod_{i=1}^3 (1 - a_i z) }
 &\text{if $n$ is odd.}
\end{cases}
$$
\end{example}

Setting $a_0 = a_1 = \dots = 0$, we recover the Pieri rule for Schur $P$-functions.

\begin{corollary}
\label{cor:nonfac-Pieri}
(Morris \cite[Theorem~1]{Morris})
For a strict partition $\mu$ and a nonnegative integer $k$, we have
$$
P_\mu(\vectx) \cdot q_r(\vectx)
 =
\sum_\lambda 2^{a(\lambda,\mu)} P_\lambda(\vectx),
$$
where $\lambda$ runs over all strict partitions such that $|\lambda| - |\mu| = r$ 
and $S(\lambda/\mu)$ contains no $2 \times 2$ block 
and $a(\lambda,\mu)$ is the number of connected components of $S(\lambda/\mu)$.
\end{corollary}

Now we use Theorem~\ref{thm:Pieri} to give a proof of Theorem~\ref{thm:fac-Pieri}.
Let $b^s_r(z|\vecta)$ be the coefficient in the expansion
\begin{equation}
\label{eq:fac-Pieri1}
(u|\vecta)^r \cdot \frac{ 1 + u z }{ 1 - u z}
 = 
\sum_{s \ge 0} b^s_r(z|\vecta) (u|\vecta)^s.
\end{equation}
The following lemma gives an explicit formula for $b^s_r(z|\vecta)$.

\begin{lemma}
\label{lem:fac-Pieri1}
For two nonnegative integers $r$ and $s$, we have
$$
b^s_r(z|\vecta)
 =
\begin{cases}
 \dfrac{ 1 + a_r z }{ 1 - a_r z } &\text{if $s=r$,} \\
 \dfrac{ 2 z^{s-r} }{ \prod_{j=r}^s (1 - a_j z) } &\text{if $s > r$,} \\
 0 &\text{otherwise.}
\end{cases}
$$
\end{lemma}

\begin{demo}{Proof}
We need to prove
\begin{equation}
\label{eq:expansion}
(u|\vecta)^r \frac{ 1 + u z }{ 1 - u z}
 =
\frac{ 1 + a_r u }{ 1 - a_r u } (u|\vecta)^r
 +
\sum_{s = r+1}^\infty
 \frac{ 2 z^{s-r} }
      { \prod_{j=r}^s (1 - a_j z) }
 (u|\vecta)^s.
\end{equation}
By dividing the both sides of (\ref{eq:expansion}) by $\prod_{i=0}^{r-1} (u-a_i)$ 
and then by shifting the indices of factorial parameters, 
we may assume $r = 0$.
This case follows from Lemma~\ref{lem:fac-dual}.
\end{demo}

We prove Theorem~\ref{thm:fac-Pieri} 
by computing the determinant given in (\ref{eq:Pieri}) of Theorem~\ref{thm:Pieri}.

\begin{demo}{Proof of Theorem~\ref{thm:fac-Pieri}}
By (\ref{eq:Pieri}), we see that nonzero $c^\lambda_\mu(z|\vecta)$ 
is equal to the determinant whose $(i,j)$ entry is equal to 
$b^{\lambda_i}_{\mu_j} = b^{\lambda_i}_{\mu_j}(z|\vecta)$.

\begin{claim}
\label{claim:1}
We have $c^\lambda_\mu(z|\vecta) = 0$ unless $S(\lambda) \supset S(\mu)$.
\end{claim}

\begin{demo}{Proof}
If then there exists an index $k$ such that $\lambda_k < \mu_k$, then
we have $\lambda_i < \lambda_k < \mu_k < \mu_j$ for $i \ge k$ and $j \le k$.
By Lemma~\ref{lem:fac-Pieri1}, we have $b^{\lambda_i}_{\mu_j} = 0$ for $i \ge k$ and $j \le k$, 
thus $c^\lambda_\mu = 0$.
\end{demo}

In what follows we assume that $S(\lambda) \supset S(\mu)$.
In this case, by Theorem~\ref{thm:Pieri}, we have
$$
c^\lambda_\mu(z|\vecta)
 =
\begin{cases}
 \det B^\lambda_\mu &\text{if $n+l(\mu)$ is even and $l(\lambda) = l(\mu)$,} \\
 \det B^{\lambda^0}_{\mu^0} &\text{if $n+l(\mu)$ is odd and $l(\lambda) = l(\mu)$, } \\
 \det B^\lambda_{\mu^0} &\text{if $l(\lambda) = l(\mu)+1$,} \\
 0 &\text{otherwise.}
\end{cases}
$$
If $n+l(\mu)$ is odd and $l(\lambda) = l(\mu) = l$, 
then $b^0_{\mu_j} = 0$ for $1 \le j \le l$ by Lemma~\ref{lem:fac-Pieri1}.
By expanding the determinant along the last row, we have
$$
\det B^{\lambda^0}_{\mu^0}
 =
b^0_0 \cdot \det B^\lambda_\mu.
$$
Hence it is enough to compute the determinant $\det B^\lambda_{\mu^*}$, where $\mu^* = \mu$ or $\mu^0$.
By abuse of the notation, we write simply $B^\lambda_\mu$ for $B^\lambda_\mu$ and $B^\lambda_{\mu^0}$ 
in the following.

First we consider the case $S(\lambda/\mu)$ is not connected.
In this case, there exists an index $k$ such that $\lambda_{k+1} < \mu_k$ or $\lambda_k = \mu_k$.

\begin{claim}
\label{claim:2}
Suppose that $S(\lambda) \supset S(\mu)$.
\begin{enumerate}
\item[(1)]
If there exists an index $k$ such that $\lambda_{k+1} < \mu_k$, 
then we have
$$
\det B^\lambda_\mu
 = 
\det B^{\lambda'}_{\mu'} \cdot \det B^{\lambda''}_{\mu''},
$$
where
\begin{gather*}
\lambda' = (\lambda_1, \dots, \lambda_k),
\quad
\mu' = (\mu_1, \dots, \mu_k),
\\
\lambda'' = (\lambda_{k+1}, \dots, \lambda_l),
\quad
\mu'' = (\mu_{k+1}, \dots, \mu_l).
\end{gather*}
\item[(2)]
If there exists an index $k$ such that $\lambda_k = \mu_k > 0$, then we have
$$
\det B^\lambda_\mu
 =
\det B^{\lambda'}_{\mu'} \cdot b^{\lambda_k}_{\mu_k} \cdot \det B^{\lambda''}_{\mu''},
$$
where
\begin{gather*}
\lambda' = (\lambda_1, \dots, \lambda_{k-1}),
\quad
\mu' = (\mu_1, \dots, \mu_{k-1}),
\\
\lambda'' = (\lambda_{k+1}, \dots, \lambda_l),
\quad
\mu'' = (\mu_{k+1}, \dots, \mu_l).
\end{gather*}
\end{enumerate}
\end{claim}

\begin{demo}{Proof}
(1)
If $i \ge k+1$ and $j \le k$, then we have $\lambda_i \le \lambda_{k+1} < \mu_k \le \mu_j$ 
and $b^{\lambda_i}_{\mu_j} = 0$, thus 
$$
\det B^\lambda_\mu
 =
\det \begin{pmatrix}
 B^{\lambda'}_{\mu'} & * \\
 O & B^{\lambda''}_{\mu''}
\end{pmatrix}
 =
\det B^{\lambda'}_{\mu'} \cdot \det B^{\lambda''}_{\mu''}.
$$

(2)
By a similar consideration, we have
$$
\det B^\lambda_\mu
 =
\det \begin{pmatrix}
 B^{\lambda'}_{\mu'} & * & * \\
 O                   & b^{\lambda_k}_{\mu_k} & * \\
 O                   & 0                     & B^{\lambda''}_{\mu''}
\end{pmatrix}
 =
\det B^{\lambda'}_{\mu'} \cdot b^{\lambda_k}_{\mu_k} \cdot \det B^{\lambda''}_{\mu''}.
$$
\end{demo}

Now we consider the case where $S(\lambda/\mu)$ is connected.

\begin{claim}
\label{claim:3}
Suppose that $S(\lambda) \supset S(\mu)$ and $S(\lambda/\mu)$ is connected.
If $S(\lambda/\mu)$ contains a $2 \times 2$ block of cells, 
then we have $\det B^\lambda_\mu = 0$.
\end{claim}

\begin{demo}{Proof}
If $S(\lambda/\mu)$ contains a $2 \times 2$ square, 
then there exists an index $k$ such that $\lambda_{k+1} > \mu_k$.
We take the minimum such index $k$.
Then we have
$$
\lambda_1 > \mu_1 = \lambda_2 > \mu_2 = \lambda_3 > \dots > \mu_{k-1} = \lambda_k > \lambda_{k+1} > \mu_k.
$$
(Since $S(\lambda/\mu)$ is connected, we have $\lambda_{i+1} \ge \mu_i$ if $\lambda_{i+1} > 0$.)
It follows from Lemma~\ref{lem:fac-Pieri1} that, if $t > s > r$, then
\begin{equation}
\label{eq:relation}
b^t_r
 =
b^s_r
 \cdot
 \frac{ z^{t-s} }
      { (1 - a_{s+1} z) \cdots (1 - a_t z) }.
\end{equation}

We proceed by induction on $k$.
If $k=1$, then the first row of $B^\lambda_\mu$ is a scalar multiple of the second row by (\ref{eq:relation}), 
and $\det B^\lambda_\mu = 0$.
If $k > 1$, then by subtracting the $(k+1)$st row multiplied 
by $z^{\lambda_k-\lambda_{k+1}}/(1-a_{\lambda_{k+1}+1}z) \cdots (1-a_{\lambda_k}z)$ from the $k$th row 
in $\det B^\lambda_\mu$, 
and then by expanding the resulting determinant along the $k$th row, 
we have
$$
\det B^\lambda_\mu
 =
(-1)^{k+(k-1)}
b^{\lambda_k}_{\mu_{k-1}}
\det
B^{\lambda'}_{\mu'},
$$
where $\lambda'$ (resp. $\mu'$) is the strict partition obtained 
from $\lambda$ (resp. $\mu$) by removing $\lambda_k$ (resp. $\mu_{k-1}$), i.e.,
$$
\lambda' = (\lambda_1, \dots, \lambda_{k-1}, \lambda_{k+1}, \dots, \lambda_l),
\quad
\mu' = (\mu_1, \dots, \mu_{k-2}, \mu_k, \dots, \mu_l).
$$
Since $\det B^{\lambda'}_{\mu'} = 0$ by the induction hypothesis, 
we obtain $\det B^\lambda_\mu = 0$.
\end{demo}

Now it remains to compute $\det B^\lambda_\mu$ when $S(\lambda/\mu)$ is a border strip.

\begin{claim}
\label{claim:4}
If $S(\lambda/\mu)$ is a border strip, i.e.,
$$
\lambda_1 > \mu_1 = \lambda_2 > \mu_2 = \lambda_3 > \dots > \mu_{l-1} = \lambda_l > \mu_l,
$$
then we have
$$
\det B^\lambda_\mu
 =
\frac{ 2 z^{\lambda_1 - \mu_l} }
     { \prod_{i=\mu_l}^{\lambda_1} (1 - a_i z) }.
$$
\end{claim}

\begin{demo}{Proof}
We proceed by induction on $l$.
A direct computation shows the cases $l=1$ and $l=2$.
So we assume $l \ge 3$.
By Lemma~\ref{lem:fac-Pieri1}, we have
\begin{align*}
b^{\lambda_i}_{\mu_l}
 &=
\frac{ z^{\mu_{l-1} - \mu_l} }
     { \prod_{j=\mu_l}^{\mu_{l-1}-1} (1 - a_j z) }
\cdot b^{\lambda_i}_{\mu_{l-1}}
\quad(1 \le i \le l-1),
\\
b^{\lambda_l}_{\mu_l}
 &=
\frac{ z^{\mu_{l-1} - \mu_l} }
     { \prod_{j=\mu_l}^{\mu_{l-1}-1} (1 - a_j z) }
\cdot
\frac{ 2 }
     { 1 - a_{\lambda_l} }
\cdot
b^{\lambda_l}_{\mu_{l-1}}.
\end{align*}
Pull out the common factor 
$z^{\mu_{l-1} - \mu_l} / \prod_{j=\mu_l}^{\mu_{l-1}-1} (1 - a_j z)$ 
from the last column of $\det B^\lambda_\mu$, 
and then subtract the $l$th column from the $(l-1)$st column.
Since we have
$$
b^{\lambda_l}_{\mu_{l-1}}
 - 
\frac{ 2 }
     { 1 - a_{\lambda_l} }
\cdot
b^{\lambda_l}_{\mu_{l-1}}
 =
-1,
$$
we expand the resulting determinant along the $(l-1)$st column to see
$$
\det B^\lambda_\mu
 =
\frac{ z^{\mu_{l-1} - \mu_l} }
     { \prod_{j=\mu_l}^{\mu_{l-1}-1} (1 - a_j z) }
\cdot
\det B^{\lambda'}_{\mu'},
$$
where $\lambda' = (\lambda_1, \dots, \lambda_{l-1})$ and $\mu' = (\mu_1, \dots, \mu_{l-1})$.
Using the induction hypothesis, we obtain Claim~\ref{claim:4}.
\end{demo}

Combining the above claims together completes the proof.
\end{demo}

Based on Theorem~\ref{thm:fac-Pieri} (3) and some experimental evidence, 
we propose the following conjecture.

\begin{conjecture}
\label{conj:fac-Pieri}
We define $F^\lambda_{\mu,\nu}(\vecta)$ by the formula
\begin{equation}
\label{eq:conj1}
P_\mu(\vectx|\vecta) P_\nu(\vectx)
 =
\sum_\lambda f^\lambda_{\mu,\nu}(\vecta) P_\lambda(\vectx|\vecta),
\end{equation}
where $\lambda$ runs over all strict partitions.
Then the coefficient $f^\lambda_{\mu,\nu}(\vecta)$ is a polynomial in 
$\vecta$ with nonnegative integer coefficients.
More generally, if we expand the product of factorial $P$-functions 
corresponding to different factorial parameters 
$\vecta = (a_0, a_1, \dots)$ and $-\vectb = (-b_0, -b_1, \dots)$ 
as a linear combination of $P_\lambda(\vectx|\vecta)$s:
\begin{equation}
\label{eq:conj2}
P_\mu(\vectx|\vecta) P_\nu(\vectx|-\vectb)
 =
\sum_\lambda f^\lambda_{\mu,\nu}(\vecta,\vectb) P_\lambda(\vectx|\vecta),
\end{equation}
then the coefficient $f^\lambda_{\mu,\nu}(\vecta,\vectb)$ is a polynomial in $\vecta$ and $\vectb$ 
with nonnegative integer coefficients.
\end{conjecture}

Cho--Ikeda \cite[Theorem~4.6]{CI} gave a combinatorial formula for the Pieri-type coefficients 
$f^\lambda_{\mu,(r)}(\vecta,-\vecta)$, which implies that $f^\lambda_{\mu,(r)}(\vecta,-\vecta)$ 
is a polynomial in $a_i \pm a_j$ with $i>j$ with nonnegative integer coefficients.

\begin{remark}
Let $s_\lambda(\vectx|\vecta)$ be the factorial Schur function with factorial parameters $\vecta$, 
and expand
$$
s_\mu(\vectx|\vecta) s_\nu(\vectx|-\vectb)
 =
\sum_\lambda m^\lambda_{\mu,\nu}(\vecta,\vectb) s_\lambda(\vectx|\vecta).
$$
Then Molev--Sagan \cite[Theorem~3.1]{MS} gave a combinatorial formula for 
the coefficient $m^\lambda_{\mu,\nu}(\vecta,\vectb)$, 
which implies that $m^\lambda_{\mu,\nu}(\vecta,\vectb)$ is a polynomial in $\vecta$ and $\vectb$ 
with nonnegative integer coefficients.
\end{remark}

\section{%
$\boldsymbol{P}$-functions associated to classical root systems
}

In this section we show that the Hall--Littlewood functions at $t = -1$ 
associated to the classical root systems can be written as 
generalized $P$-functions associated to certain polynomial sequences.

\subsection{%
Hall--Littlewood function associated root systems
}

Macdonald \cite{Macdonald00} generalized the definition of Hall--Littlewood functions 
to any root system.
Let $\Phi$ be a root system in a Euclidean vector space $V$ 
and fix a positive system $\Phi^+$.
We denote by $\Lambda$ and $\Lambda^+$ the weight lattice and the set of dominant weights respectively.
Let $K = \Rat(t)$ be the rational function field in an indeterminate $t$ 
and $K[\Lambda]$ the group algebra of $\Lambda$ with basis $\{ e^\lambda : \lambda \in \Lambda \}$ over $K$.
Let $W$ be the Weyl group of $\Phi$. 
Then the Hall--Littlewood functions associated to the root system $\Phi$ are defined as follows:

\begin{definition}
\label{def:HL}
The Hall--Littlewood function $\PP^\Phi_\lambda \in K[\Lambda]$ corresponding 
to a dominant weight $\lambda \in \Lambda^+$ 
is defined by
\begin{equation}
\label{eq:def_HL}
\PP^\Phi_\lambda
 = 
\frac{1}{W_\lambda(t)}
\sum_{w \in W}
 w \left(
  e^\lambda \prod_{\alpha \in R^+} \frac{1 - t e^{-\alpha}}{1 - e^{-\alpha}}
 \right),
\end{equation}
where $W_\lambda = \{ w \in W : w \lambda = \lambda \}$ is the stabilizer of $\lambda$ in $W$, 
and $W_\lambda(t) = \sum_{w \in W^\lambda} t^{l(w)}$ is the Poincar\'e polynomial of $W_\lambda$.
\end{definition}

In this section we consider the root systems of types $X_n = B_n$, $C_n$ and $D_n$.
Let $V$ be the $n$-dimensional Euclidean vector space with orthonormal basis $\ep_1, \dots, \ep_n$.
We put $x_i = e^{\ep_i}$ for $1 \le i \le n$ 
and write $\PP^\Phi_\lambda = \PP^\Phi_\lambda(\vectx;t)$.
Let $\Phi(X_n) \subset V$ be the root system of type $X_n$ with the positive system $\Phi^+(X_n)$ given by
\begin{align*}
\Phi^+(B_n)
 &=
\{ \ep_i \pm \ep_j : 1 \le i < j \le n \} \cup \{ \ep_i : 1 \le i \le n \},
\\
\Phi^+(C_n)
 &=
\{ \ep_i \pm \ep_j : 1 \le i < j \le n \} \cup \{ 2 \ep_i : 1 \le i \le n \},
\\
\Phi^+(D_n)
 &=
\{ \ep_i \pm \ep_j : 1 \le i < j \le n \}. 
\end{align*}
Then the set $\Lambda^+(X_n)$ of dominant weights is given by
\begin{align*}
\Lambda^+(B_n)
 &=
\left\{
 \sum_{i=1}^n \lambda_i \ep_i : 
 (\lambda_i)_{i=1}^n \in \Int^n \cup (\Int+1/2)^n, \ \lambda_1 \ge \lambda_2 \ge \dots \ge \lambda_n \ge 0 
\right\},
\\
\Lambda^+(C_n)
 &=
\left\{
 \sum_{i=1}^n \lambda_i \ep_i : 
 (\lambda_i)_{i=1}^n \in \Int^n, \ \lambda_1 \ge \lambda_2 \ge \dots \ge \lambda_n \ge 0 
\right\},
\\
\Lambda^+(D_n)
 &=
\left\{
 \sum_{i=1}^n \lambda_i \ep_i :
 (\lambda_i)_{i=1}^n \in \Int^n \cup (\Int+1/2)^n, \ \lambda_1 \ge \lambda_2 \ge \dots \lambda_{n-1} \ge |\lambda_n|
\right\},
\end{align*}
where $\Int+1/2 = \{ r+1/2 : r \in \Int \}$.
We identify a partition $\lambda$ of length $l \le n$ with a dominant weight 
$\lambda_1 \ep_1 + \dots + \lambda_l \ep_l \in \Lambda^+(X_n)$.

We note that a reciprocal Laurent polynomial $g(z) = \sum_{i=-d}^d a_i z^i$ with $a_i = a_{-i}$ and $a_d \neq 0$ 
can be written as $g(z) = f(z+z^{-1})$ for some polynomial of degree $d$.
We use the notation $\vectx + \vectx^{-1} = (x_1 + x_1^{-1}, \dots, x_n + x_n^{-1})$.
The following is the main theorem of this section.

\begin{theorem}
\label{thm:N=HL-BCD}
\begin{enumerate}
\item[(1)]
Let $\mathcal{F}^B = \{ f^B_d \}_{d=0}^\infty$ be the sequence of polynomials defined by
$$
f^B_0 = 1,
\quad
f^B_d(x+x^{-1})
=
(x^d - x^{-d}) \frac{ x^{1/2} + x^{-1/2} }{ x^{1/2} - x^{-1/2} }
\quad(d \ge 1).
$$
For a strict partition $\lambda$ of length $l \le n$, we have
\begin{equation}
\label{eq:N=HL-B}
\PP^{\Phi(B_n)}_\lambda (\vectx;-1)
 =
P^{\mathcal{F}^B}_\lambda(\vectx+\vectx^{-1}).
\end{equation}
\item[(2)]
Let $\mathcal{F}^C = \{ f^C_d \}_{d=0}^\infty$ be the sequence of polynomials defined by
$$
f^C_0 = 1,
\quad
f^C_d(x+x^{-1})
=
(x^d - x^{-d}) \frac{ x + x^{-1} }{ x - x^{-1} }
\quad(d \ge 1).
$$
For a strict partition $\lambda$ of length $l \le n$, we have
\begin{equation}
\label{eq:N=HL-C}
\PP^{\Phi(C_n)}_\lambda (\vectx;-1)
 =
P^{\mathcal{F}^C}_\lambda(\vectx+\vectx^{-1}).
\end{equation}
\item[(3)]
Let $\mathcal{F}^D = \{ f^D_d \}_{d=0}^\infty$ be the sequence of polynomials defined by
$$
f^D_0 = 1,
\quad
f^D_d(x+x^{-1})
 =
x^d + x^{-d}
\quad(d \ge 1).
$$
For a strict partition $\lambda$ of length $l < n$, we have
\begin{equation}
\label{eq:N=HL-D1}
\PP^{\Phi(D_n)}_\lambda (\vectx;-1)
 =
P^{\mathcal{F}^D}_\lambda(\vectx+\vectx^{-1}),
\end{equation}
and, for a strict partition $\lambda$ of length $n$, we have
\begin{equation}
\label{eq:N=HL-D2}
\PP^{\Phi(D_n)}_\lambda (\vectx;-1)
 +
\PP^{\Phi(D_n)}_{\lambda'} (\vectx;-1)
 =
P^{\mathcal{F}^D}_\lambda(\vectx+\vectx^{-1}),
\end{equation}
where $\lambda' = \lambda_1 \ep_1 + \dots + \lambda_{n-1} \ep_{n-1} - \lambda_n \ep_n$.
\end{enumerate}
\end{theorem}

The first few terms of the sequences $\mathcal{F}^B$, $\mathcal{F}^C$ and $\mathcal{F}^D$ are
\begin{alignat*}{3}
f^B_1(u) &= u + 2, 
&\quad
f^B_2(u) &= u^2 + 2 u, 
&\quad
f^B_3(u) &= u^3 + 2 u^2 - u -2,
\\
f^C_1(u) &= u,
&\quad
f^C_2(u) &= u^2, 
&\quad
f^C_3(u) &= u^3 - u,
\\
f^D_1(u) &= u, 
&\quad
f^D_2(u) &= u^2 - 2, 
&\quad
f^D_3(u) &= u^3 - 3 u. 
\end{alignat*}
See Lemma~\ref{lem:GF-BCD0} for the generating functions of $f^X_d(x+x^{-1})$.

\begin{definition}
\label{def:P-BCD}
Let $\vectx = (x_1, \dots, x_n)$ be indeterminates and $\lambda$ a strict partition of length $l \le n$.
We define symmetric Laurent polynomials $P^X_\lambda(\vectx)$ and 
$Q^X_\lambda(\vectx)$ of type $X_n$, where $X \in \{ B, C, D \}$, by putting
$$
P^X_\lambda(\vectx) = P^{\FF^X}_\lambda(\vectx+\vectx^{-1}),
\quad
Q^X_\lambda(\vectx) = 2^l P^{\FF^X}_\lambda(\vectx+\vectx^{-1}),
$$
where $\FF^X$ is the polynomial sequence given in Theorem~\ref{thm:N=HL-BCD}.
We call $Q^B_\lambda(\vectx)$, $Q^C_\lambda(\vectx)$ and $Q^D_\lambda(\vectx)$ 
the \emph{odd orthogonal $Q$-function}, 
\emph{symplectic $Q$-function} and \emph{even orthogonal $Q$-function} respectively.
\end{definition}

Note that $Q^X_\lambda(\vectx)$ is obtained as the generalized $P$-function $P^\GG_\lambda(\vectx+\vectx^{-1})$ 
associated to the sequence $\GG^X = \{ g^X_d \}_{d=0}^\infty$ given by
$$
g^X_d(u)
 =
\begin{cases}
 1 &\text{if $d=0$,} \\
 2 f^X(u) &\text{if $d \ge 1$.}
\end{cases}
$$

In order to prove Theorem~\ref{thm:N=HL-BCD}, we recall the structure of the Weyl group 
of type $B_n$, $C_n$ and $D_n$.
Let $T_n$ be the abelian group of order $2^n$ generated by $t_1, \dots, t_n$ subject to the relations 
$t_i^2 = 1$ ($1 \le i \le n$) and $t_i t_j = t_j t_i$ ($1 \le i, j \le n$), 
and $W_n = T_n \rtimes S_n$ the semidirect product of $T_n$ with the symmetric group $S_n$, 
where $S_n$ acts on $T_n$ by permuting $t_1, \dots, t_n$.
Put $T'_n = \{ t_1^{u_1} \dots t_n^{u_n} : \sum_{i=1}^n u_i \equiv 0 \bmod 2 \}$ 
and $W'_n = T'_n \rtimes S_n$ the semidirect product of $T'_n$ with $S_n$.
Then $W_n$ is isomorphic to the Weyl group of type $B_n$ and $C_n$, 
and $W'_n$ is isomorphic to the Weyl group of type $D_n$.
The natural action of $S_n$ on $V$ and $K[x_1^{\pm 1}, \dots, x_n^{\pm 1}]$ 
is extended to $W_n$ by
$$
t_i \ep_k
 =
\begin{cases}
 - \ep_i &\text{if $k=i$}, \\
 \ep_k &\text{if $k \neq i$,}
\end{cases}
\quad
t_i x_k
 = 
\begin{cases}
 x_i^{-1} &\text{if $k = i$,} \\
 x_k &\text{if $k \neq i$.}
\end{cases}
$$
If $\lambda$ is a strict partition of length $l$, 
then the stabilizer $W_\lambda$ of $\lambda_1 \ep_1 + \dots + \lambda_l \ep_l$ 
is isomorphic to $W_{n-l}$ for types $B_n$ and $C_n$, 
and to $W'_{n-l}$ for type $D_n$.

\begin{lemma}
\label{lem:HL-BCD}
For a strict partition $\lambda$ of length $l$, we have
\begin{align*}
\PP^{\Phi(B_n)}_{[\lambda]}(\vectx;t)
 &=
\sum_{w \in W_n/W_{n-l}}
 w \left(
  \prod_{i=1}^l x_i^{\lambda_i}
  \prod_{i=1}^l
   \frac{ 1 - t x_i^{-1} }{ 1 - x_i^{-1} }
  \prod_{\substack{1 \le i < j \le n \\ i \le l}}
   \frac{ (1 - t x_i^{-1} x_j) (1 - t x_i^{-1} x_j^{-1}) }
        { (1 - x_i^{-1} x_j) (1 - x_i^{-1} x_j^{-1}) }
 \right),
\\
\PP^{\Phi(C_n)}_{[\lambda]}(\vectx;t)
 &=
\sum_{w \in W_n/W_{n-l}}
 w \left(
  \prod_{i=1}^l x_i^{\lambda_i}
  \prod_{i=1}^l
   \frac{ 1 - t x_i^{-2} }{ 1 - x_i^{-2} }
  \prod_{\substack{1 \le i < j \le n \\ i \le l}}
   \frac{ (1 - t x_i^{-1} x_j) (1 - t x_i^{-1} x_j^{-1}) }
        { (1 - x_i^{-1} x_j) (1 - x_i^{-1} x_j^{-1}) }
 \right),
\\
\PP^{\Phi(D_n)}_{[\lambda]}(\vectx;t)
 &=
\sum_{w \in W'_n/W'_{n-l}}
 w \left(
  \prod_{i=1}^l x_i^{\lambda_i}
  \prod_{\substack{1 \le i < j \le n \\ i \le l}}
   \frac{ (1 - t x_i^{-1} x_j) (1 - t x_i^{-1} x_j^{-1}) }
        { (1 - x_i^{-1} x_j) (1 - x_i^{-1} x_j^{-1}) }
 \right).
\end{align*}
\end{lemma}

\begin{demo}{Proof}
For a general root system $\Phi$ with Weyl group $W$, 
we have (see \cite[Theorem~2.8]{Macdonald72})
$$
\sum_{w \in W}
 w \left(
  \prod_{\alpha \in \Phi^+}
   \frac{ 1 - t e^{-\alpha} }
        { 1 - e^{-\alpha} }
 \right)
 =
\sum_{w \in W}
 t^{l(w)}.
$$
We can use the same argument as in the proof of Lemma~\ref{lem:HL} to prove this lemma.
\end{demo}

\begin{demo}{Proof of Theorem~\ref{thm:N=HL-BCD}}
(1)
We can take $\{ w t : w \in S_n/S_{n-l}, t \in T_l \}$ as 
a complete set of coset representatives of $W_n/W_{n-l}$, 
where $T_l = \langle t_1, \dots, t_l \rangle$.
Since the product
$$
\prod_{\substack{1 \le i < j \le n \\ i \le l}}
\frac{ (1 + x_i^{-1} x_j) (1 + x_i^{-1} x_j^{-1}) }
     { (1 - x_i^{-1} x_j) (1 - x_i^{-1} x_j^{-1}) }
$$
is invariant under $T_l$, we see that
\begin{multline*}
\sum_{w \in T_l}
 w \left(
  \prod_{i=1}^l
   x_i^{\lambda_i}
   \frac{ 1 + x_i^{-1} }{ 1 - x_i^{-1} }
  \prod_{\substack{1 \le i < j \le n \\ i \le l}}
   \frac{ (1 + x_i^{-1} x_j) (1 + x_i^{-1} x_j^{-1}) }
        { (1 - x_i^{-1} x_j) (1 - x_i^{-1} x_j^{-1}) }
 \right)
\\
 =
  \prod_{i=1}^l
   \left(
    x_i^{\lambda_i} \frac{ 1 - s x_i^{-1} }{ 1 - x_i^{-1} }
    +
    x_i^{-\lambda_i} \frac{ 1 - s x_i^{1} }{ 1 - x_i^{1} }
   \right)
  \prod_{\substack{1 \le i < j \le n \\ i \le l}}
   \frac{ (1 + x_i^{-1} x_j) (1 + x_i^{-1} x_j^{-1}) }
        { (1 - x_i^{-1} x_j) (1 - x_i^{-1} x_j^{-1}) }.
\end{multline*}
By using
$$
\frac{ (1 + x_i^{-1} x_j) (1 + x_i^{-1} x_j^{-1}) }
     { (1 - x_i^{-1} x_j) (1 - x_i^{-1} x_j^{-1}) }
 =
\frac{ (x_j + x_j^{-1}) + (x_i + x_i^{-1}) }
     { (x_j + x_j^{-1}) - (x_i + x_i^{-1}) },
$$
we have
\begin{multline*}
\PP^{\Phi(B_n)}_\lambda(\vectx;-1)
\\
 =
\sum_{w \in S_n/S_{n-l}}
 w \left(
  \prod_{i=1}^l
   \left( x_i^{\lambda_i} - x_i^{-\lambda_i} \right)
   \frac{ x_i^{1/2} + x_i^{-1/2} }
        { x_i^{1/2} - x_i^{-1/2} }
  \prod_{\substack{1 \le i < j \le n \\ i \le l}}
   \frac{ (x_j + x_j^{-1}) + (x_i + x_i^{-1}) }
        { (x_j + x_j^{-1}) - (x_i + x_i^{-1}) }
 \right).
\end{multline*}
Comparing this with (\ref{eq:HL3}), we obtain (\ref{eq:N=HL-B}).

(2) can be shown in the same manner as (1).

(3)
Suppose that $l < n$.
Since $W_n = W'_n \sqcup W'_n t_n$, $W_{n-l} = W'_{n-l} \sqcup W'_{n-l} t_n$ and 
$$
\prod_{\substack{1 \le i < j \le n \\ i \le l}}
\frac{ (1 + x_i^{-1} x_j) (1 + x_i^{-1} x_j^{-1}) }
     { (1 - x_i^{-1} x_j) (1 - x_i^{-1} x_j^{-1}) }
$$
is invariant under $t_n$, we see that
$$
\PP^{\Phi(D_n)}_\lambda(\vectx;-1)
 =
\sum_{w \in W_n/W_{n-l}}
 w \left(
  \prod_{i=1}^l x_i^{\lambda_i}
  \prod_{\substack{1 \le i < j \le n \\ i \le l}}
   \frac{ (x_j + x_j^{-1}) + (x_i + x_i^{-1}) }
        { (x_j + x_j^{-1}) - (x_i + x_i^{-1}) }
 \right).
$$
If $l = n$, then the stabilizer $W_\lambda$ is trivial and we have
\begin{align*}
&
\PP^{\Phi(D_n)}_\lambda(\vectx;-1) + \PP^{\Phi(D_n)}_{\lambda'}(\vectx;-1)
\\
 &\quad
=
\sum_{w \in W'_n}
 w \left(
  \prod_{i=1}^{n-1} x_i^{\lambda_i} x_n^{\lambda_n}
  \prod_{1 \le i < j \le n} 
   \frac{ (1 + x_i^{-1} x_j) (1 + x_i^{-1} x_j^{-1}) }
        { (1 - x_i^{-1} x_j) (1 - x_i^{-1} x_j^{-1}) }
 \right)
\\
&\quad\quad
+
\sum_{w \in W'_n}
 w \left(
  \prod_{i=1}^{n-1} x_i^{\lambda_i} x_n^{-\lambda_n}
  \prod_{1 \le i < j \le n} 
   \frac{ (1 + x_i^{-1} x_j) (1 + x_i^{-1} x_j^{-1}) }
        { (1 - x_i^{-1} x_j) (1 - x_i^{-1} x_j^{-1}) }
 \right)
\\
 &\quad
=
\sum_{w \in W_n}
 w \left(
  \prod_{i=1}^{n} x_i^{\lambda_i}
  \prod_{1 \le i < j \le n} 
   \frac{ (x_j + x_j^{-1}) + (x_i + x_i^{-1}) }
        { (x_j + x_j^{-1}) - (x_i + x_i^{-1}) }
 \right).
\end{align*}
The rest of the proof is the same as (1).
\end{demo}

\subsection{%
Generating functions
}

Since the $Q$-functions $Q^X_\lambda(\vectx)$ of type $X$ are special cases of 
generalized $P$-functions, we have the Schur-type Pfaffian formula:

\begin{prop}
\label{prop:Schur-BCD}
For a strict partition $\lambda$ of length $l$, we have
$$
Q^X_\lambda(\vectx)
 =
\Pf \left(
 Q^X_{(\lambda_i,\lambda_j)}(\vectx)
\right)_{1 \le i, j \le r},
$$
where $r = l$ or $l+1$ according whether $l$ is even or odd, 
and we use the convention (\ref{eq:convention2}).
\end{prop}

Hence, in order to obtain $Q^X_\lambda(\vectx)$ for a general strict partition $\lambda$, 
we need to know $Q^X_{(r)}(\vectx)$ and $Q^X_{(r,s)}(\vectx)$.
We compute the generating functions for them.
To state formulas, we introduce formal power series $\varphi^X(z)$ and $\psi^X(z)$ 
by putting
$$
\varphi^X(z)
 =
\begin{cases}
 \dfrac{ (1+z)^2 }{ 1 + z^2 } &\text{if $X=B$,} \\
 1                            &\text{if $X=C$,} \\
 \dfrac{ 1 - z^2 }{ 1 + z^2 } &\text{if $X=D$,}
\end{cases}
\quad
\psi^X(z)
 =
\begin{cases}
 \dfrac{ 2z }{ 1 + z^2 } &\text{if $X=B$,} \\
 0                       &\text{if $X=C$,} \\
 - \dfrac{ 2 z^2}{1 + z^2} &\text{if $X=D$.}
\end{cases}
$$
Note that $\varphi^X(z) - \psi^X(z) = 1$.
And we put
$$
\tilde{\Pi}_z(\vectx)
 =
\prod_{i=1}^n
 \frac{ (1 + x_i z) (1 + x_i^{-1} z) }
      { (1 - x_i z) (1 - x_i^{-1} z) }.
$$
Then we have

\begin{prop}
\label{prop:GF-BCD12}
\begin{enumerate}
\item[(1)]
The generating function of $Q^X_{(r)}(\vectx)$ is given by
$$
\sum_{r=0}^\infty Q^X_{(r)}(\vectx) z^r
 =
\varphi^X(z) \tilde{\Pi}_z(\vectx) - (-1)^n \psi^X(z).
$$
\item[(2)]
The generating function of $Q^X_{(r,s)}(\vectx)$ is given by
\begin{align*}
\sum_{r,s \ge 0} Q^X_{(r,s)}(\vectx) z^r w^s
&=
\frac{ (z-w)(1-zw) }
     { (z+w)(1+zw) }
\cdot
\varphi^X(z) \varphi^X(w)
\left(
 \tilde{\Pi}_z(\vectx) \tilde{\Pi}_w(\vectx) - 1 
\right)
\\
&\quad
+ (-1)^n
\left(
 \varphi^X(z) \psi^X(w) \tilde{\Pi}_z(\vectx)
 - \varphi^X(w) \psi^X(z) \tilde{\Pi}_w(\vectx)
\right)
\\
&\quad
+ \psi^X(z) - \psi^X(w).
\end{align*}
\end{enumerate}
\end{prop}

By a straightforward case-by-case computation, we can show the following lemma:

\begin{lemma}
\label{lem:GF-BCD0}
We have
$$
1 + 2 \sum_{r=1}^\infty f^X_r(x+x^{-1}) z^r
 =
\varphi^X(z)
 \cdot
\frac{ (1 + xz) (1 + x^{-1}z) }
     { (1 - xz) (1 - x^{-1}z) }
 +
\psi^X(z).
$$
\end{lemma}

\begin{demo}{Proof of Proposition~\ref{prop:GF-BCD12}}
By Theorem~\ref{thm:N=HL-BCD} and (\ref{eq:Nimmo}), we see that 
the $Q$-functions $Q^X_\lambda(\vectx)$ corresponding to a strict partition $\lambda$ of length $l$ 
is expressed as
$$
Q^X_\lambda(\vectx)
 =
\frac{ 1 }
     { \tilde{\Delta}(\vectx) }
\Pf \begin{pmatrix}
 \tilde{A}(\vectx) & \tilde{V}^X_{\lambda^*}(\vectx) \\
 -\trans \tilde{V}^X_{\lambda^*}(\vectx) & O
\end{pmatrix},
$$
where $\lambda^* = \lambda$ or $\lambda^0$ according whether $n+l$ is even or odd, 
and
$$
\tilde{V}^X_{(\alpha_1, \dots, \alpha_r)}(\vectx)
 =
\left(
 \chi(\alpha_j) f^X_{\alpha_j}(x_i+x_i^{-1})
\right)_{1 \le i \le n, 1 \le j \le r},
\quad
\chi(d)
 = 
\begin{cases}
 1 &\text{if $d=0$,} \\
 2 &\text{if $d \ge 1$.}
\end{cases}
$$
Now by an argument similar to the proof of Proposition~\ref{prop:GF}, 
we can prove this proposition by using the Pfaffian evaluations in Proposition~\ref{prop:Schur-Pf-BCD}.
The details are left to the readers.
\end{demo}

\appendix
\section{%
Pfaffian formulas
}

In this appendix, we collect several useful Pfaffian identities.

\subsection{%
Pfaffians
}

Recall the definition and some properties of Pfaffians 
(see \cite{IO} for some expositions).
Let $X = \bigl( x_{ij} \bigr)_{1 \le i, j \le 2m}$ be a skew-symmetric matrix of order $2m$.
The \emph{Pfaffian} of $X$, denoted by $\Pf (X)$, is defined by
\begin{equation}
\label{eq:def-Pf}
\Pf (X)
 =
\sum_{\sigma \in F_{2m}} \sgn(\sigma) \prod_{i=1}^m x_{\sigma(2i-1), \sigma(2i)},
\end{equation}
where $F_{2m}$ is the set of permutations $\sigma \in S_{2m}$ satisfying 
$\sigma(1) < \sigma(3) < \dots < \sigma(2m-1)$ and $\sigma(2i-1) < \sigma(2i)$ for $1 \le i \le m$.

Pfaffians are multilinear and alternating in the following sense.
Let $X = \bigl( x_{ij} \bigr)_{1 \le i, j \le n}$ be a skew-symmetric matrix 
and fix a row/column index $k$.
If the entries of the $k$th row and $k$th column of $X$ are written as 
$x_{i,j} = \alpha x'_{i,j} + \beta x''_{i,j}$ for $i=k$ or $j=k$, 
then
$$
\Pf X = \alpha \Pf X' + \beta \Pf X'',
$$
where $X'$ (resp. $X''$) is the skew-symmetric matrix obtained from $X$ by 
replacing the entries $x_{ij}$ for $i=k$ or $j=k$ with $x'_{ij}$ (resp. $x''_{ij}$).
And, for a permutation $\sigma \in S_n$, we have
$$
\Pf \bigl( x_{\sigma(i),\sigma(j)} \bigr)_{1 \le i, j \le n}
 =
\sgn(\sigma) \Pf \bigl( x_{i,j} \bigr)_{1 \le i, j \le n}.
$$
It follows that, if $Y$ is the skew-symmetric matrix obtained from $X$ 
by adding the $k$th row multiplied by a scalar $\alpha$ to the $l$th row and then 
adding the $k$th column multiplied by $\alpha$ to the $l$th column,
the we have $\Pf Y = \Pf X$.

We use the following notations for submatrices.
For a positive integer $n$, we put $[n] = \{ 1, 2, \dots, n \}$.
Given a subset $I \subset [n]$, we put $\Sigma(I) = \sum_{i \in I} i$.
For an $M \times N$ matrix $X = \bigl( x_{i,j} \bigr)_{1 \le i \le M, 1 \le j \le N}$ 
and subsets $I \subset [M]$ and $J \subset [N]$, 
we denote by $X(I;J)$ the submatrix of $X$ obtained by picking up rows indexed by $I$ and 
columns indexed by $J$.
If $X$ is a skew-symmetric matrix, then we write $X(I)$ for $X(I;I)$.
We use the convention that $\det X(\emptyset;\emptyset) = 1$ and $\Pf X(\emptyset) = 1$.

For an $n \times n$ skew-symmetric matrix $X = \bigl( x_{i,j} \bigr)_{1 \le i, j \le n}$, 
we have the following expansion formula along the $k$th row/column:
\begin{equation}
\label{eq:Pf-expansion}
\Pf X
 = 
\sum_{i=1}^{k-1} (-1)^{k+i-1} x_{i,k} \Pf X([n] \setminus \{ i, k \})
 +
\sum_{i=k+1}^n (-1)^{k+i-1} x_{k,i} \Pf X([n] \setminus \{ k, i \}).
\end{equation}

\subsection{%
Schur's Pfaffian evaluation and its variations
}

Recall that
$$
A(\vectx) = \left( \frac{ x_j - x_i }{ x_j + x_i } \right)_{1 \le i, j \le n},
\quad
\Delta(\vectx) = \prod_{1 \le i < j \le n} \frac{ x_j - x_i }{ x_j + x_i }
$$
for a sequence $\vectx = (x_1, \dots, x_n)$ of indeterminates.
The following evaluation of the Pfaffian (\ref{eq:Schur-Pf1}) originates form \cite{Schur}, 
and its simple proof can be found in \cite{Knuth}.
Equation (\ref{eq:Schur-Pf2}) is derived from (\ref{eq:Schur-Pf1}) 
by specializing the last indeterminate to $0$.

\begin{prop}
\label{prop:Schur-Pf}
If $n$ is even, then we have
\begin{equation}
\label{eq:Schur-Pf1}
\Pf A(\vectx) = \Delta(\vectx).
\end{equation}
If $n$ is odd, then we have
\begin{equation}
\label{eq:Schur-Pf2}
\Pf \begin{pmatrix}
A(\vectx) & \vectone \\
-\vectone & 0
\end{pmatrix}
 =
\Delta(\vectx),
\end{equation}
where $\vectone$ is the all-one column vector.
\end{prop}

For two sequences $\vectx = (x_1, \dots, x_n)$ and $\vecty = (y_1, \dots, y_p)$ of indeterminates, 
we put
$$
B(\vectx;\vecty) = \left( \frac{ 1 + x_i y_j }{ 1 - x_i y_j } \right)_{1 \le i \le n, 1 \le j \le p},
\quad
\Pi(\vectx;\vecty) = \prod_{i=1}^n \prod_{j=1}^l \frac{ 1 + x_i y_j }{ 1 - x_i y_j }.
$$
Let $B_z(\vectx) = B(\vectx;(z))$ be the column vector with $i$th entry $(1+x_iz)/(1-x_iz)$ and 
set $\Pi_z(\vectx) = \Pi(\vectx;(z)) = \prod_{i=1}^n (1 + x_i z)/(1 - x_i z)$.
Then we have the following variations of Schur's Pfaffian evaluation.

\begin{prop}
\label{prop:Schur-Pf-var}
\begin{enumerate}
\item[(1)]
If $n+p$ is even, then we have
\begin{equation}
\label{eq:Schur-Pf3}
\Pf \begin{pmatrix}
 A(\vectx) & B(\vectx;\vecty) \\
 -\trans B(\vectx;\vecty) & - A(\vecty)
\end{pmatrix}
 =
(-1)^{\binom{p}{2}}
\Delta(\vectx) \Delta(\vecty) \Pi(\vectx;\vecty).
\end{equation}
\item[(2)]
If $n$ is even, then we have
\begin{equation}
\label{eq:Schur-Pf4}
\Pf \begin{pmatrix}
 A(\vectx) & B_z(\vectx) & B_w(\vectx) \\
 -\trans B_z(\vectx) & 0 & 0 \\
 -\trans B_w(\vectx) & 0 & 0
\end{pmatrix}
 =
\Delta(\vectx) \cdot \frac{ z-w }{ z+w }
\big( \Pi_z(\vectx) \Pi_w(\vectx) - 1 \big).
\end{equation}
\item[(3)]
If $n$ is odd, then we have
\begin{multline}
\label{eq:Schur-Pf5}
\Pf \begin{pmatrix}
 A(\vectx) & B_z(\vectx) & B_w(\vectx) & \vectone \\
 -\trans B_z(\vectx) & 0 & 0 & 0 \\
 -\trans B_w(\vectx) & 0 & 0 & 0 \\
 -\trans\vectone & 0 & 0 & 0
\end{pmatrix}
\\
 =
\Delta(\vectx)
\cdot
\left\{
 \frac{ z-w }{ z+w }
 \big( \Pi_z(\vectx) \Pi_w(\vectx) -1 \big)
 - \Pi_z(\vectx) + \Pi_w(\vectx)
\right\}.
\end{multline}
\end{enumerate}
\end{prop}

\begin{demo}{Proof}
(1)
Apply (\ref{eq:Schur-Pf1}) to the indeterminates 
$(x_1, \dots, x_n, -1/y_1, \dots, \allowbreak -1/y_p)$.

(2)
By applying (1) with $p=2$ and $(y_1, y_2) = (z,w)$, we have
$$
\Pf \begin{pmatrix}
 A(\vectx) & B_z(\vectx) & B_w(\vectx) \\
 -\trans B_z(\vectx) & 0 & \dfrac{z-w}{z+w} \\
 -\trans B_w(\vectx) & -\dfrac{z-w}{z+w} & 0
\end{pmatrix}
 =
\Delta(\vectx) \Pi_z(\vectx) \Pi_w(\vectx) \cdot \frac{z-w}{z+w}.
$$
By using the multilinearity of Pfaffians, we obtain
\begin{multline*}
\Pf \begin{pmatrix}
 A(\vectx) & B_z(\vectx) & B_w(\vectx) \\
 -\trans B_z(\vectx) & 0 & \dfrac{z-w}{z+w} \\
 -\trans B_w(\vectx) & -\dfrac{z-w}{z+w} & 0
\end{pmatrix}
\\
 =
\Pf \begin{pmatrix}
 A(\vectx) & B_z(\vectx) & B_w(\vectx) \\
 -\trans B_z(\vectx) & 0 & 0 \\
 -\trans B_w(\vectx) & 0 & 0
\end{pmatrix}
 +
\frac{z-w}{z+w}
\cdot
\Pf \begin{pmatrix}
 A(\vectx) & B_z(\vectx) & 0 \\
 -\trans B_z(\vectx) & 0 & 1 \\
 0 & -1 & 0
\end{pmatrix}.
\end{multline*}
The last Pfaffian is shown to be equal to $\Delta(\vectx)$ 
by expanding along the last row/column and using (\ref{eq:Schur-Pf1}), 
and thus we obtain (\ref{eq:Schur-Pf4}).

(3)
Applying (1) with $p=3$ and $(y_1, y_2, y_3) = (z,w,0)$, we obtain
$$
\Pf \begin{pmatrix}
 A(\vectx) & B_z(\vectx) & B_w(\vectx) & -\vectone \\
 -\trans B_z(\vectx) & 0 & \dfrac{z-w}{z+w} & -1 \\
 -\trans B_w(\vectx) & -\dfrac{z-w}{z+w} & 0 & -1 \\
 \trans\vectone & 1 & 1 & 0
\end{pmatrix}
 =
\Delta(\vectx) \Pi_z(\vectx) \Pi_w(\vectx) \cdot \frac{ z-w }{ z+w } \cdot (-1)^{n+2}.
$$
By using the multilinearity we see that
\begin{align*}
&
\Pf \begin{pmatrix}
 A(\vectx) & B_z(\vectx) & B_w(\vectx) & -\vectone \\
 -\trans B_z(\vectx) & 0 & \dfrac{z-w}{z+w} & -1 \\
 -\trans B_w(\vectx) & -\dfrac{z-w}{z+w} & 0 & -1 \\
 \trans\vectone & 1 & 1 & 0
\end{pmatrix}
\\
&\quad
=
-
\Pf \begin{pmatrix}
 A(\vectx) & B_z(\vectx) & B_w(\vectx) & \vectone \\
 -\trans B_z(\vectx) & 0 & 0 & 0 \\
 -\trans B_w(\vectx) & 0 & 0 & 0 \\
 -\trans\vectone & 0 & 0 & 0
\end{pmatrix}
-
\frac{z-w}{z+w}
\Pf \begin{pmatrix}
 A(\vectx) & B_z(\vectx) & 0 & \vectone \\
 -\trans B_z(\vectx) & 0 & 1 & 0 \\
 -\trans B_w(\vectx) & -1 & 0 & 0 \\
 -\trans\vectone & 0 & 0 & 0
\end{pmatrix}
\\
&\quad\quad
-
\Pf \begin{pmatrix}
 A(\vectx) & B_z(\vectx) & B_w(\vectx) & 0 \\
 -\trans B_z(\vectx) & 0 & \dfrac{z-w}{z+w} & 1 \\
 -\trans B_w(\vectx) & -\dfrac{z-w}{z+w} & 0 & 0 \\
 0 & -1 & 0 & 0
\end{pmatrix}
-
\Pf \begin{pmatrix}
 A(\vectx) & B_z(\vectx) & B_w(\vectx) & 0 \\
 -\trans B_z(\vectx) & 0 & \dfrac{z-w}{z+w} & 0 \\
 -\trans B_w(\vectx) & -\dfrac{z-w}{z+w} & 0 & 1 \\
 0 & 0 & -1 & 0
\end{pmatrix}.
\end{align*}
The last three Pfaffians can be evaluated by expanding them along a row/column and then 
by using (\ref{eq:Schur-Pf2}) and (\ref{eq:Schur-Pf3}) with $p=1$ as follows: 
\begin{align*}
\Pf \begin{pmatrix}
 A(\vectx) & B_z(\vectx) & 0 & \vectone \\
 -\trans B_z(\vectx) & 0 & 1 & 0 \\
 -\trans B_w(\vectx) & -1 & 0 & 0 \\
 -\trans\vectone & 0 & 0 & 0
\end{pmatrix}
&=
\Pf \begin{pmatrix}
 A(\vectx) & \vectone \\
 -\trans\vectone & 0
\end{pmatrix}
 =
\Delta(\vectx),
\\
\Pf \begin{pmatrix}
 A(\vectx) & B_z(\vectx) & B_w(\vectx) & 0 \\
 -\trans B_z(\vectx) & 0 & \dfrac{z-w}{z+w} & 1 \\
 -\trans B_w(\vectx) & -\dfrac{z-w}{z+w} & 0 & 0 \\
 0 & -1 & 0 & 0
\end{pmatrix}
 &=
-
\Pf \begin{pmatrix}
 A(\vectx) & B_w(\vectx) \\
 -\trans B_w(\vectx) & 0
\end{pmatrix}
 =
- \Delta(\vectx) \Pi_w(\vectx),
\\
\Pf \begin{pmatrix}
 A(\vectx) & B_z(\vectx) & B_w(\vectx) & 0 \\
 -\trans B_z(\vectx) & 0 & \dfrac{z-w}{z+w} & 0 \\
 -\trans B_w(\vectx) & -\dfrac{z-w}{z+w} & 0 & 1 \\
 0 & 0 & -1 & 0
\end{pmatrix}
 &=
\Pf \begin{pmatrix}
 A(\vectx) & B_z(\vectx) \\
 -\trans B_z(\vectx) & 0
\end{pmatrix}
 =
\Delta(\vectx) \Pi_z(\vectx).
\end{align*}
Combining these evaluations completes the proof of (\ref{eq:Schur-Pf5}).
\end{demo}

The following Pfaffian evaluations are used in Section~7.
For $\vectx = (x_1, \dots, x_n)$ and $\vecty = (y_1, \dots, y_p)$, we put
\begin{gather*}
\tilde{A}(\vectx)
 =
\left(
 \frac{ (x_j + x_j^{-1}) - (x_i + x_i^{-1}) }
      { (x_j + x_j^{-1}) + (x_i + x_i^{-1}) }
\right)_{1 \le i, j \le n}
 =
\left(
 \frac{ (x_j - x_i) (1 - x_i x_j) }
      { (x_j + x_i) (1 + x_i x_j) }
\right)_{1 \le i, j \le n},
\\
\tilde{\Delta}(\vectx)
 =
\prod_{1 \le i < j \le n}
 \frac{ (x_j + x_j^{-1}) - (x_i + x_i^{-1}) }
      { (x_j + x_j^{-1}) + (x_i + x_i^{-1}) }
 =
\prod_{1 \le i < j \le n}
 \frac{ (x_j - x_i) (1 - x_i x_j) }
      { (x_j + x_i) (1 + x_i x_j) },
\\
\tilde{B}(\vectx;\vecty)
 =
\left(
 \frac{ (1 - x_i y_j) (1 - x_i^{-1} y_j) }
      { (1 + x_i y_j) (1 + x_i^{-1} y_j) }
\right)_{1 \le i \le n, 1 \le j \le p},
\\
\tilde{\Pi}(\vectx;\vecty)
 =
\prod_{i=1}^n \prod_{j=1}^p
 \frac{ (1 - x_i y_j) (1 - x_i^{-1} y_j) }
      { (1 + x_i y_j) (1 + x_i^{-1} y_j) }.
\end{gather*}
We write $\tilde{B}_z(\vectx) = \tilde{B}(\vectx;(z))$ 
and $\tilde{\Pi}_z(\vectx) = \tilde{\Pi}(\vectx;(z))$.
Then we have

\begin{prop}
\label{prop:Schur-Pf-BCD}
\begin{enumerate}
\item[(1)]
If $n$ is even, then we have
\begin{equation}
\label{eq:Schur-Pf-BCD1}
\Pf \tilde{A}(\vectx) = \tilde{\Delta}(\vectx).
\end{equation}
\item[(2)]
If $n$ is odd, then we have
\begin{equation}
\label{eq:Schur-Pf-BCD2}
\Pf \begin{pmatrix}
 \tilde{A}(\vectx) & \vectone \\
 -\trans\vectone & 0
\end{pmatrix}
 = 
\tilde{\Delta}(\vectx).
\end{equation}
\item[(3)]
If $n+p$ is even, then we have
\begin{equation}
\label{eq:Schur-Pf-BCD3}
\Pf \begin{pmatrix}
 \tilde{A}(\vectx) & \tilde{B}(\vectx;\vecty) \\
 -\trans \tilde{B}(\vectx;\vecty) & - \tilde{A}(\vecty)
\end{pmatrix}
 =
\tilde{\Delta}(\vectx) \tilde{\Delta}(\vecty) \tilde{\Pi}(\vectx;\vecty).
\end{equation}
\item[(4)]
If $n$ is even, then we have
\begin{equation}
\label{eq:Schur-Pf-BCD4}
\Pf \begin{pmatrix}
 \tilde{A}(\vectx) & \tilde{B}_z(\vectx) & \tilde{B}_w(\vectx) \\
 -\trans \tilde{B}_z(\vectx) & 0 & 0 \\
- \trans \tilde{B}_w(\vectx) & 0 & 0
\end{pmatrix}
 =
\tilde{\Delta}(\vectx)
\frac{ (z - w) (1 - z w) }
     { (z + w) (1 + z w) }
\big( \tilde{\Pi}_z(\vectx) \tilde{\Pi}_w(\vectx) - 1 \big).
\end{equation}
\item[(5)]
If $n$ is odd, then we have
\begin{multline}
\label{eq:Schur-Pf-BCD5}
\Pf \begin{pmatrix}
 \tilde{A}(\vectx) & \tilde{B}_z(\vectx) & \tilde{B}_w(\vectx) & \vectone \\
 -\trans \tilde{B}_z(\vectx) & 0 & 0 & 0 \\
 -\trans \tilde{B}_w(\vectx) & 0 & 0 & 0 \\
 -\trans\vectone & 0 & 0 & 0
\end{pmatrix}
\\
 =
\tilde{\Delta}(\vectx)
\left\{
 \frac{ (z - w) (1 - z w) }
      { (z + w) (1 + z w) }
 \big( 
  \tilde{\Pi}_z(\vectx) \tilde{\Pi}_w(\vectx) - 1
 \big)
 - \tilde{\Pi}_z(\vectx) + \tilde{\Pi}_w(\vectx)
\right\}.
\end{multline}
\end{enumerate}
\end{prop}

\begin{demo}{Proof}
(1) and (2) are obtained by replacing $x_i$ with $x_i+x_i^{-1}$ 
in (\ref{eq:Schur-Pf1}) and (\ref{eq:Schur-Pf2}) respectively.
(3) is obtained by applying (\ref{eq:Schur-Pf1}) 
with $(x_1+x_1^{-1}, \dots, x_n+x_n^{-1}, -(y_1+y_1^{-1}), \dots, -(y_p+y_p^{-1}))$.
(4) and (5) are derived from (3) by the similar argument to the proof of (\ref{eq:Schur-Pf4}) 
and (\ref{eq:Schur-Pf5}) respectively.
\end{demo}

\subsection{%
Useful formulas for Pfaffians
}

The following propositions are Pfaffian analogues of 
the Sylvester identity, 
the Laplace expansion formula, 
and the Cauchy--Binet formula for determinants.

\begin{prop}
\label{prop:Pf-Sylvester}
(\cite[(2.5)]{Knuth})
Let $n$ and $m$ be even integers.
If $X$ is an $(n+m) \times (n+m)$ skew-symmetric matrix such that $\Pf X([n]) \neq 0$, 
then we have
\begin{equation}
\label{eq:Pf-Sylvester}
\Pf \left(
 \frac{ \Pf X([n] \cup \{ n+i, n+j \}) }
      { \Pf X([n]) }
\right)_{1 \le i, j \le m}
 =
\frac{ \Pf X }
     { \Pf X([n]) }.
\end{equation}
\end{prop}

\begin{prop}
\label{prop:Pf-Laplace}
(\cite[Corollary~2.4 (1)]{Okada})
Let $n$ and $l$ be nonnegative integers with the same parity.
If $Z$ is an $m \times m$ skew-symmetric matrix and $W$ is an $m \times n$ matrix, then we have
\begin{equation}
\label{eq:Pf-Laplace1}
\Pf \begin{pmatrix}
 Z & W \\
 -\trans W & O_{n,n}
\end{pmatrix}
 =
\begin{cases}
\displaystyle\sum_{I}
 (-1)^{\Sigma(I) + \binom{m}{2}}
 \Pf Z(I) \det W( [m] \setminus I ; [n] )
 &\text{if $m > n$,}
\\
(-1)^{\binom{m}{2}} \det W
 &\text{if $m = n$,}
\\
0
 &\text{if $m < n$,}
\end{cases}
\end{equation}
where $I$ runs over all $(m-n)$-element subsets of $[n]$.
\end{prop}

\begin{prop}
(\cite[Theorem~3.2]{Okada})
\label{prop:Pf-CB}
Let $m$ and $n$ be nonnegative integers with the same parity, .
Let $A$ and $B$ be $m \times m$ and $n \times n$ skew-symmetric matrices,
and let $S$ and $T$ be $m \times l$ and $n \times l$ matrices.
Then we have
\begin{multline}
\label{eq:Pf-CB1}
\sum_I (-1)^{ \binom{\# I}{2} }
 \Pf \begin{pmatrix} A & S([m];I) \\ -\trans S([m];I) & O \end{pmatrix} 
 \Pf \begin{pmatrix} B & T([n];I) \\ -\trans T([n];I) & O \end{pmatrix} 
\\
 =
\Pf \begin{pmatrix}
 A & S \trans T \\
 - T \trans S & B
\end{pmatrix},
\end{multline}
\begin{multline}
\label{eq:Pf-CB2}
\sum_I 
 \Pf \begin{pmatrix} A & S([m];I) \\ -\trans S([m];I) & O \end{pmatrix} 
 \Pf \begin{pmatrix} B & T([n];I) \\ -\trans T([n];I) & O \end{pmatrix} 
\\
=
(-1)^{ \binom{n}{2} }
\Pf \begin{pmatrix}
 A & S \trans T \\
 - T \trans S & -B
\end{pmatrix},
\end{multline}
where $I$ runs over all subsets of $[l]$ with $\# I \equiv m \equiv n \bmod 2$.
\end{prop}




\begin{thebibliography}{WW}

\bibitem{CI}
S.~Cho and T.~Ikeda,
Pieri rule for the factorial Schur $P$-functions,
in ``Schubert Varieties, Equivariant Cohomology and Characteristic Classes'', 
eds. J.~Buczy\'nski, M.~Micha{\l}ek, and E.~Postinghel,
EMS Ser. Congr. Rep., Eur. Math. Soc., 2018, pp.~25--48.

\bibitem{IN09}
T.~Ikead and H.~Naruse,
Excited Young diagrams and equivariant Schubert calculus,
Trans. Amer. Math. Soc. {\bfseries 361} (2009), 5193--5221.

\bibitem{IN13}
T.~Ikead and H.~Naruse,
$K$-theoretic analogues of factorial Schur $P$- and $Q$-functions,
Adv. Math. {\bfseries 243} (2013), 22--66.

\bibitem{IO}
M.~Ishikawa and S.~Okada,
Identities for determinants and Pfaffians, and their applications,
Sugaku Expositions {\bf 27} (2014), 85--116.

\bibitem{Ivanov1}
V.~N.~Ivanov,
Combinatorial formula for factorial Schur $Q$-functions,
J. Math. Sci. (N.Y.) {\bfseries 107} (2001), 4195--4211.

\bibitem{Ivanov2}
V.~N.~Ivanov,
Interpolation analogues of Schur $Q$-functions,
J. Math. Sci. (N.Y.) {\bfseries 131} (2005), 5495--5507.

\bibitem{Korotkikh}
S.~Korotkikh,
Dual multiparameter Schur $Q$-functions,
J. Mathematical Sciences {\bfseries 224} (2) (2017), 263--268.
 
\bibitem{Knuth}
D.~E.~Knuth,
Overlapping Pfaffians,
Electron. J. Combin., {\bf 3} (no.~2, The Foata Festschrift) (1996), \#R5.

\bibitem{Macdonald72}
I.~G.~Macdonald,
The Poincar\'e series of a Coxeter group, 
Math. Ann. {\bfseries 199} (1972), 161--174.

\bibitem{Macdonald92}
I.~G.~Macdonald,
Schur functions: Theme and variations,
S\'em. Lothar. Combin. {\bfseries 28} (1992), 5--39.

\bibitem{Macdonald95}
I.~G.~Macdonald,
``Symmetric Functions and Hall Polynomials, 2nd edition'',
Oxford Univ. Press, 1995.

\bibitem{Macdonald00}
I.~G.~Macdonald,
Orthogonal polynomials associated with root systems,
S\'em. Lothar. Combin. {\bfseries 45} (2000/01), Art. B45a.

\bibitem{MS}
A.~I.~Molev and B.~E.~Sagan,
A Littlewood--Richardson rule for factorial Schur functions,
Trans. Amer. Math. Soc. {\bfseries 351} (1999), 4429--4443.

\bibitem{Morris}
A.~O.~Morris,
A note on the multiplication of Hall functions,
J. London Math. Soc. {\bfseries 39} (1964), 481--488.

\bibitem{NNSY}
J.~Nakagawa, M.~Noumi, M.~Shirakawa, and Y.~Yamada
Tableau representation for Macdonald's ninth variation of Schur functions, 
in ``Physics and Combinatorics 2000, Proceedings of 
the Nagoya 2000 International Workshop'' (Eds. A.~N.~ Kirillov and N.~Liskova), 
World Scientific, 2001, 
pp. 180--195,

\bibitem{NN}
M.~Nakagawa and H.~Naruse,
Generalized (co)homology of the loop spaces of classical groups 
and the universal factorial Schur $P$- and $Q$-functions,
in ``Schubert Calculus --- Osaka 2012'' (Eds. H.~Naruse, T.~Ikeda, M.~Masuda and T.~Tanisaki), 
Adv. Stud. Pure Math. {\bfseries 171} (2016), 
pp. 337--417.

\bibitem{Nimmo}
J.~J.~C.~Nimmo,
Hall--Littlewood symmetric functions and the BKP equation,
J. Phys. A {\bfseries 23} (1990), 751--760.

\bibitem{Okada}
S.~Okada,
Pfaffian formulas and Schur $Q$-function identities,
arXiv:1706.01029.

\bibitem{Pragacz}
P.~Pragacz,
Algebro-geometric applications of Schur $S$- and $Q$-polynomials, 
in ``Topics in Invariant Theory'', S\'eminaire d'Alg\`ebre Dubreil-Malliavin 1989-90, 
Lecture Notes in Math. {\bfseries 1478}, 
Springer-Verlag, 1991,
pp. 130--191.

\bibitem{JP}
P.~Pragacz and T.~J\'ozefiak,
A determinantal formula for skew $Q$-functions, 
J. London Math. Soc. (2) {\bf 43} (1991), 76--90.

\bibitem{Schur}
I.~Schur,
\"Uber die Darstellung der symmetrischen und der alternierenden Gruppe durch gebrochene lineare Substitutionen, 
J. Reine Angew. Math. {\bfseries 139} (1911), 155--250.

\bibitem{Stembridge}
J.~R.~Stembridge, 
Nonintersecting paths, pfaffians, and plane partitions,
Adv. Math. {\bfseries 83} (1990), 96--131.

\end{thebibliography}
\end{document}